\documentclass{arxiv}

\usepackage[utf8]{inputenc}
\usepackage[T1]{fontenc}
\usepackage[english]{babel}
\usepackage{enumitem}
\usepackage[hidelinks=true]{hyperref}

\usepackage{xfrac}
\usepackage{cancel}
\usepackage{bbm}

\usepackage{subcaption}
\usepackage{xcolor}

\OneAndAHalfSpacedXI
\usepackage{natbib}
\bibpunct[, ]{(}{)}{,}{a}{}{,}%
\def\bibfont{\small}%
\TheoremsNumberedThrough     % Preferred (Theorem 1, Lemma 1, Theorem 2)
\ECRepeatTheorems

\EquationsNumberedThrough    % Default: (1), (2), ...
\MANUSCRIPTNO{}

\allowdisplaybreaks % To allow page breaks inside equation environment

\newcommand{\ds}{\displaystyle}
\newcommand{\dto}{\downarrow}
\newcommand{\tth}{{}^{\text{th}}}
\newcommand{\defn}{\stackrel{\triangle}{=}}
\renewcommand{\Longrightarrow}{\Rightarrow}

\newcommand{\cA}{\mathcal{A}}
\newcommand{\cC}{\mathcal{C}}
\newcommand{\cK}{\mathcal{K}}

\newcommand{\bR}{\mathbb{R}}

\newcommand{\vy}{\boldsymbol{y}}
\newcommand{\vx}{\boldsymbol{x}}
\newcommand{\ve}{\boldsymbol{e}}
\newcommand{\vone}{\boldsymbol{1}}
\newcommand{\vzero}{\boldsymbol{0}}

\newcommand{\amax}{A_{\max}}
\newcommand{\smax}{S_{\max}}
\newcommand{\plus}{+}

\newcommand{\ind}[1]{\mathbbm{1}_{\left\{#1 \right\}}}
\newcommand{\E}[1]{E\left[#1 \right]}

\newcommand{\Evq}[1]{E_{\vq}\left[#1 \right]}
\newcommand{\Var}[1]{\text{Var}\left[#1 \right]}
\newcommand{\Cov}[1]{\text{Cov}\left[#1 \right]}
\newcommand{\Prob}[1]{P\left[#1\right]}

\def\ba#1\ea{\begin{align*}#1\end{align*}}
\def\ban#1\ean{\begin{align}#1\end{align}}

\newcommand{\qbar}{\overline{q}}
\newcommand{\abar}{\overline{a}}
\newcommand{\sbar}{\overline{s}}
\newcommand{\ubar}{\overline{u}}

\newcommand{\peps}{{(\epsilon)}}
\newcommand{\theteps}{\theta\epsilon}

\newcommand{\vq}{\boldsymbol{q}}
\newcommand{\va}{\boldsymbol{a}}
\renewcommand{\vs}{\boldsymbol{s}}
\newcommand{\vu}{\boldsymbol{u}}
\newcommand{\vqbar}{\overline{\vq}}
\newcommand{\vqpar}{\vq_{\parallel}}
\newcommand{\vqperp}{\vq_{\perp}}
\newcommand{\vqparbar}{\vqbar_{\parallel}}
\newcommand{\vqperpbar}{\vqbar_{\perp}}
\newcommand{\vsbar}{\overline{\vs}}
\newcommand{\vabar}{\overline{\va}}

\newcommand{\pl}{{(\ell)}}

\newcommand{\vlambda}{\boldsymbol{\lambda}}
\newcommand{\vcl}{\boldsymbol{c}^\pl}
\newcommand{\vc}{\boldsymbol{c}}
\newcommand{\bl}{b^{(\ell)}}
\newcommand{\cS}{\mathcal{S}}

\newcommand{\cF}{\mathcal{F}}
\newcommand{\cFl}{\cF^\pl}

\newcommand{\vubar}{\overline{\vu}}

\newcommand{\vr}{\boldsymbol{r}}
\newcommand{\vrl}{\vr^\pl}

\newcommand{\vchi}[1]{\boldsymbol{\chi}^{(#1)}}
\newcommand{\vchitilde}[1]{\widetilde{\boldsymbol{\chi}}^{(#1)}}
\newcommand{\chii}[1]{\chi^{(#1)}}

\newcommand{\pt}{{(t)}}

\newcommand{\cT}{\mathcal{T}}
\newcommand{\Tbar}{\overline{T}}
\newcommand{\Et}[1]{E_t\left[#1 \right]}
\newcommand{\btl}{b^{(t,\ell)}}
\newcommand{\nuq}{\nu_{\vqbar}}
\newcommand{\numin}{\nu_{\min}}

\newcommand{\Bbar}{\overline{B}}

\begin{document}
	
\TITLE{Transform Methods for Heavy-Traffic Analysis}

\ARTICLEAUTHORS{%
	\AUTHOR{Daniela Hurtado-Lange}
	\AFF{Department of Industrial and Systems Engineering, Georgia Institute of Technology\\
		765 Ferst Drive NW, Atlanta, GA 30332,
	\EMAIL{d.hurtado@gatech.edu}} %, \URL{}}
	\AUTHOR{Siva Theja Maguluri}
	\AFF{Department of Industrial and Systems Engineering, Georgia Institute of Technology, \\
		755 Ferst Drive NW, Atlanta, GA 30332
	\EMAIL{siva.theja@gatech.edu}}
} % end of the block

\ABSTRACT{
	The Drift method was recently developed to study queueing systems in steady-state. It was successfully used to obtain bounds on the moments of the scaled queue lengths, that are asymptotically tight in heavy-traffic, in a wide variety of systems including generalized switches \citep{atilla}, input-queued switches \citep{MagSri_SSY16_Switch,QUESTA_switch}, bandwidth sharing networks \citep{Weina_bandwidth_journal}, etc. In this paper we develop the use of transform techniques for heavy-traffic analysis, with a special focus on the use of moment generating functions. This approach simplifies the proofs of the Drift method, and provides a new perspective on the Drift method. We present a general framework and then use the MGF method to obtain the stationary distribution of scaled queue lengths in heavy-traffic in queueing systems that satisfy the Complete Resource Pooling condition. In particular, we study load balancing systems and generalized switches under general settings.
}

\KEYWORDS{Drift method, Heavy-traffic analysis, State Space Collapse, Complete Resource Pooling}

\maketitle

\section{Introduction}\label{sec:introduction}
Exact analysis of queueing systems that arise in the study of Stochastic Processing Networks (SPNs) is usually intractable, so it is common to analyze them in various asymptotic regimes to get insights on their behavior. A very popular regime in the literature is the heavy-traffic regime, where the system is loaded very close to its maximum capacity. This regime is sometimes called the classical or conventional heavy-traffic regime.
One of the advantages of the heavy-traffic limit is that many queueing systems behave as if they live in a much lower dimensional subspace of the state space in the limit. This phenomenon is known as State Space Collapse (SSC). If the heavy-traffic limit is taken such that exactly one resource constraint is made tight, then the system is said to satisfy the Complete Resource Pooling (CRP) condition \citep{harlop_state_space,Williams_CRP,dai2008max_pressure}.

Over the past decades, several queueing systems that satisfy the CRP condition have been successfully and extensively studied using diffusion limits and Brownian Motion processes. This approach was first developed by \cite{kingman1962_brownian}, where a $G/G/1$ queue was studied in heavy-traffic. Later, it was successfully applied in a variety of systems that satisfy the CRP condition \citep{harrison1988brownian,har_state_space,Williams_state_space,Williams_CRP,harlop_state_space,stolyar2004maxweight,gamarnik2006validity}. In this approach, a scaled version of the queue lengths process is considered, and it is shown that it converges to a Reflected Brownian Motion (RBM) process. SSC is then established to show that this RBM process lives in a lower dimensional subspace. Since the queue lengths cannot be negative, they `reflect' at the origin, so this lower dimensional Brownian Motion process is called a Reflected Brownian Motion process. Such a result is called process level convergence, and may be useful in approximating transient behavior. The next step is to obtain the stationary distribution of this RBM, which is usually the same as the heavy-traffic limiting stationary distribution of the original (unscaled) queueing system. However, this must be formally established by proving that the limit to steady-state and the limit to heavy-traffic (equivalently, limit to the RBM) can be interchanged.
Showing this interchange of limits is challenging in many systems, because one needs to establish tightness of a sequence of probability measures.
Even though this method has been successfully used to study a wide variety of problems that satisfy the CRP condition, it is challenging to study systems when the  CRP condition is not satisfied.

In addition, three different `direct methods' were developed to study queueing systems in heavy-traffic without considering the scaled process and the diffusion limits \citep{jim_achievement_lecture}. Therefore, none of these direct methods require the interchange of limits step. They are the Drift method \citep{atilla,MagSri_SSY16_Switch,QUESTA_switch,Weina_bandwidth_journal,zhou2018flexible}, Stein's method \citep{gurvich2014diffusion,braverman2017stein,braverman2017stein2} and the BAR method \citep{braverman_BAR}. We briefly describe each of them below.

The main idea in the Drift method is to carefully choose a test function, and to equate the expected value of the test function in steady-state to the same value at the following time step. Equating the expected value of the test function in two different time steps, is also known as `setting to zero the drift of the test function' (see Definition \ref{def.drift} for a formal definition of this expression). Since this method does not involve the use of diffusion limits, SSC must be established independently, and this is done using the Lyapunov drift arguments and the moment bounds developed by \cite{hajek_drift} and \cite{bertsimas_momentbound}. When selecting a test function, one needs to keep in mind that one of the reasons to perform heavy-traffic analysis is SSC. Therefore, test functions that depend on the geometry of the region where SSC occurs yield bounds that are tight in heavy-traffic. Usually, if quadratic test functions are used, bounds on the mean of the queue lengths are obtained. To obtain bounds on the $m\tth$ moments, polynomial test functions of degree $(m+1)$ are used.
The complete steady-state distribution in heavy-traffic is obtained once all the moments are obtained inductively, under some mild conditions (see Section 4.10 in \cite{gut2012probability} for a formal discussion of these conditions). For example, in the case of a single server queue, the test functions $q,q^2,q^3,\ldots$ are used inductively, where $q$ denotes the queue length.

This approach was first used to reprove known heavy-traffic results in a class of queueing systems that satisfy the CRP condition \citep{atilla}, and include a load balancing system and an ad hoc wireless network in presence of interference and fading (time-varying channel conditions). The Drift method was later successfully applied to obtain the heavy-traffic mean of the sum queue lengths even in systems that do not satisfy the CRP condition such as the input-queued switch \citep{MagSri_SSY16_Switch, QUESTA_switch} and bandwidth sharing networks \citep{Weina_bandwidth_journal}. However, it was recently shown that, when the CRP condition is not satisfied, the Drift method with polynomial test functions does
not have all the information needed to obtain  all the higher moments and the distribution of the queue lengths \citep{Hurtado_gen-switch_temp}.

In this paper we develop the Moment Generating Function (MGF) method in systems that satisfy the CRP condition, by generalizing the Drift method to directly study the stationary distribution (as opposed to the moments) in heavy-traffic. The key insight is that, instead of using the polynomial test functions of increasing degrees inductively as in the Drift method, all the polynomials can be combined in Taylor series to obtain an exponential test function. For example, in the case of a single server queue, combining $q,q^2,q^3,\ldots$ in Taylor series (with appropriate coefficients), we obtain $e^{\theta q}$ for some constant $\theta$, and $\E{e^{\theta q}}$ is the MGF of $q$. The MGF method is similar to the Drift method in the sense that we use the same notion of SSC, and that we set to zero the drift of a carefully chosen test function in steady-state. However, in the Drift method one needs to perform an inductive argument to compute the stationary distribution, whereas the MGF method immediately yields the stationary distribution.

While the Drift method is based on setting the drift of carefully chosen polynomial test functions to zero, the BAR method uses carefully chosen exponential functions. The focus in the BAR method  is to choose the exponential functions to get a handle on the jumps in a continuous time system under general arrivals and services. In this paper, we illustrate how the MGF method can be thought of as a natural generalization of the Drift method using exponential test functions, and in that sense is similar in spirit to the BAR method. Using the BAR method, it was shown by  \cite{braverman_BAR},  that in heavy-traffic, the  stationary distribution of a Generalized Jackson Network is identical to that of an appropriately defined RBM. In contrast, the focus in the current paper is to incorporate SSC and to evaluate the closed form stationary distribution in heavy-traffic in a variety of systems under the CRP condition. Moreover, while the BAR method was developed to study continuous time systems, we focus on studying discrete time systems in this paper.

The Drift method and the BAR method are focused on computing the stationary distribution of the scaled queue lengths in heavy-traffic. On the other hand, Stein's method is focused on computing rates of convergence to the limiting distribution. Stein's method for studying queueing systems was first introduced by \cite{gurvich2014diffusion}. Erlang-A and Erlang-C queueing models were studied using Stein's method by \cite{braverman2017stein}, and $M/Ph/n+M$ systems by \cite{braverman2017stein2}.
Similar to the MGF method, a key step in using Stein's method for some results is in establishing SSC.
Stein's Method was used to study load balancing systems in mean field regime
\citep{Lei_Stein_SIGM16,Lei_steinHT_SIGM17},  in Halfin-Whitt regime in \citep{braverman2018jsq}, and in sub-Halfin-Whitt regimes in \citep{liu2018simple}. Universal approximations for queues with abandonment were obtained using Stein's method by \cite{huang_gurvich2018universality}. More recently, a single server queue in heavy-traffic was studied using Stein's method by \cite{Walton_SteinHT}. \cite{gurvich2013excursion} studies Erlang-A system and obtains universal approximations using excursion-based analysis, as opposed to using Stein's method.

In this paper, we develop the MGF method and illustrate its power to study a variety of queueing systems that satisfy the CRP condition. In order to introduce the method, and to showcase its simplicity, we first present a sketch of the MGF method in the case of a single server queue operating in discrete time in Section \ref{sec:single.server.queue}. We show that the stationary distribution of scaled queue length in heavy-traffic limit converges to an exponential distribution. This is of course a classic result first proved by \cite{kingman1962_brownian} using the diffusion limit method, and later by \cite{atilla} using the Drift method.

We then develop the MGF method framework and  apply it to load balancing systems and generalized switches. In both cases we study the queueing systems under some general conditions and we exemplify with specific systems that satisfy those conditions. In Section \ref{sec:load.balancing} we study load balancing systems and identify that the  Join the Shortest Queue (JSQ) \citep{JSQ_HT_optimality,winston_JSQ_1977} and power-of-two choices \citep{dobrushin_po2,mitzenmacher_po2,mitzenmacher_po2_2} routing policies satisfy the assumptions. In Section \ref{sec:generalized.switch} we study generalized switches \citep{stolyar2004maxweight} under the CRP condition, operating under MaxWeight scheduling algorithm \citep{TasEph_92}. We also show that ad hoc wireless networks operating under MaxWeight scheduling algorithm satisfy our assumptions. All these systems are assumed to satisfy the CRP condition, and they are operated under algorithms that ensure that SSC occurs into a one-dimensional subspace. We show that the stationary distribution of this one-dimensional component is exponential. In addition to Moment Generating Functions, which are the two-sided Laplace transforms of the probability distribution, one may use other transforms such as one-sided Laplace transforms and characteristic functions. We present a brief discussion about other transform methods in Remark \ref{remark.other.transforms}, at the end of Section \ref{sec:details.MGFmethod}.

The primary contribution of this paper is the development of the MGF method, which is a simple framework to compute the stationary distribution of the scaled vector of queue lengths in heavy-traffic. This is done  by considering the above mentioned set of systems. The paper also shows how the MGF method can be thought of as a generalization of the Drift method by considering a richer class of test functions. This class of test functions leads to substantially different proofs, that are much simpler than in the Drift method, as will be illustrated in the following sections. However, unlike the Drift method, the MGF method does not involve an art of picking a test function, since the test function is essentially the MGF. Even though most of the results that we present have already been established in the literature using diffusion limit and drift methods, the purpose of this paper is to develop a framework based on transform techniques and illustrate its power and simplicity.
A secondary contribution is that the load balancing system we consider is allowed to have  correlated servers and the generalized switch is allowed to have correlated arrival processes. Under the CRP condition and control algorithms that ensure SSC to a one-dimensional subspace, we show that even under correlated arrivals or services, the heavy-traffic scaled stationary distribution continues to be exponential (Theorems \ref{thm.load.balancing} and \ref{thm.generalized.switch}, respectively). It is possible to allow for this generalization using other methods, but we illustrate the simplicity of such generalizations using the MGF method.

The focus of this paper is on queueing systems that satisfy the CRP condition. However, the long-term goal of developing the MGF method is to characterize the heavy-traffic stationary distribution  of systems that do not satisfy the CRP condition, such as input-queued switches \citep{MagSri_SSY16_Switch,QUESTA_switch}. This will form the basis for future work on input-queued switches, which is briefly discussed in Section \ref{sec:future.work}. This approach is similar to the one taken in the development of the Drift method, which was first proposed by \cite{atilla} to prove known results in systems under the CRP condition. The Drift method was later generalized to  study the input-queued switch when CRP condition is not satisfied \citep{MagSri_SSY16_Switch,QUESTA_switch}.

\subsection{Notation}

In this section we introduce the notation that we will use along the paper. We use $\Prob{A}$ to denote the probability of the event $A$,  $\E{X}$ to denote the expected value of the random variable $X$, $\Cov{X,Y}$ to denote the covariance between the random variables $X$ and $Y$ and $\Var{X}$ to denote the variance of the random variable $X$. The indicator function of an event $A$ is $\ind{A}$, i.e., $\ind{A}$ is one if $A$ is true and 0 otherwise. Convergence in distribution is denoted by $\Longrightarrow$.

We use $\bR$ to denote the set of real numbers and $\bR^n$ to denote the set of $n$-dimensional vectors with real components. We use $\bR_+$ and $\bR_+^n$ to denote the set of nonnegative numbers and the set of $n$-dimensional vectors with nonnegative elements, respectively. Vectors are written in bold letters and we use the same letter, but not bold and with a subindex, to denote their elements. For example, for a positive integer $n$, the vector $\vx\in\bR^n$ has elements $x_i\in\bR$ for $i\in \{1,\ldots,n\}$. We use $\vone$ to denote a vector of ones and $\vzero$ to denote a vector of zeroes, i.e., if $\vx=\vone$, then $x_i=1$ for all $i\in\{1,\ldots,n\}$ and if $\vx=\vzero$, then $x_i=0$ for all $i\in\{1,\ldots,n\}$. The dot product of two vectors $\vx$ and $\boldsymbol{y}$ is denoted by $\langle \vx,\boldsymbol{y}\rangle$ and the Euclidean norm is denoted by $\|\vx\|$, i.e., $\|\vx\|=\sqrt{\langle \vx,\vx\rangle}$. For each $i\in\{1,\ldots,n\}$ we use $\ve^{(i)}$ to denote the $i\tth$ canonical vector, i.e., a vector with elements $e^{(i)}_i=1$ and $e^{(i)}_j=0$ for all $j\neq i$. Given a fixed vector $\vc\in\bR^n$ and a parameter $b\in\bR$, the set $\{\vx\in\bR^n: \langle\vc,\vx\rangle = b\}$ is a hyperplane and the set $\{\vx\in\bR^n: \langle\vc,\vx\rangle \leq b\}$ is a half-space.

We say $f(x)$ is $O\left(g(x)\right)$ if $\ds\lim_{x\dto 0}\left|\dfrac{f(x)}{g(x)} \right|$ is finite and we say that $f(x)$ is $o\left(g(x)\right)$ if $\ds\lim_{x\dto 0}\dfrac{f(x)}{g(x)}=0$.

\section{Related Work}\label{sec:related.work}

In this section, we present an overview of related work on heavy-traffic analysis of queueing systems in general, as well as the different systems that we will study in particular.

Moment Generating Functions have been used in the literature to study queueing systems such as the classical analysis of $M/G/1$ queue \citep{harrison1992performance}. However, here we use the MGF to study heavy-traffic scaled queue lengths, since the queue lengths go to infinity in the heavy-traffic limit. There has been only a little work in the literature that uses Transform Methods for heavy-traffic analysis.
Characteristic Functions were used by \cite{kollerstrom1974heavy} and \cite{kingman1961_charfunction} to study heavy-traffic queueing systems, and moment generating functions were used  by \cite{lehoczky1996real,lehoczky1997using}. In contrast, the primary focus of this work is to develop transform methods for heavy-traffic analysis that incorporate SSC.

The single server queue was first studied in heavy-traffic by \cite{kingman1961_charfunction} using Characteristic Functions and tools from complex analysis. \cite{kollerstrom1974heavy} also used
Characteristic Functions to study single server queue.
The diffusion limit method to study queueing systems was developed by studying the single server queue \citep{kingman1962_brownian}. The well known Kingman bound for the expected waiting time in a single server queue was developed in the 60's \citep{kingman}, and later \cite{marshall_ineq_drift} used similar arguments to compute bounds on the second moment. These formed the basis for the Drift method, that was developed by \cite{atilla}. The single server queue was also presented as an illustrative example of the BAR method \citep{braverman_BAR}. Most of these papers study the delay in $G/G/1$ queue in continuous time, which evolves according to Lindley's equation \citep{lindley_equation}. Similar to \cite{atilla}, in this paper we study the queue length in discrete time. The queue lengths process evolves according to  \eqref{eq.qu}, which is equivalent to Lindley's equation for the waiting time of $(k+1)\tth$ customer in a $G/G/1$ queue. Consequently, the results established for queue lengths in discrete time can be easily extended to delay in continuous time.

The load balancing system (also known as the supermarket checkout model) has been widely studied since the 70's. It was shown that the JSQ policy minimizes the mean delay among the class of policies that do not know the job durations \citep{winston_JSQ_1977,weber1978optimal,ephremides1980simple}. Heavy-traffic optimality of the JSQ policy in a system with two servers was established by \cite{JSQ_HT_optimality} using the diffusion limit method, where they also introduced the notion of SSC. Since then, the load balancing system has been extensively studied both to improve performance and to decrease the complexity of the algorithms \citep{chen2012po2_journal,li2018loadbalancing,braverman2018jsq,zhou2018flexible,Lei_Stein_SIGM16,Lei_steinHT_SIGM17, Gamarnik_JSQ,lu_JIQ,stolyar_JIQ,ying_power_more_than2}. One lower complexity algorithm that has received attention is the power-of-two choices algorithm \citep{dobrushin_po2,mitzenmacher_po2,mitzenmacher_po2_2,chen2012po2_journal}. An exhaustive survey of literature on load balancing is presented by \cite{load_balancing_survey}. The most relevant work for our purposes is the study of the JSQ policy under the Drift method by \cite{atilla} and that of the power-of-two choices algorithm by \cite{magsriyin_itc12_journal}.

MaxWeight algorithm was first proposed by \cite{TasEph_92} in the context of scheduling for down-links in wireless base stations. This algorithm was later adapted to be used in a variety of systems including ad hoc wireless networks, input-queued switches \citep{mckeown96walrand}, cloud computing \citep{magsriyin_itc12_journal}, was generalized into the back-pressure algorithm \citep{TasEph_92} in networks, and was extensively studied by \cite{stolyar2004maxweight}, \cite{gupta2010delay}, and \cite{Mey_08}. The generalized switch model subsumes many of these systems, and has been studied under the CRP condition using the diffusion limit method \citep{stolyar2004maxweight}, and the Drift method \citep{atilla}. We use the notion of SSC as developed by \cite{atilla}. \cite{dai2008max_pressure} generalizes the results in \cite{stolyar2004maxweight} to SPNs where the jobs can join a queue after being served.

\section{The MGF method}\label{sec:MGFmethod}
In this section we introduce the MGF method to compute the distribution of scaled queue lengths in heavy-traffic. This section is organized as follows. In Section \ref{sec:MGF.general.model} we define a general queueing model; in Section \ref{sec:single.server.queue} we introduce the method with a single server queue, as a simple example; and in Section \ref{sec:details.MGFmethod} we describe the  MGF method as a step by step procedure, so that it can be applied in the context of a variety of queueing systems.

\subsection{A general queueing model}\label{sec:MGF.general.model}

We first introduce a general queueing model for an SPN that includes the single server queue, the load balancing system and the generalized switch as special cases. We provide the details of each system in the corresponding section.

Consider a single hop queueing system operating in discrete time, with $n$ separate servers. Each server has an infinite buffer, where jobs line up if the server is busy. For $k\geq 1$ and $i\in\{1,\ldots,n\}$ let $q_i(k)$ be the number of jobs in the $i\tth$ queue at the beginning of time slot $k$, i.e., the number of jobs waiting to be served and the job that is being served (if any). Let $\vq(k)$ be an $n$-dimensional vector with elements $q_i(k)$ for $i\in\{1,\ldots,n\}$. Given that the vector of queue lengths in time slot $k$ is $\vq(k)$, let $a_i\left(\vq(k)\right)$ be the number of arrivals to the $i\tth$ queue in time slot $k$ and $s_i\left(\vq(k)\right)$ be the potential number of jobs that can be served from queue $i$ in time slot $k$. We say $s_i\left(\vq(k)\right)$ is potential service because, if there are not enough jobs in line, then less than $s_i\left(\vq(k)\right)$ jobs are processed. For ease of exposition, and with a slight abuse of notation, from now on we write $\va(k)$ and $\vs(k)$ instead of $\va\left(\vq(k)\right)$ and $\vs\left(\vq(k)\right)$, respectively.
We assume that $a_i(k)$ and $s_i(k)$ are upper bounded by constants. Specifically, let $\amax$ and $\smax$ be finite constants such that $a_i(k)\leq \amax$ and $s_i(k)\leq \smax$ with probability 1 for all $i\in\{1,\ldots,n\}$ and all $k\geq 1$. The difference between potential and actual service is called unused service, which we denote $u_i\left(\vq(k)\right)$. We also write $\vu(k)$ instead of $\vu\left(\vq(k) \right)$ from now on, for ease of exposition. In some queueing systems, the control problem is to decide the vector $\va(k)$ in each time slot (e.g. the load balancing system) and, in others the vector $\vs(k)$ (e.g. the generalized switch). We give more details about these selection processes in the systems that we study in Sections \ref{sec:load.balancing} and \ref{sec:generalized.switch}, respectively.

In each time slot, the order of events is as follows. First, queue lengths are observed. Second, given the vector of queue lengths $\vq(k)$, the control problem is solved. Then, arrivals occur and, at the end of each time slot, jobs are processed by the servers. Therefore, the dynamics of the queues are as follows
\begin{align}\label{eq.dynamics.bracketplus}
q_i(k+1)=& \max\big\{q_i(k)+a_i(k)-s_i(k),0\big\}\qquad\forall i\in\{1,\ldots,n\},\; \forall k\geq 1.
\end{align}
For each $i\in\{1,\ldots,n\}$ the variables $a_i(k)$ and $s_i(k)$ depend only on $\vq(k)$, (or they are independent of $\vq(k)$), then \eqref{eq.dynamics.bracketplus} implies that the process $\{\vq(k):\;k\geq 1\}$ is a Markov chain.

We can also describe the dynamics of the queues using unused service instead of the maximum, as follows
\begin{align}\label{eq.dynamics}
q_i(k+1)=& q_i(k)+a_i(k)-s_i(k)+u_i(k)\qquad\forall i\in\{1,\ldots,n\},\; \forall k\geq 1.
\end{align}
Observe that \eqref{eq.dynamics} implies
\begin{align}\label{eq.qu}
q_i(k+1)u_i(k)=0\qquad\forall i\in\{1,\ldots,n\},\;\forall k\geq 1
\end{align}
because the unused service is nonzero only when the potential service is larger than the number of jobs available to be served (queue length and arrivals), and in this case the queue is empty in the next time slot. If $i\neq j$, then we do not necessarily have $q_i(k+1)u_j(k)=0$ because the fact that queue $j$ is empty at the end of time slot $k$ does not imply that queue $i$ will be empty at the beginning of time slot $k+1$, and vice versa.
It turns out that getting a handle on the unused service plays an important role in heavy-traffic analysis and \eqref{eq.qu} will be an important tool in the analysis. Equation \eqref{eq.qu} can be thought of as a defining property of the queueing process and is analogous to the Skorohod map \citep{Skorohod_map}.

In this paper we add a line on top of the variables and vectors to denote steady-state. Specifically, let $\vqbar$, $\vabar\defn\va(\vqbar)$, $\vsbar\defn\vs(\vqbar)$ and $\vubar\defn\vu(\vqbar)$ 
be steady-state vectors that represent the queue lengths at the beginning of a time slot, and arrivals, potential service and unused service in one time slot in steady-state, respectively. 
Let $\vqbar^+ \defn \vqbar+\vabar-\vsbar+\vubar$ denote the queue length at time $k+1$ in terms of the queue length, arrival and service at time $k$, assuming the system is in steady-state. The precise definition of each of these steady-state vectors depends on the control problem, so we provide them in Section \ref{sec:single.server.queue} for the single server queue, in Section \ref{sec:load.balancing} for the load balancing system and in Section \ref{sec:generalized.switch} for the generalized switch.

The MGF method will be used to compute the joint distribution of the scaled vector of queue lengths in heavy-traffic, so before introducing the framework we specify what we mean by heavy-traffic and how we parametrize the queueing systems to obtain the limit. The heavy-traffic limit is the limit as the arrival rate vector approaches the boundary of the capacity region of the system. The capacity region of an SPN is the set of arrival rate vectors such that the system can be positive recurrent. In other words, if the arrival rate vector is in the interior of the capacity region, there exists an algorithm that solves the control problem and is such that the queue length process is positive recurrent; if the vector of arrival rates is outside the capacity region, no algorithm can ensure positive recurrence. We use $\cC$ to denote the capacity region and we parametrize the heavy-traffic limit as follows. Take $\epsilon>0$ and consider a set of queueing systems parametrized by $\epsilon$. The parametrization is such that $\epsilon$ represents how far away the vector of arrival rates is from a fixed point $\vr$ in the boundary of $\cC$. Then, the heavy-traffic limit is the limit as $\epsilon\dto 0$. In this paper we add a superscript $\peps$ when we refer to the parametrized queueing system. More details on the parametrization of each queueing system will be provided once the models are completely specified, i.e., in Section \ref{sec:single.server.queue} for the single server queue, in Section \ref{sec:load.balancing} for the load balancing system and in Section \ref{sec:generalized.switch} for the generalized switch.

Before introducing the MGF framework in the context of a single server queue we formally define the drift of a function and we explain what `set the drift to zero' means. 

\begin{definition}[Drift of a function]\label{def.drift}
	Let $V:\bR^n\to \bR_+$ be a function. We define the drift of $V$ at $\vq$ as
	\begin{align*}
		\Delta V(\vq) \defn \left(V\big(\vq(k+1)\big)- V\big(\vq(k)\big) \right)\ind{\vq(k)=\vq}.
	\end{align*}
	If $\E{V\big(\vq(k)\big)}<\infty$ for $k$ such that the Markov chain $\left\{\vq(k):k\geq 1 \right\}$ is in steady-state, we say that we set the drift of $V$ to zero when we use the property
	\begin{align*}
		\E{\Delta V\left(\vq(k)\right)}=0,
	\end{align*}
	where the expected value is taken with respect to the stationary distribution.
\end{definition}

Observe that we can set to zero the drift of any function with finite expected value, by definition of steady-state.

\subsection{MGF method in the single server queue}\label{sec:single.server.queue}

Before presenting the details of the MGF framework, we use it in the simplest queueing system: a single server queue. We provide a proof of the well-known result that the scaled queue length of a single server queue has an exponential distribution in heavy-traffic to illustrate the method and to show its simplicity. We do not provide all the details of our proofs, since the single server queue is a special case of the load balancing system ($n=1$) and this system is studied in detail in Section \ref{sec:load.balancing}.

Consider a single server queue operating in discrete time. Arrivals and potential service in each time slot are assumed to be independent sequences of i.i.d. random variables. Since they are also assumed to be finite with probability 1 (as specified in Section \ref{sec:MGF.general.model}), their MGFs $\E{e^{\theta a(1)}}$ and $\E{e^{\theta s(1)}}$ exist for all $\theta\in\bR$.

Let $\lambda\defn\E{a(1)}$ and $\mu\defn\E{s(1)}$. Observe that $\lambda$ and $\mu$ are the rates of arrival and service, respectively, since they are the expected number of arrival/services in one time slot. Then, the capacity region of the single server queue is $\cC=\{\lambda\in\bR_+:\; \lambda\leq \mu\}$. We consider a set of single server queues parametrized by $\epsilon$ with a fixed service process of rate $\mu$ and arrival rate $\lambda^\peps\defn \mu-\epsilon$.

Let $\abar^\peps$ and $\sbar$ be steady-state random variables that have the same distribution as $a^\peps(1)$ and $s(1)$, respectively.  Then, $\lambda^\peps=\E{\abar^\peps}$ and $\mu=\E{\sbar}$. Let $\left(\sigma_a^\peps\right)^2=\Var{\abar^\peps}$ and $\sigma_s^2=\Var{\sbar}$. 

In the rest of this section we prove Theorem \ref{thm.single.server.queue}. This is a well-known result and there are proofs using diffusion limits \citep{kingman1962_brownian} and the Drift method \citep{atilla} in the literature. We present an alternate proof which is simpler than the two proofs mentioned above, and will serve as a template for the MGF method.

\begin{theorem}\label{thm.single.server.queue}
	Let $\epsilon\in(0,\mu)$ and consider a set of single server queues parametrized by $\epsilon$ as described above. Let $\qbar^\peps$ be a steady-state random variable such that $\{q^\peps(k):\;k\geq 1\}$ converges in distribution to $\qbar^\peps$ as $k\uparrow\infty$. Further, assume $\lim_{\epsilon\dto 0}\sigma_a^\peps=\sigma_a$. Then, $\epsilon\qbar^\peps\Rightarrow \Upsilon$ as $\epsilon\dto 0$, where $\Upsilon$ is an exponential random variable with mean $\frac{\sigma_a^2+\sigma_s^2}{2}$.
\end{theorem}

It is well-known that for all $\epsilon\in(0,\mu)$, the Markov chain $\{q^\peps(k):k\geq 1 \}$ is positive recurrent. For instance, the reader can find a proof using Foster-Lyapunov theorem in \citep[Theorem 3.4.2]{srikantleibook}. Then, $\qbar^\peps$ is well defined.

Before presenting the proof, we prove two lemmas. The first lemma is a different version of \eqref{eq.qu} and is key in the MGF method.	For other queueing systems we use a weaker version of this lemma, that is sufficient for the MGF method (see Step 1 in Section \ref{sec:details.MGFmethod} for more details).

\begin{lemma}\label{lemma.ssq.expoqu}
	Consider a single server queue parametrized by $\epsilon$ as described above. Then, for all $\alpha,\beta\in\bR$ and each $k\geq 1$ we have
	\begin{align*}
	\left(e^{\alpha q^\peps(k+1)}-1 \right)\left(e^{-\beta u^\peps(k)}-1\right)=0.
	\end{align*}
\end{lemma}
\proof{Proof of Lemma \ref{lemma.ssq.expoqu}.}
It follows from \eqref{eq.qu} and because $e^x-1=0$ if and only if $x=0$. \Halmos
\endproof

The next Lemma gives the first moment of the unused service in steady-state, and it will be used at the end of the proof of Theorem \ref{thm.single.server.queue}.

\begin{lemma}\label{lemma.ssq.Eu}
	Consider a single server queue parametrized by $\epsilon\in(0,\mu)$ as described above. Then,
	\begin{align*}
		\E{\ubar^\peps}=\epsilon.
	\end{align*}
\end{lemma}
\proof{Proof of Lemma \ref{lemma.ssq.Eu}.}

We set to zero the drift of the linear test function $V_1(q)=q$, and we obtain
\begin{align*}
0=& \E{\left(\qbar^\peps\right)^\plus-\qbar^\peps} \\
=& \E{(\qbar^\peps+\abar^\peps-\sbar+\ubar^\peps)-\qbar^\peps},
\end{align*}
where the last equality holds by definition of $\left(\qbar^\peps\right)^+$. Rearranging terms we obtain
\begin{align*}
\E{\ubar^\peps}=& \E{\sbar-\abar^\peps} = \mu-(\mu-\epsilon)= \epsilon.
\end{align*}

\Halmos
\endproof

Now we prove Theorem \ref{thm.single.server.queue}.

\proof{Proof of Theorem \ref{thm.single.server.queue}.}

If we expand the product in Lemma \ref{lemma.ssq.expoqu} and rearrange terms we obtain
\begin{align}\label{eq.ssq.expo.qu.expanded}
	e^{\theteps q^\peps(k+1)}- e^{\theteps\left(q^\peps(k)+a^\peps(k)-s(k) \right)}=& 1-e^{-\theteps u^\peps(k)}
\end{align}

Observe that \eqref{eq.ssq.expo.qu.expanded} holds for all $k\geq 1$. In particular, it holds in steady-state. Also, it can be shown that $\E{e^{\theteps\qbar^\peps}}<\infty$ in an interval around 0. We omit the proof because in Lemma \ref{lemma:mgf-jsq} we provide a proof for the load balancing system, which is a more general case. Therefore, $\E{e^{\theteps\left(\qbar^\peps \right)^+}}= \E{e^{\theteps \qbar^\peps}}$. Taking expected value with respect to the stationary distribution in \eqref{eq.ssq.expo.qu.expanded} we obtain
\begin{align*}
	\E{e^{\theteps \qbar^\peps}\left(1- e^{\theteps\left(\abar^\peps-\sbar \right)}\right)}=& 1-\E{e^{-\theteps \ubar^\peps}}.
\end{align*}
Since $\abar^\peps$ and $\sbar$ are independent of the queue length, rearranging terms we obtain
\begin{align}\label{eq.ssq.fraction}
\E{e^{\theteps \qbar^\peps}}=\dfrac{1-\E{e^{-\theteps \ubar^\peps}}}{1-\E{e^{\theteps \left(\abar^\peps-\sbar\right)}}}
\end{align}

Now we take the heavy-traffic limit. Observe that the right hand side yields a $\frac{0}{0}$ form in the limit as $\epsilon\dto 0$. Then, we take Taylor series of each term with respect to $\theta$, around $\theta=0$. The technical details of why this expansion can be done are established in Lemma \ref{lemma.taylor}, which is presented in Section \ref{sec:details.MGFmethod} . For the numerator we obtain
\begin{align}
	1-\E{e^{-\theteps \ubar^\peps}} =& \theteps\E{\ubar^\peps}-\dfrac{(\theteps)^2}{2}\E{\left(\ubar^\peps\right)^2}+O(\epsilon^3) \nonumber \\
	=& \theteps^2+O(\epsilon^3), \label{eq.ssq.taylor.u}
\end{align}
where the last equality holds by Lemma \ref{lemma.ssq.Eu} and because $\E{\left(\ubar^\peps\right)^2}$ is $O(\epsilon)$. Details of this argument will be provided in Section \ref{sec:load.balancing} for the load balancing system (see Claim \ref{claim.load.balancing}), but the main idea is that $\ubar^\peps$ is a bounded random variable. For the denominator we obtain
\begin{align}
1-\E{e^{\theteps(\abar^\peps-\sbar)}} \nonumber
=& -\theteps\E{\abar^\peps-\sbar}-\dfrac{(\theteps)^2}{2}\E{(\abar^\peps-\sbar)^2}+O(\epsilon^3) \nonumber  \\
=& \theteps^2-\dfrac{(\theteps)^2}{2}\left(\left(\sigma_a^\peps\right)^2+\sigma_s^2+\epsilon^2 \right) +O(\epsilon^3), \label{eq.ssq.taylor.as}
\end{align}
where the last step holds because $\E{\abar^\peps}=\mu-\epsilon$ and by definition of variance.

If we replace \eqref{eq.ssq.taylor.u} and \eqref{eq.ssq.taylor.as} in \eqref{eq.ssq.fraction}, and cancel out $\theteps^2$ from numerator and denominator we obtain
\begin{align*}
\E{e^{\theteps\qbar^\peps}}= \dfrac{1+O(\epsilon) }{1-\dfrac{\theta}{2}\left(\left(\sigma_a^\peps\right)^2+\sigma_s^2\right)+O(\epsilon)}
\end{align*}

Therefore, taking the heavy-traffic limit we obtain
\begin{align}\label{eq.ssq.limit}
\lim_{\epsilon\dto 0}\E{e^{\theteps \qbar^\peps}}=\dfrac{1}{1-\theta\left(\frac{\sigma_a^2+\sigma_s^2}{2} \right)}
\end{align}

Since the right hand side is the MGF of an exponential random variable with mean $\frac{\sigma_a^2+\sigma_s^2}{2}$, Equation \eqref{eq.ssq.limit} implies $\epsilon\qbar^\peps$ converges in distribution to such an exponential random variable \cite[Theorem 9.5 in Section 5]{gut2012probability}. \Halmos
\endproof

In this section we exemplified the MGF method in an intuitive fashion for the simplest queueing system. In the next subsection we generalize these steps for other queueing systems that satisfy the CRP condition.

\subsection{General framework}\label{sec:details.MGFmethod}

In the last subsection we proved a well-known result using the MGF method in the case of the simplest queueing system, i.e., the single server queue. In this subsection we describe the method in detail for more general queueing systems that satisfy the CRP condition. Before presenting the framework, we present a formal definition of the CRP condition. We use the definition provided by \cite{stolyar2004maxweight}.

\begin{definition}[CRP condition]\label{def.CRP}
	Consider a set of queueing systems parametrized by $\epsilon$ as described in Section \ref{sec:MGF.general.model}, where the  capacity region is $\cC$. Suppose that in heavy-traffic (i.e., as $\epsilon\dto 0$), the vector of arrival rates approaches a point $\vr$ in the boundary of $\cC$. We say that the queueing system satisfies the Complete Resource Pooling (CRP) condition if the outer normal vector to $\cC$ at $\vr$ is unique up to a scalar coefficient.
\end{definition}

This implies that the system can be operated such that all the servers pool together in the heavy-traffic limit \citep{harlop_state_space,dai2008max_pressure,Williams_CRP}. Intuitively, this means that the queueing system behaves as a one-dimensional queueing system (i.e. as a single server queue) if it is operated under a `good' control algorithm. Therefore, the MGF method is essentially similar to the proof of Theorem \ref{thm.single.server.queue} after one establishes SSC on a one-dimensional subspace of the state space.

In order to use the MGF method, one needs to make sure that two prerequisites are satisfied. We state them before presenting the framework.

\subsubsection*{Prerequisite 1. Positive recurrence.}
Prove that the Markov chain $\{\vq^\peps(k):\;k\geq 1\}$ is positive recurrent for $\epsilon>0$.

Positive recurrence is a requirement to make sure there exists a steady-state random vector $\vqbar^\peps$ such that the queue lengths process $\{\vq^\peps(k):\;k\geq 1\}$ converges in distribution to $\vqbar^\peps$ as $k\uparrow \infty$.

\subsubsection*{Prerequisite 2. State Space Collapse.} Prove SSC into a one-dimensional subspace.
	
Let $\vc\geq \vzero$ be the direction into which SSC occurs. For simplicity, we assume $\|\vc\|=1$. Then $\cK = \left\{\vy\in\bR^n: \vy = \alpha \vc\,,\,\alpha\geq 0 \right\}$ is the cone where the state space collapses in heavy-traffic. For any $n$-dimensional vector $\vx$, let $\vx_\parallel\defn \langle\vx,\vc\rangle\vc$ be the projection of $\vx$ on $\cK$ and let $\vx_\perp\defn \vx-\vx_\parallel$. In this step it should be proved that $\E{\left\|\vqperpbar^\peps \right\|^2}$ is $o\left(\frac{1}{\epsilon^2}\right)$, which is equivalent to proving that $\epsilon^2 \E{\left\|\vqperpbar^\peps \right\|^2}$ is $o(1)$.

The queueing systems that we study in this paper actually exhibit a stronger form of SSC, where $\E{\left\|\vqperpbar^\peps \right\|^m}$ is $O(1)$ for all $m=1,2,\ldots$ However, a weaker form of SSC is studied by \cite{Weina_bandwidth_journal} and \cite{javidi_optical_ht}.

From this notion of SSC, we conclude that
\begin{align*}
	\lim_{\epsilon\dto 0}\epsilon^2 \E{\left\|\vqperpbar^\peps \right\|^2}=0,
\end{align*}
i.e., $\epsilon \left\|\vqperpbar^\peps \right\|$ converges to zero in the mean squares sense and, therefore, in probability.

In the case of the single server queue we did not have to verify Prerequisite 2, because the state space is already one-dimensional. Now we present the MGF method.

\subsubsection*{Step 1. Prove an equation of the form
\begin{align}\label{eq.MGF.method.exp.qu}
	\E{\left(e^{\theteps \langle\vc,\left(\vqbar^\peps\right)^\plus \rangle} -1\right)\left(e^{-\theteps \langle\vc,\vubar^\peps\rangle}-1 \right)} \quad\text{is $o(\epsilon^2)$}
\end{align}
and compute an expression for the MGF of $\epsilon\langle\vc,\vqbar^\peps\rangle$.}

The key in the MGF method is to handle unused service and its interaction with the queue lengths, arrivals and potential service. In the Drift method, the unused service is handled with \eqref{eq.qu}. However, in this case we want to work with an exponential transform of the queue lengths, so we need to write \eqref{eq.qu} in a different way. In the case of the single server queue, we used Lemma \ref{lemma.ssq.expoqu} which, in fact, it is much stronger than what we actually use in the MGF method. For more general queueing systems we use \eqref{eq.MGF.method.exp.qu}.

To prove an equation of the form of \eqref{eq.MGF.method.exp.qu} it is essential to use SSC. After proving \eqref{eq.MGF.method.exp.qu}, we need to obtain an expression for the MGF of $\epsilon\langle\vc,\vqbar^\peps\rangle$ that is valid for all traffic. Below we sketch some algebraic steps that are useful to do it. Expanding the product in the left hand side of \eqref{eq.MGF.method.exp.qu} we obtain
\begin{align}
& \E{\left(e^{\theteps \langle\vc,\left(\vqbar^\peps\right)^\plus \rangle} -1\right)\left(e^{-\theteps \langle\vc,\vubar^\peps\rangle}-1 \right)} \nonumber \\
=& \E{e^{\theteps \langle\vc, \left(\vqbar^\peps\right)^+-\vubar^\peps\rangle}}- \E{e^{\theteps \langle\vc,\left(\vqbar^\peps\right)^+\rangle}} + 1-\E{e^{-\theteps \langle\vc,\vubar^\peps\rangle}} \label{eq.MGFmethod.expo.first.step} \\
\stackrel{(a)}{=}& \E{e^{\theteps \langle\vc, \vqbar^\peps+\vabar^\peps-\vsbar^\peps\rangle}}- \E{e^{\theteps \langle\vc,\left(\vqbar^\peps\right)^+\rangle}}+ 1 -\E{e^{-\theteps \langle\vc,\vubar^\peps\rangle}} \nonumber \\
\stackrel{(b)}{=}& \E{e^{\theteps \langle\vc, \vqbar^\peps+\vabar^\peps-\vsbar^\peps\rangle}}- \E{e^{\theteps \langle\vc,\vqbar^\peps\rangle}}+ 1 -\E{e^{-\theteps \langle\vc,\vubar^\peps\rangle}}, \label{eq.MGFmethod.expo.expanded}
\end{align} 
where $(a)$ holds by the dynamics of the queues described in \eqref{eq.dynamics} and by definition of $\left(\vqbar^\peps\right)^+$; and $(b)$ holds if the MGF of $\epsilon\langle\vc,\vqbar^\peps\rangle$ exists in an interval around 0 (this must be proved). In such case, by definition of steady-state we have $\E{e^{\theteps \langle\vc,\left(\vqbar^\peps\right)^+\rangle}}=\E{e^{\theteps \langle\vc,\vqbar^\peps\rangle}}$, which is equivalent to setting to zero the drift of the test function $V(\vq)=e^{\theteps \langle\vc,\vq\rangle}$.

Observe that when we first expand the product in \eqref{eq.MGFmethod.expo.first.step}, we obtain two terms that are related to the unused service (the first and the last term). We use \eqref{eq.dynamics} to deal with the first one, and we write $\left(\vqbar^\peps\right)^+-\vubar^\peps$ in terms of $\vqbar^\peps$, $\vabar^\peps$ and $\vsbar^\peps$. The last term is the MGF of $\epsilon\langle\vc,\vubar^\peps\rangle$, and we deal with it in the second step of the MGF method.

Using \eqref{eq.MGFmethod.expo.expanded} in \eqref{eq.MGF.method.exp.qu} and reorganizing terms we obtain
\begin{align}\label{eq.MGFmethod.fraction}
	\E{e^{\theteps\langle\vc,\vqbar^\peps\rangle}\left(1- e^{\theteps\langle\vc, \vabar^\peps-\vsbar^\peps\rangle}\right)} = 1- \E{e^{-\theteps\langle\vc,\vubar^\peps\rangle}} + o(\epsilon^2)
\end{align}

From \eqref{eq.MGFmethod.fraction} we can obtain an expression for the MGF of $\epsilon\langle\vc,\vqbar^\peps\rangle$ which is valid for all traffic. However, the steps to obtain it depend on the properties of each queueing system. For example, in the case of the single server queue we know that the arrival and potential service processes are independent of the queue lengths. Then, we can separate the product on the left hand side and we obtain \eqref{eq.ssq.fraction}.

\subsubsection*{Step 2. Bound unused service and take heavy-traffic limit.}

Observe that the MGF of $\epsilon\langle\vc,\vabar^\peps\rangle$ and $\epsilon\langle\vc,\vsbar^\peps\rangle$ exist for all $\theta\in\bR$, because the random variables are bounded by assumption. Further, by definition of unused service, we have $\vzero \leq \vubar^\peps\leq \vsbar^\peps$ component-wise. Then, the MGF of $\epsilon\langle\vc,\vubar^\peps\rangle$ exists for all $\theta\in\bR$. Also, in Step 1 (before obtaining \eqref{eq.MGFmethod.expo.expanded}) it was proved that the MGF of $\epsilon\langle\vc,\vqbar^\peps\rangle$ exists in an interval around zero. Therefore, as $\epsilon\dto 0$, Equation \eqref{eq.MGFmethod.fraction} yields $0=0$. As mentioned above, depending on the queueing system we will use different approaches to obtain an expression for the MGF of $\epsilon\langle\vc,\vqbar^\peps\rangle$ that is valid for all traffic from \eqref{eq.MGFmethod.fraction}. For example, in the case of the single server queue we obtained \eqref{eq.ssq.fraction}, which yields a $\frac{0}{0}$ form in the limit as $\epsilon\dto 0$. Therefore, to compute the heavy-traffic limit we take Taylor series of each term around $\theta=0$, except for the MGF of $\epsilon\langle\vc,\vqbar^\peps\rangle$. To do that, we use the following lemma.

\begin{lemma}\label{lemma.taylor}
	Let $X^\peps$ be a set of random variables indexed by $\epsilon>0$. Assume $X^\peps$ is bounded for all $\epsilon$, i.e., there exists a constant $K_{\max}$ (that does not depend on $\epsilon$) such that $X^\peps\leq K_{\max}$ with probability 1. Define $f_{\epsilon,X}(\theta)\defn e^{\theteps X^\peps}$. Then,
	\begin{align*}
	\left|\E{f_{\epsilon,X}(\theta)}-1-\theteps \E{X^\peps}-\dfrac{(\theteps)^2}{2}\E{\left(X^\peps\right)^2} \right| \leq C\epsilon^3,
	\end{align*}
	where $C$ is a finite constant. With a slight abuse of notation, we write the inequality above as follows
	\begin{align}\label{eq.taylor.order3}
	\E{f_{\epsilon,X}(\theta)}= 1+\theteps \E{X^\peps}+\dfrac{(\theteps)^2}{2}\E{\left(X^\peps\right)^2}+O(\epsilon^3).
	\end{align}
\end{lemma}

We present the proof of Lemma \ref{lemma.taylor} in Appendix \ref{app.taylor}.

\begin{remark}
	Since we are working with a bounded random variable, the proof that we presented of Lemma \ref{lemma.taylor} was straightforward. However, in general, one needs an assumption on the existence of the MGF.
\end{remark}

Expanding each term on the right hand side of \eqref{eq.MGFmethod.fraction} in Taylor series according to Lemma \ref{lemma.taylor} will yield terms related to the moments of the unused service. As illustrated in the case of the single server queue, it suffices to handle the first moment. To do that, we set to zero the drift of the linear test function $V_1(\vq)=\langle\vc,\vq\rangle$, i.e., we set $\E{\langle\vc,\left(\vqbar^\peps\right)^+\rangle}= \E{\langle\vc,\vqbar^\peps\rangle}$ (which is finite because in Step 1 it was proved that the MGF of $\epsilon\langle\vc,\vqbar^\peps\rangle$ exists in an interval around 0). For example, see Lemma \ref{lemma.ssq.Eu} in the case of the single server queue, which is used in \eqref{eq.ssq.taylor.u}.

From this step we obtain an expression for the limit as $\epsilon\dto 0$ of the MGF of $\epsilon\langle\vc,\vqbar^\peps\rangle$. This proves convergence in distribution of $\epsilon\langle\vc,\vqbar^\peps\rangle$ to a random variable $Y$, which turns out to be exponential in the cases we study in this paper. Then, $\epsilon\vqparbar^\peps=\epsilon\langle\vc,\vqbar^\peps\rangle\vc \Rightarrow Y\vc$ as $\epsilon\dto 0$ because $\vc$ is a fixed vector. We also know from SSC in Prerequisite 2 that $\epsilon\vqperpbar^\peps \to 0$ in probability as $\epsilon\dto 0$. Then, by Slutsky's theorem \cite[Theorem 11.4 in Section 5]{gut2012probability},  we obtain that $\epsilon\vqbar^\peps=\epsilon\vqparbar^\peps+\epsilon\vqperp^\peps\Rightarrow Y\vc$ as $\epsilon\dto 0$. 

\begin{remark}\label{remark.other.transforms}
	In order to set $\E{e^{\theteps \langle\vc,\left(\vqbar^\peps\right)^+\rangle}}=\E{e^{\theteps \langle\vc,\vqbar^\peps\rangle}}$ in Step 1, one must first prove the existence of  the MGF of $\epsilon\langle\vc,\vqbar^\peps\rangle$ in an interval around zero. An alternative approach (where this difficulty does not arise), is to use characteristic functions, because they always exist. However, working with characteristic functions involve the use of  complex analysis. Another way to overcome this difficulty is to use one-sided Laplace transform, i.e., to consider $\theta<0$. One-sided Laplace transform of $\epsilon\langle\vc,\vq\rangle$ always exists because $\epsilon$, $\vc$ and $\vq$ are nonnegative. If one chooses to work with other transforms such as the characteristic function or one-sided Laplace transform to get around the issue of the existence of the MGF, then one needs to assume that certain moments exist in a counterpart of Lemma \ref{lemma.taylor}. For instance, Theorem 2.3.3. in \citep{lukacs1970characteristic} can be used when one is working with characteristic functions.
\end{remark}

\section{Load balancing systems}\label{sec:load.balancing}
In this section we use the MGF method in the context of load balancing systems, also known as supermarket checkout systems. We first define the model and then we use the MGF method to prove that the steady-state distribution of the scaled vector of queue lengths is exponential in heavy-traffic.

\subsection{Load balancing model}\label{sec:load.balancing.model}

Consider a system with $n$ separate queues, as described in Section \ref{sec:MGF.general.model}. For each $i\in\{1,\ldots,n\}$, $\{s_i(k):\;k\geq 1\}$ is a sequence of i.i.d. random variables with $\mu_i\defn\E{s_i(1)}$, and let $\mu_{\Sigma}\defn\ds\sum_{i=1}^n \mu_i$. We consider this system in a general setting, so we do not assume independence of the servers. For $i,j\in\{1,\ldots,n\}$, let $\Cov{s_i,s_j}$ be the covariance between $s_i(1)$ and $s_j(1)$. There is a single stream of arrivals, that we model as a sequence $\{a(k):\;k\geq 1\}$ of i.i.d. random variables such that $a(k)$ is the number of arrivals to the system in time slot $k$. In this queueing system the control problem is to route the arrivals to one of the $n$ queues in each time slot. We assume the routing policy is fixed for all $k\geq 1$, but we do not assume any particular policy. After routing, $a_i(k)$ is the number of arrivals routed to the $i\tth$ queue in time slot $k$, for $i\in\{1,\ldots,n\}$. We assume $a(k)\leq\amax$ with probability 1 for all $k\geq 1$, and that the arrival process is independent of the queue length and service processes. The dynamics of the queues are according to \eqref{eq.dynamics}.
It is well known that the capacity region of the load balancing system is $\cC=\{\lambda\in\bR_+:\; \lambda\leq \mu_\Sigma \}$. A proof can be found in Appendix A of \citep{atilla}.

To study the heavy-traffic limit of this queueing system, we parametrize the arrival process as follows. For $\epsilon\in(0,\mu_\Sigma)$ we consider a load balancing system as described above, where the arrival process $\{a^\peps(k):\;k\geq 1 \}$ is such that $\E{a^{\peps}(1)}= \mu_{\Sigma}-\epsilon$ and $\Var{a^{\peps}(1)}=\left(\sigma_a^\peps \right)^2$. In other words, the arrival rate approaches the point $r=\mu_\Sigma$ in the boundary of $\cC$ as $\epsilon\dto 0$. Since the capacity region $\cC$ of the load balancing system is one-dimensional, the CRP condition (as defined in Definition \ref{def.CRP}) is trivially satisfied.

\subsection{MGF method applied to load balancing systems}\label{sec:load.balancing.theorem}

In this subsection we state the main theorem of this section and provide some examples, and in the next subsection we will prove the theorem using the MGF method as developed in Section \ref{sec:details.MGFmethod}. Before presenting the formal statement of the result we introduce the following definitions.

\begin{definition}[Throughput optimality] \label{def.load.balancing.throughput}
	A routing algorithm $\cA$ is throughput optimal for the load balancing system described in Section \ref{sec:load.balancing.model} if the Markov chain $\left\{\vq^\peps(k):\;k\geq 1 \right\}$ operating under $\cA$ is positive recurrent for all $\epsilon\in\left(0,\mu_\Sigma\right)$.
\end{definition}

\begin{definition}[State Space Collapse]\label{def.load.balancing.ssc}
	Consider a routing algorithm $\cA$ and let
	\begin{align*}
		\cK=\left\{\vx\in\bR^n:\; x_i=x_j\quad\forall i,j\in\{1,\ldots,n\} \right\},
	\end{align*}
	i.e., $\vc=\frac{1}{\sqrt{n}}\vone$. For any vector $\vy\in\bR^n$, let $\vy_{\parallel}$ be the projection of $\vy$ on $\cK$ and let $\vy_\perp\defn \vy-\vy_\parallel$. We say that the algorithm $\cA$ satisfies SSC if the load balancing system described in Section \ref{sec:load.balancing.model} operating under $\cA$ satisfies the following property.
	\begin{align*}
		\E{\left\|\vqperpbar^\peps \right\|^2}\;\text{is}\; o\left(\frac{1}{\epsilon^2}\right)
	\end{align*}
	where $\vqbar^\peps$ is a steady-state random vector such that  $\left\{\vq^\peps(k):\;k\geq 1\right\}$ converges in distribution to $\vqbar^\peps$ if it is positive recurrent.
\end{definition}

Observe that if an algorithm $\cA$ satisfies SSC (as defined above), then SSC occurs into the one-dimensional space $\cK$. Therefore, a load balancing system operating under such $\cA$ behaves as a single server queue in the heavy-traffic limit.

Now we formally present the result that we will prove using the MGF method.

\begin{theorem}\label{thm.load.balancing}
	Let $\epsilon\in(0,\mu_\Sigma)$ and consider a set of load balancing systems parametrized by $\epsilon$, as described in Section \ref{sec:load.balancing.model}. Suppose that the routing algorithm is throughput optimal and that it satisfies SSC. For each $\epsilon\in(0,\mu_{\Sigma})$, let $\vqbar^\peps$ be a steady-state random vector such that the queue length process $\{\vq^\peps(k):\;k\geq 1\}$ converges in distribution to $\vqbar^\peps$. Assume the MGF of $\epsilon\sum_{i=1}^n \qbar_i$ exists, i.e., $\E{e^{\theteps \sum_{i=1}^n \qbar_i^\peps}}<\infty$ for $\theta\in[-\Theta,\Theta]$ where $\Theta>0$ is a finite number, and that $\lim_{\epsilon\dto 0}\sigma_a^\peps=\sigma_a$. Then $\epsilon\vqbar^\peps\Longrightarrow \widetilde{\Upsilon} \vone$ as $\epsilon\dto 0$, where $\widetilde{\Upsilon}$ is an exponential random variable with mean $\dfrac{1}{2n}\left(\ds\sigma_a^2+\sum_{i=1}^n\sum_{j=1}^n \text{\emph{Cov}}\left[s_i,s_j\right]\right)$.
\end{theorem}

Now we introduce two routing policies that satisfy SSC as defined above. We first define the policies.

\begin{definition}[JSQ and Power-of-two choices]
	Consider a load balancing system as described in Section \ref{sec:load.balancing.model}. Then, for each $k\geq 1$, given the vector of queue lengths $\vq^\peps(k)$, a routing policy selects $i^*(k)$ and sends arrivals according to the following formula.
	\begin{align*}
		a^\peps_i(k)=\begin{cases}
			a^\peps(k) & \text{, if }i=i^*(k) \\
			0 & \text{, otherwise.}
		\end{cases}
	\end{align*}
	
	\begin{enumerate}[label=(\alph*)]
		\item\label{def.load.balancing.jsq} The routing policy Join the Shortest Queue (JSQ) sends all arrivals in time slot $k$ to the queue with the least number of jobs, breaking ties at random. Formally, under JSQ routing policy
		\begin{align*}
			i^*(k)\in \argmin_{i\in\{1,\ldots,n\}}\left\{q^\peps_i(k) \right\},
		\end{align*}
		breaking ties at random.
		
		\item\label{def.load.balancing.po2} The routing policy power-of-two choices selects two queues uniformly at random, say $i_1,i_2\in\{1,\ldots,n\}$ and sends all arrivals in time slot $k$ to the queue with the least number of jobs between those two, breaking ties at random. Formally, under power-of-two choices, if queues $i_1$ and $i_2$ are selected, then
		\begin{align*}
			i^*(k)\in\argmin_{i\in\{i_1,i_2\}}\left\{q^\peps_i(k) \right\},
		\end{align*}
		breaking ties at random.
	\end{enumerate}
\end{definition}

In the following two corollaries we show that these routing policies satisfy the assumptions of Theorem \ref{thm.load.balancing} and, therefore, the scaled vector of queue lengths in a load balancing system operating under any of these policies has an exponential distribution in heavy-traffic.

\begin{corollary}\label{cor.load.balancing.jsq}
	Consider a set of load balancing systems parametrized by $\epsilon\in(0,\mu_\Sigma)$ as described in Section \ref{sec:load.balancing.model}, operating under JSQ routing policy. Then, $\epsilon\vqbar^\peps\Longrightarrow \widetilde{\Upsilon}_1 \vone$ as $\epsilon\dto 0$, where $\widetilde{\Upsilon}_1$ is an exponential random variable with mean $\dfrac{1}{2n}\left(\ds\sigma_a^2+\sum_{i=1}^n\sum_{j=1}^n \text{\emph{Cov}}\left[s_i,s_j\right]\right)$.
\end{corollary}

A particular case of the queueing system described in Corollary \ref{cor.load.balancing.jsq} is the load balancing system operating under JSQ with independent servers. In this case, $\sum_{i=1}^n\sum_{j=1}^n \Cov{s_i,s_j}$ reduces to the sum of variances of the servers. This is one of the systems studied by \cite{atilla}.

\proof{Proof of Corollary \ref{cor.load.balancing.jsq}.}
	We only need to show that JSQ is throughput optimal, that it satisfies SSC, and that there exists $\Theta>0$ such that $\E{e^{\theteps\sum_{i=1}^n \qbar_i^\peps}}<\infty$ for all $\theta\in[-\Theta,\Theta]$. \cite{atilla} prove throughput optimality and SSC in the case of independent servers. However, their proofs hold for correlated servers. The proof of throughput optimality can be found in Appendix A of \cite{atilla}. 
	
	The SSC result proved by \cite{atilla} is stronger than the property presented in Definition \ref{def.load.balancing.ssc}. In fact, they prove that $\E{\left\|\vqperpbar^\peps \right\|^m}$ is upper bounded by a constant for each $m=1,2,\ldots$. This clearly implies that Definition \ref{def.load.balancing.ssc} is satisfied. We provide a sketch of their proof of SSC in Appendix \ref{app.jsq.ssc}.
	
	The existence of MGF of $\epsilon\sum_{i=1}^n \qbar_i^\peps$ in an interval around 0 is proved in Appendix \ref{app.jsq.mgf}. \Halmos
\endproof

\begin{corollary}\label{cor.load.balancing.po2}
	Consider a set of load balancing systems parametrized by $\epsilon$ as described in Section \ref{sec:load.balancing.model}, operating under Power-of-two choices and where all the servers are identical. Then, $\epsilon\vqbar^\peps\Longrightarrow \widetilde{\Upsilon}_2 \vone$ as $\epsilon\dto 0$, where $\widetilde{\Upsilon}_2$ is an exponential random variable with mean $\dfrac{1}{2n}\left(\ds\sigma_a^2+\sum_{i=1}^n\sum_{j=1}^n \text{\emph{Cov}}\left[s_i,s_j\right]\right)$.
\end{corollary}

\proof{Proof of Corollary \ref{cor.load.balancing.po2}.}
	Similar to the proof of Corollary \ref{cor.load.balancing.jsq}, we check throughput optimality, SSC and existence of MGF. \cite{magsriyin_itc12_journal} prove SSC in the case of independent servers in Section 4.3 of the article, but their proof holds true if this assumption is dropped. Their proof is along the lines of the proof for JSQ in Appendix \ref{app.jsq.ssc}, so we do not present it here. Throughput optimality can be proved using Foster-Lyapunov theorem and the calculations that \cite{magsriyin_itc12_journal} develop in the proof of SSC, and existence of MGF is similar to the case of JSQ. We omit these proofs in this paper, since our goal is to introduce the MGF method. \Halmos
\endproof

Observe that the assumption of identical servers is essential for the power-of-two choices algorithm to be throughput optimal. The case when the servers are not identical was studied by 
\cite{chen2012po2_journal} using the diffusion limits approach. The routing policy there randomly selects $d$ servers in each time slot, where the probability of choosing server $i$ is proportional to its service rate $\mu_i$, for all $i\in\{1,\ldots,n\}$. Then, the arrivals are sent to the server with the shortest queue among the $d$ selected servers. They prove that this queueing system satisfies the CRP condition and that the distribution of the scaled vector of queue lengths is exponential.
A similar result can be obtained using the MGF method once the SSC as stated in Definition \ref{def.load.balancing.ssc} is established. This is straightforward extension, and we do not present the details here because the focus is on illustrating the MGF approach.

In this subsection we presented the main theorem of this section, and two examples where the assumptions of the theorem are satisfied. Observe that in both cases we only needed to check that the conditions of the theorem are satisfied. In fact, if we want to prove that the scaled vector of queue lengths of the load balancing system operating under any other routing policy has an exponential distribution, we only need to check these three assumptions.

\subsection{Proof of Theorem \ref{thm.load.balancing}}\label{sec:load.balancing.proof}

In the rest of this section we prove Theorem \ref{thm.load.balancing} using the MGF method. Before presenting the proof we specify notation.

Let $\abar^\peps$ be a steady-state random variable with the same distribution as $a^\peps(1)$ and let $\vabar^\peps\defn \va^\peps(\vqbar)$ be the vector of arrivals to each queue after routing in steady-state. The vector $\vubar^\peps$ is defined as in Section \ref{sec:MGF.general.model}. Observe that in this case the vector $\vsbar$ is independent of $\vqbar^\peps$ and it has the same distribution as $\vs(1)$, because the potential service sequences $\{s_i(k):k\geq 1\}$ are i.i.d. and independent of the queue length processes for each $i\in\{1,\ldots,n\}$.

\proof{Proof of Theorem \ref{thm.load.balancing}.}

For ease of exposition, we omit the dependence on $\epsilon$ of the variables in this proof. We use the MGF method. Before applying the steps, we need to verify that the prerequisites are satisfied, i.e., we need to check positive recurrence and SSC. In fact, one of the assumptions of the theorem is that the routing policy is throughput optimal. Therefore, for any $\epsilon>0$ the Markov chain $\left\{\vq^\peps(k):\,k\geq 1 \right\}$ is positive recurrent. Also, SSC is satisfied by assumption. Now we go through the steps of the MGF method.

\subsubsection*{Step 1. Prove an equation of the form of \eqref{eq.MGF.method.exp.qu} and compute an expression for the MGF of $\epsilon\langle\vc,\vqbar^\peps\rangle$.}

We first prove the following lemma.
	
\begin{lemma}\label{lemma.load.balancing.expoqu}
	Consider a load balancing system parametrized by $\epsilon$ as described in Theorem \ref{thm.load.balancing}. Then, there exists $\theta_{\max}>0$ finite such that for any real number $\theta\in[-\theta_{\max},\theta_{\max}]$ we have
	\begin{align*}
	\E{\left(e^{\theteps \sum_{i=1}^n \left(\qbar_i^\peps\right)^\plus}-1 \right)\left(e^{-\theteps \sum_{i=1}^n \ubar_i^\peps}-1\right)}\quad\text{is $o(\epsilon^2)$}
	\end{align*}
\end{lemma}

We present the proof of Lemma \ref{lemma.load.balancing.expoqu} in Appendix \ref{app.load.balancing.expo.qu}.

Since $\ds\langle\vc,\vq\rangle=\dfrac{1}{\sqrt{n}}\sum_{i=1}^n q_i$, proving an equation of the form of \eqref{eq.MGF.method.exp.qu} is equivalent to Lemma \ref{lemma.load.balancing.expoqu} using $\frac{\theta}{\sqrt{n}}$ instead of $\theta$. For ease of exposition, we work with $\theta$ in the rest of this proof.

Note that  $\Prob{\abar-\sum_{i=1}^n \sbar_i\neq 0}>0$ whenever $\epsilon>0$. If we expand the product in the expression of Lemma \ref{lemma.load.balancing.expoqu} and we follow the steps sketched after Step 1 in Section \ref{sec:details.MGFmethod} we obtain
\begin{align}\label{eq.load.balancing.before.fraction}
& \E{e^{\theteps\sum_{i=1}^n\qbar_i }\left(1-e^{\theteps \sum_{i=1}^n(\abar_i-\sbar_i)}\right)} = 1-\E{e^{-\theteps\sum_{i=1}^n \ubar_i }}+o(\epsilon^2).
\end{align}

Recall $\sum_{i=1}^n \abar_i=\abar$ and that $\abar$, $\vsbar$ are independent of $\vqbar$, by definition. Therefore, reorganizing terms we obtain
\begin{align}\label{eq.load.balancing.fraction}
\E{e^{\theteps\sum_{i=1}^n \qbar_i}}= \dfrac{1-\E{e^{-\theteps\sum_{i=1}^n \ubar_i }} +o(\epsilon^2)}{1-\E{e^{\theteps\left(\abar-\sum_{i=1}^n \sbar_i\right)} } },
\end{align}
which gives an expression for the MGF of $\epsilon\sum_{i=1}^n\qbar_i$ that is valid for all traffic.

\subsubsection*{Step 2. Bound unused service and take heavy-traffic limit.} 

Equation \eqref{eq.load.balancing.fraction} yields a $\frac{0}{0}$ form in the limit as $\epsilon\dto 0$, just like \eqref{eq.ssq.fraction} in the case of the single server queue. Equivalently, we can observe that \eqref{eq.load.balancing.before.fraction} yields $0=0$ in the limit as $\epsilon\dto 0$. Then, we take Taylor series of the numerator and the denominator of \eqref{eq.load.balancing.fraction} at $\theta=0$ to obtain the limit. To take Taylor expansion we use Lemma \ref{lemma.taylor}.

In order to bound the numerator we need to compute $\E{\sum_{i=1}^n \ubar_i}$, so we start with a lemma.

\begin{lemma}\label{lemma.load.balancing.Eu}
	Consider a load balancing system parametrized by $\epsilon\in(0,\mu_\Sigma)$ as described in Section \ref{sec:load.balancing.model}, operating under a throughput optimal routing policy. Then,
	\begin{align*}
	\E{\sum_{i=1}^n \ubar_i^\peps}=\epsilon.
	\end{align*}
\end{lemma}

\proof{Proof of Lemma \ref{lemma.load.balancing.Eu}.}

We set to zero the drift of $V_1(\vq)=\langle\vc,\vq\rangle$ in steady-state. In this case, from the definition of $\cK$ in Definition \ref{def.load.balancing.ssc} we have $\vc=\dfrac{1}{\sqrt{n}}\vone$. Then, we obtain
\begin{align*}
0=& \E{V_1\left(\vqbar^\plus\right)-V_1\left(\vqbar\right)} \\
=& \dfrac{1}{\sqrt{n}} \E{\sum_{i=1}^n \qbar_i^\plus - \sum_{i=1}^n \qbar_i} \\
\stackrel{(a)}{=}& \dfrac{1}{\sqrt{n}}\E{\sum_{i=1}^n \left(\qbar_i+\abar_i-\sbar_i+\ubar_i\right)-\sum_{i=1}^n \qbar_i} \\
\stackrel{(b)}{=}& \dfrac{1}{\sqrt{n}}\E{\abar-\sum_{i=1}^n \sbar_i + \sum_{i=1}^n \ubar_i}
\end{align*}
where $(a)$ holds by definition of $\vqbar^\plus$; and $(b)$ holds because $\abar=\sum_{i=1}^n \abar_i$ by definition of $\abar$ and $\abar_i$. Rearranging terms and canceling $\frac{1}{\sqrt{n}}$, we obtain
\begin{align*}
\E{\sum_{i=1}^n \ubar_i}=& \sum_{i=1}^n \E{\sbar_i}-\E{\abar} \\
\stackrel{(a)}{=}& \sum_{i=1}^n \mu_i-(\mu_{\Sigma}-\epsilon) \\
\stackrel{(b)}{=}& \epsilon,
\end{align*}
where $(a)$ holds because $\E{\abar}=\mu_\Sigma-\epsilon$; and $(b)$ holds by definition of $\mu_\Sigma$. \Halmos
\endproof 

Now we expand the numerator and denominator of \eqref{eq.load.balancing.fraction} in Taylor series. We start with the numerator, and we obtain
\begin{align}
1-\E{e^{-\theteps\sum_{i=1}^n \ubar_i }}=& 1-\E{f_{\epsilon,-\sum_{i=1}^n \ubar_i}(\theta)} \nonumber\\
=& \theteps\E{\sum_{i=1}^n \ubar_i}-\dfrac{(\theteps)^2}{2}\E{\left(\sum_{i=1}^n \ubar_i\right)^2} +O(\epsilon^3) \nonumber \\
=& \theteps^2 -\dfrac{(\theteps)^2}{2}\E{\left(\sum_{i=1}^n \ubar_i\right)^2} +O(\epsilon^3), \label{eq.load.balancing.num.partial}
\end{align}
where the last equality holds by Lemma \ref{lemma.load.balancing.Eu}. Now we need to bound the second moment of the sum of unused services.

\begin{claim}\label{claim.load.balancing}
	Consider a load balancing system as described in Theorem \ref{thm.load.balancing}. Then,
	\begin{align*}
	\dfrac{(\theteps)^2}{2}\E{\left(\sum_{i=1}^n \ubar_i^\peps\right)^2}\; \text{is }O(\epsilon^3).
	\end{align*}
\end{claim}

We prove the claim in Appendix \ref{app.claim.load.balancing}. Using the Claim in \eqref{eq.load.balancing.num.partial} we obtain
\begin{align}
1-\E{e^{-\theteps\sum_{i=1}^n \ubar_i }}=& \theteps^2+O(\epsilon^3), \label{eq.load.balancing.numerator}
\end{align}

For the denominator, we obtain
\begin{align}
1-\E{e^{\theteps\left(\abar-\sum_{i=1}^n \sbar_i\right)}} =& 1-\E{f_{\epsilon,\left(\abar-\sum_{i=1}^n \sbar_i\right)}(\theta)} \nonumber \\
=& -\theteps\E{\abar-\sum_{i=1}^n \sbar_i}-\dfrac{(\theteps)^2}{2}\E{\left(\abar-\sum_{i=1}^n \sbar_i\right)^2}+O(\epsilon^3) \nonumber \\
=& \theteps^2-\dfrac{(\theteps)^2}{2}\left(\left(\sigma_a^\peps\right)^2+\sum_{i=1}^n \sum_{j=1}^n \Cov{s_i,s_j} +\epsilon^2\right)+O(\epsilon^3), \label{eq.load.balancing.denominator}
\end{align}
where the last step holds because $\E{\abar}=\mu_\Sigma-\epsilon$, $\E{\sum_{i=1}^n \sbar_i}=\mu_{\Sigma}$ and by definition of covariance.

Using \eqref{eq.load.balancing.numerator} and \eqref{eq.load.balancing.denominator} in \eqref{eq.load.balancing.fraction}, and since $O(\epsilon^3)$ is $o(\epsilon^2)$, we obtain
\begin{align*}
\E{e^{\theteps\sum_{i=1}^n \qbar_i}}= \dfrac{\theteps^2+o(\epsilon^2)}{\ds\theteps^2-\dfrac{(\theteps)^2}{2}\left(\left(\sigma_a^\peps\right)^2+\sum_{i=1}^n \sum_{j=1}^n \Cov{s_i,s_j} +\epsilon^2\right)+O(\epsilon^3)}
\end{align*}

Canceling $\theteps^2$ from the numerator and denominator, we obtain
\begin{align*}
\E{e^{\theteps\sum_{i=1}^n \qbar_i}}=\dfrac{1+o(1)}{1-\dfrac{\theta}{2}\left(\left(\ds\sigma_a^\peps\right)^2+\sum_{i=1}^n\sum_{j=1}^n \Cov{ s_i,s_j} \right)+O(\epsilon)}.
\end{align*}

Therefore, taking the limit we obtain
\begin{align*}
\lim_{\epsilon\dto 0}\E{e^{\theteps\sum_{i=1}^n \qbar_i}}=\dfrac{1}{1-\frac{\theta}{2}\left(\sigma_a^2+\sum_{i=1}^n \sum_{j=1}^n \Cov{s_i,s_j}\right)},
\end{align*}
which is the MGF of an exponential random variable with mean $\dfrac{1}{2}\left(\sigma_a^2+\sum_{i=1}^n \sum_{j=1}^n \Cov{s_i,s_j}\right)$. Then, $\epsilon\langle\vc,\vqbar\rangle\vc = \epsilon\left( \frac{1}{n}\sum_{i=1}^n \qbar_i\right)\vone\Rightarrow \widetilde{\Upsilon} \vone$ as $\epsilon\dto 0$, where $\widetilde{\Upsilon}$ is an exponential random variable with mean $\frac{1}{2n}\left(\sigma_a^2+\sum_{i=1}^n\sum_{j=1}^n \Cov{s_i,s_j}\right)$.

Therefore, we conclude that $\epsilon\vqbar^\peps=\epsilon\vqparbar^\peps+\epsilon\vqperpbar^\peps \Rightarrow \widetilde{\Upsilon}\vone$ as $\epsilon\dto 0$. This proves Theorem \ref{thm.load.balancing}. \Halmos
\endproof

\section{Generalized switch}\label{sec:generalized.switch}
In this section we apply the MGF method in the context of a generalized switch operating under MaxWeight. We compute the distribution of the scaled vector of queue lengths in heavy-traffic under the assumption that CRP is satisfied. The generalized switch is a model that was first introduced by \cite{stolyar2004maxweight}, and it represents a generalization of a variety of queueing systems, such as the input-queued switch \citep{mckeown96walrand}, cloud computing \citep{magsriyin_itc12_journal},  down-links in wireless base stations \citep{TasEph_92}, etc.

\subsection{Generalized switch model}\label{sec:generalized.switch.model}

Consider a system with $n$ separate queues, as described in Section \ref{sec:MGF.general.model}. For each $i\in\{1,\ldots,n\}$, let $\{a_i(k):\;k\geq 1\}$ be a sequence of i.i.d. random variables such that $a_i(k)$ is the number of arrivals to queue $i$ in time slot $k$. For $i,j\in\{1,\ldots,n\}$, let $\Cov{a_i,a_j}$ be the covariance between $a_i(1)$ and $a_j(1)$. The servers interfere with each other. Then, the vector of service rates must satisfy feasibility constraints in each time slot. Additionally, there are conditions of the environment that affect these constraints. We group all these conditions in a random variable called channel state. For each $k\geq 1$, let $T(k)$ be the channel state in time slot $k$. The sequence of random variables $\{T(k):\;k\geq 1\}$ is i.i.d. and it is independent of the queue length and the arrival processes. We assume that the state space of the channel state is a finite set $\cT$ and we let $\boldsymbol{\psi}$ be the probability mass function of $T(1)$, i.e., for each $t\in \cT$ the probability of observing state $t$ is $\psi_t\defn \Prob{T(1)=t}$. For each $t\in\cT$, let $\cS^\pt$ be the set of feasible service rate vectors under channel state $t$. We also assume that if $\vx\in\cS^\pt$ for some $t\in\cT$, then all vectors that are strictly dominated by $\vx$ are feasible. In other words, if $\vy$ is a nonnegative vector that satisfies $\vy\leq \vx$ component-wise, then $\vy$ is also a feasible service rate vector if the channel state is $t$. In particular, the projection of $\vx\in\cS^\pt$ on each of the coordinate axes is a feasible service rate vector as well. We assume that $\cS^\pt$ is finite for each $t\in\cT$, so we only consider maximal feasible schedules and their projection on the coordinate axes in $\cS^\pt$. With this assumption we do not lose much generality because the vector $\vs(k)$ is the potential (not actual) service rate vector and we are interested in the heavy-traffic limit.

In this queueing system the control problem (which is a scheduling problem), is to select $\vs(k)$ in each time slot after realizing the channel state. Let $\vs(k)$ be the solution of the scheduling problem in time slot $k$. Since $\cS^\pt$ is finite for each $t\in\cT$ and $\cT$ is also finite, there exists a constant $\smax$ such that $s_i(k)\leq \smax$ for all $i\in\{1,\ldots,n\}$ and all $k\geq 1$. 

It is known \citep{atilla} that the capacity region of this queueing system is
\begin{align}\label{eq.gen.switch.capacity.region}
	\cC=\sum_{t\in\cT}\psi_t\, \text{ConvexHull}\left\{\cS^\pt\right\}.
\end{align}
Providing a formal proof of \eqref{eq.gen.switch.capacity.region} is beyond the scope of this paper, but we intuitively explain why it holds. First suppose that the channel state is fixed and the set of feasible service rate vectors is $\cS^{(1)}$. Then, the capacity region should have all vectors $\vx$ that satisfy $\vx\leq \vs$ for all $\vs\in \cS^{(1)}$. Since $\cS^{(1)}$ contains the projection of its elements on the coordinate axis, the set of such vectors $\vx$ is $\text{ConvexHull}\left\{\cS^{(1)} \right\}$. Now, if we consider the channel state as a random variable, recall that $\psi_t$ is the probability that the channel state is $t$, and if the channel state is $t$ then the set of feasible service rate vectors is $\cS^\pt$. Then, \eqref{eq.gen.switch.capacity.region} just gives the capacity region associated to each channel state, weighted by the probability that each channel state is observed.

Recall that, by assumption, each set $\cS^\pt$ is finite. Then, for each $t\in\cT$ the set $\text{ConvexHull}\left\{\cS^\pt \right\}$ is the convex hull of finitely many points. Therefore, $\text{ConvexHull}\left\{\cS^\pt \right\}$ is a polytope, i.e., a bounded polyhedron. Also, the state space of the channel state $\cT$ is finite by assumption. Then, \eqref{eq.gen.switch.capacity.region} is the weighted sum of finitely many polytopes. This implies that $\cC$ is also a polytope. In order to exploit this structure, we describe it as the intersection of a finite number of half-spaces, where each half-space defines a facet of $\cC$. Let $L$ be the minimal number of half-spaces that are required to describe $\cC$, and for each $\ell\in\{1,\ldots,L\}$ let $\vcl\in\bR^n$ and $\bl\in\bR$ be the parameters that define each facet of the polytope. In other words, we describe $\cC$ as follows
\begin{align}\label{eq.gen.switch.capacity.region2}
	\cC= \left\{\vx\in\bR^n:\;\langle \vcl,\vx\rangle\leq \bl\;\text{ for all }\ell\in\{1,\ldots,L\} \right\}.
\end{align}

Without loss of generality we can assume $\vcl\geq \vzero$, $\|\vcl\|=1$ and $\bl>0$ for all $\ell\in\{1,\ldots,L\}$, because we assumed that the sets $\cS^\pt$ contain the projection on the coordinate axes of all their feasible vectors. Therefore, the capacity region is coordinate convex. For each $\ell\in\{1,\ldots,L\}$, let $\cFl\defn \left\{\vx\in\cC:\;\langle\vcl,\vx\rangle=\bl \right\}$ be the $\ell\tth$ facet of the polytope $\cC$.

In this paper we assume that the scheduling problem is solved using MaxWeight algorithm in each time slot, i.e., if the channel state is $t$, then the selected schedule satisfies
\begin{align}\label{eq.gs.MW}
\vs(k)\in \argmax_{\vx\in\cS^\pt}\langle\vx,\vq(k)\rangle
\end{align} 
and ties are broken at random.

From \eqref{eq.gen.switch.capacity.region} and \eqref{eq.gs.MW}, observe that the service rate vector $\vs(k)$ does not necessarily belong to the capacity region $\cC$ because $\psi_t\leq 1$ for all $t\in\cT$. To overcome this difficulty we define the following random variable.
For each $\ell\in\{1,\ldots,L\}$ and each $t\in\cT$, define the \textit{maximum $\vcl$-weighted service rate available in channel state $t$} \citep{atilla} as
\begin{align}\label{eq.gen.switch.def.bjl}
	\btl= \max_{\vs\in\cS^\pt}\langle \vcl,\vs\rangle.
\end{align}
In other words, given that the channel state is $t$, $\btl$ is a real number such that the hyperplane $\mathcal{H}^{(t,\ell)}=\left\{\vx\in\bR^n:\;\langle \vcl,\vx\rangle=\btl\right\}$ is tangent to the boundary of $\text{ConvexHull}\left\{\cS^\pt\right\}$. Let $\{B_\ell(k):\;k\geq 1\}$ be a sequence of i.i.d. random variables such that $\Prob{B_\ell(k)=\btl}=\psi_t$ and $\sigma_{B_\ell}^2\defn\Var{B_\ell(k)}$. In the next lemma we present the relation between the random variable $B_\ell(1)$ and the parameter $\bl$ for each $\ell\in\{1,\ldots,L\}$.

\begin{lemma}\label{gen.switch.lemma.B.bl}
	Consider a generalized switch as described above. Then, for each $\ell\in\{1,\ldots,L\}$
	\begin{align*}
		\E{B_\ell(1)}=\bl.
	\end{align*}
\end{lemma}
\proof{Proof of Lemma \ref{gen.switch.lemma.B.bl}.}

By definition of the random variable $B_\ell(1)$, we have
\begin{align*}
\E{B_\ell(1)}=& \sum_{t\in\cT} \psi_t\btl \\
\stackrel{(a)}{=}& \sum_{t\in\cT} \psi_t \max_{\vs\in\cS^\pt}\langle\vcl,\vs\rangle \\
\stackrel{(b)}{=}& \max_{\vs\in \cC}\langle\vcl,\vs\rangle \\
\stackrel{(c)}{=}& \bl
\end{align*}
where $(a)$ holds by the definition of $\btl$ given in \eqref{eq.gen.switch.def.bjl}; $(b)$ holds by definition of the capacity region $\cC$ given in \eqref{eq.gen.switch.capacity.region}; and $(c)$ holds by definition of the $\ell\tth$ facet and because the objective function in the maximization problem is linear. \Halmos
\endproof 

To perform heavy-traffic analysis, we fix a facet $\cFl$ and we study a set of generalized switches where the vector of arrival rates approaches a fixed point in the relative interior of $\cFl$. Formally, we fix $\vrl$ in the relative interior of $\cFl$ and we let $\epsilon\in(0,1)$. Then, the system parametrized by $\epsilon$ is such that $\E{\va^\peps(k)}=\vrl-\epsilon\vcl$ and $\Cov{a_i^\peps,a_j^\peps}$ is the covariance between $a_i^\peps(1)$ and $a_j^\peps(1)$ for each $i,j\in\{1,\ldots,n\}$. In this case, since the point $\vr=\vrl$ of the boundary of the capacity region $\cC$ is in the relative interior of the facet $\cFl=\left\{\vx\in\cC: \langle\vc^\pl,\vx\rangle=b^\pl \right\}$, the unique outer normal vector to the capacity region $\cC$ at $\vr$ is the outer normal vector to the facet $\cFl$, i.e., it is $\vc^\pl$. Therefore, the CRP condition as defined in Definition \ref{def.CRP} is satisfied. Observe that if $\vr$ is in the intersection of two (or more) facets, then the CRP condition is not satisfied because there is a range of vectors that are normal to $\cC$ at $\vr$.

\subsection{MGF method applied to generalized switches}\label{sec:generalized.switch.theorem}

In this subsection we state the main theorem of this section and we provide some examples. In the next subsection we prove the theorem. 

\begin{theorem}\label{thm.generalized.switch}
	Let $\epsilon\in(0,1)$. Given the $\ell\tth$ facet of $\cC$, $\cFl$, and a vector $\vrl$ in the relative interior of $\cFl$, consider a set of generalized switches operating under MaxWeight algorithm, parametrized by $\epsilon$ as described in Section \ref{sec:generalized.switch.model}. For each $\epsilon$, let $\vqbar^\peps$ be a steady-state vector such that the queue length process $\left\{\vq^\peps(k):
	\; k\geq 1 \right\}$ converges in distribution to $\vqbar^\peps$. Further, 	let $\lim_{\epsilon\dto0} \text{\emph{Cov}}\left[a^\peps_i,a^\peps_j\right]=\text{\emph{Cov}}\left[a_i,a_j\right]$ for each $i,j\in\{1,\ldots,n\}$. Then, $\epsilon\vqbar^\peps\Rightarrow \overline{\Upsilon}\vcl$ as $\epsilon\dto 0$, where $\overline{\Upsilon}$ is an exponential random variable with mean $\dfrac{1}{2}\left(\ds\sum_{i=1}^n \sum_{j=1}^n c^\pl_i c^\pl_j \text{\emph{Cov}}\left[a_i,a_j\right]+ \sigma_{B_\ell}^2\right)$, where $c^\pl_i$ is the $i\tth$ element of $\vcl$, for each $i\in\{1,\ldots,n\}$.
\end{theorem}

In the next corollary we present a particular example of a generalized switch operating under MaxWeight.

\begin{corollary}\label{cor.gen.switch.adhoc}
	Consider a set of generalized switches parametrized by $\epsilon$, as described in Section \ref{sec:generalized.switch.model}, operating under MaxWeight algorithm. Suppose that $\cT$ has only one element, i.e. the channel state is fixed over time. Then, $\epsilon\vqbar^\peps\Rightarrow \overline{\Upsilon}_2\vcl$, where $\overline{\Upsilon}_2$ is an exponential random variable with mean $\dfrac{1}{2}\left(\ds\sum_{i=1}^n \sum_{j=1}^n c^\pl_i c^\pl_j \text{\emph{Cov}}\left[a_i,a_j\right]\right)$.
\end{corollary}

The queueing system described in Corollary \ref{cor.gen.switch.adhoc} is also known as ad hoc wireless network. In an ad hoc wireless network we have $\sigma_{B_\ell}^2=0$ because the channel state is not a random variable anymore. The input-queued switch or a cross bar switch \citep{srikantleibook,MagSri_SSY16_Switch,QUESTA_switch} is yet another system that is well studied. When only one port of the switch is saturated, it satisfies the CRP condition \citep{stolyar2004maxweight}, and forms a special case of Corollary \ref{cor.gen.switch.adhoc}. In the next subsection we present the model and we formalize this result.

\subsection{MGF method applied to the input-queued switch}

An input-queued switch is a generalized switch where $n$ is a perfect square, i.e., there exists an integer $N$ such that $n=N^2$. Then, it can be represented as a square matrix, where the rows are input ports and the columns are output ports. The feasibility constraints are that, in each time slot, at most one queue can be served from each input and output port, and all jobs take exactly one time slot to be processed. Therefore, the set of feasible service rate vectors is analogous to permutation matrices of $N\times N$. 
	
For each $i\in\{1,\ldots,N\}$ let $\vchi{i}$ be the normalized indicator vector of row $i$, i.e., it is such that for each $i'\in\{1,\ldots,n\}$ we have $\chii{i}_{i'}=\frac{1}{\sqrt{N}}$ if queue $i'$ corresponds to row $i$ of the switch and $\chii{i}_{i'}=0$ otherwise. Similarly, for each $j\in\{1,\ldots,N\}$ let $\vchitilde{j}$ be the normalized indicator vector of column $j$. With this notation, we can write the capacity region of the input-queued switch as
\begin{align*}
	\cC_{\text{switch}}\defn \left\{\vx\in\bR_+^n:\; \langle \vchi{i},\vx\rangle\leq 1 ,\; \langle\vchitilde{j},\vx\rangle\leq 1\;,\; \forall i,j\in\{1,\ldots,N\} \right\},
\end{align*}
which is the intersection of $L=2N$ half-spaces.

Only one port can be saturated in heavy-traffic to ensure that CRP condition is satisfied. Without loss of generality, assume input port 1 is saturated, i.e., we consider a vector $\vr^{(1)}\in\cF^{(1)}$, where $\cF^{(1)}\defn \left\{\vx\in\cC_{\text{switch}}: \langle\vchi{1},\vx\rangle=1 \right\}$. For simplicity, we let $\vr^{(1)}=\vchi{1}$. Then, the heavy-traffic parametrization for $\epsilon\in(0,1)$ is such that $\vlambda^\peps= (1-\epsilon)\vchi{1}$. Unlike the generalized switch, for the input-queued switch we do not give the scheduling algorithm. Instead, we write the result in terms of the conditions that this algorithm must satisfy (similar to the load balancing case).

Similar to the case of the load balancing system, we say that an algorithm $\cA$ is throughput optimal for the input-queued switch if $\{\vq^\peps(k):k\geq 1\}$ is positive recurrent for all $\epsilon\in(0,1)$. Also, defining $\vx_\parallel\defn \langle\vchi{1},\vx\rangle\vchi{1}$ and $\vx_\perp\defn \vx-\vx_\parallel$ for any vector $\vx$, we say that the switch operating under a scheduling algorithm $\cA$ satisfies SSC if 
\begin{align*}
	\E{\left\|\vqbar_{\perp}^\peps \right\|^2}\text{ is }o\left(\frac{1}{\epsilon^2}\right)
\end{align*}

In the next proposition we compute the distribution of the scaled vector of queue lengths in heavy-traffic.

\begin{proposition}\label{gs.prop:switch}
	Let $\epsilon\in(0,1)$ and consider a set of input-queued switches parametrized by $\epsilon$, as described above. Suppose that the scheduling algorithm is throughput optimal and it satisfies SSC. For each $\epsilon\in(0,1)$, let $\vqbar^\peps$ be a steady-state random vector such that the queue length process $\left\{\vq^\peps(k):k\geq 1 \right\}$ converges in distribution to $\vqbar^\peps$. Assume the MGF of $\epsilon\langle\vchi{1},\vqbar^\peps\rangle$ exists, and that $\lim_{\epsilon\dto 0}\Sigma_a^\peps=\Sigma_a$ component-wise. Then, $\epsilon\vqbar^\peps\Longrightarrow \overline{\Upsilon}_s \vchi{1}$ as $\epsilon\dto 0$, where $\overline{\Upsilon}_s$ is an exponential random variable with mean $\frac{1}{2}\sum_{i=1}^n\sum_{j=1}^n \chii{1}_i\chii{1}_j \Cov{a_i,a_j}$.
\end{proposition}
\proof{Sketch of proof of Proposition \ref{gs.prop:switch}.}
	For ease of exposition we do not write the dependence on $\epsilon$ of the variables. We use the MGF method. We only present a sketch of this proof, since it is similar to the proofs of Theorems \ref{thm.load.balancing} and \ref{thm.generalized.switch}. We only show the main differences. 
	
	Both prerequisites are satisfied by assumption. Now we go through the steps.
	\subsubsection*{Step 1. Prove an equation of the form of \eqref{eq.MGF.method.exp.qu} and compute an expression for the MGF of $\epsilon\langle\vchi{1},\vqbar^\peps\rangle$.} Proving an equation of the form of \eqref{eq.MGF.method.exp.qu} is similar to the proof of Lemmas \ref{lemma.load.balancing.expoqu} and \ref{lemma.gen.switch.expo.qu}. Then, following the steps sketched in Step 1 in Section \ref{sec:details.MGFmethod} we obtain 
	\begin{align*}
		\E{e^{\theteps\langle\vchi{1},\vqbar\rangle}\left(1- e^{\theteps\langle\vchi{1}, \vabar-\vsbar\rangle}\right)} = 1- \E{e^{-\theteps\langle\vchi{1},\vubar\rangle}} + o(\epsilon^2).
	\end{align*}
	Since $\vsbar$ is a function of the queue lengths that is obtained through the scheduling problem, $\vsbar$ is not independent of $\vqbar$. However, $\langle\vchi{1},\vsbar\rangle=\frac{1}{\sqrt{N}}$ because all the feasible schedules $\vsbar$ are analogous to permutation matrices. Then, the sum of all the elements of $\vsbar$ corresponding to the first input port (row 1 of the switch) is 1. Then, $\langle\vchi{(1)},\vsbar\rangle$ is independent of the vector of queue lengths $\vqbar$. Also, the vector of arrivals is independent of $\vqbar$. Therefore, reorganizing terms we obtain
	\begin{align*}
		\E{e^{\theteps \langle\vchi{1},\vqbar\rangle}}=\dfrac{1-\E{e^{-\theteps\langle \vchi{1},\vubar\rangle}} + o\left(\epsilon^2\right)}{1-\E{e^{\theteps \langle\vchi{1},\vabar-\vsbar\rangle}}}.
	\end{align*}
	
	\subsubsection*{Step 2. Bound unused service and take heavy-traffic limit.} This step is equivalent to Step 2 in the proof of Theorems \ref{thm.load.balancing} and \ref{thm.generalized.switch}, so we omit the details. \Halmos
\endproof

In the case of a generalized switch, one of the difficulties is to handle the dependence on the queue lengths of the potential service vector. In the case of an input-queued switch this difficulty does not arise because, even though $\vsbar^\peps$ depends on the queue lengths, the projection $\vs_\parallel^\peps\defn \langle\vchi{1},\vsbar^\peps\rangle\vchi{1}$ is independent of $\vqbar^\peps$. Therefore, we do not need to assume that the scheduling problem is solved with MaxWeight. In general, for any special case of the generalized switch such that $\vs^\peps_\parallel$ is independent of the queue lengths, we can obtain a result similar to Proposition \ref{gs.prop:switch}, i.e., where we assume properties of the scheduling algorithm but not a specific algorithm.

\subsection{Proof of Theorem \ref{thm.generalized.switch}}

In the rest of this section we prove Theorem \ref{thm.generalized.switch} using the MGF method. Before presenting the proof, we introduce some notation.

Let $\Tbar$ and $\Bbar$ be steady-state random variables with the same distribution as $T(1)$ and $B_\ell(1)$, respectively.

\proof{Proof of Theorem \ref{thm.generalized.switch}.}

For ease of exposition we omit the dependence on $\epsilon$ of the variables in this proof. We use the MGF method. Similarly to the proof of Theorem \ref{thm.load.balancing}, we first need to verify that the prerequisites are satisfied.

\subsubsection*{Prerequisite 1. Positive recurrence.} 

In fact, MaxWeight algorithm is throughput optimal \citep{stolyar2004maxweight,atilla}. Then, for each $\epsilon>0$ the Markov chain $\left\{\vq^\peps(k):k\geq 1 \right\}$ is positive recurrent.

\subsubsection*{Prerequisite 2. SSC.}

Let $\cK=\left\{\vx\in\bR^n_+: \vx=\alpha\vcl\; ,\; \alpha\geq 0 \right\}$. Using the notation introduced in Prerequisite 2 in Section \ref{sec:details.MGFmethod}, we have $\vc=\vcl$, $\vqbar_\parallel^\peps=\langle\vcl,\vqbar^\peps\rangle\vcl$ and $\vqbar_\perp^\peps= \vqbar^\peps-\vqbar^\peps_{\parallel}$. \cite{atilla} proved that $\E{e^{\theta^* \left\|\vqbar_\perp\right\|}}$ is bounded for some finite $\theta^*$\footnote{In fact, the exponential moment bound is not part of the SSC statement of \cite{atilla}, but their proof of Proposition 2 implies it.}. Then, for each $m=1,2,\ldots$ there exists a constant $M_m$ such that $\E{\left\|\vqperpbar^\peps \right\|^m}\leq M_m$. Therefore, SSC as defined in Section \ref{sec:details.MGFmethod} is satisfied, and it occurs into the one-dimensional subspace $\cK$. In fact, in this case $\E{\left\|\vqperpbar^\peps \right\|^m}$ is $O(1)$, which is stronger.

Now we go through the steps of the MGF method.

\subsubsection*{Step 1. Prove an equation of the form of \eqref{eq.MGF.method.exp.qu} and compute an expression for the MGF of $\epsilon\langle\vc,\vqbar^\peps\rangle$.} 

We first prove Lemma \ref{lemma.gen.switch.expo.qu}.

\begin{lemma}\label{lemma.gen.switch.expo.qu}
	Consider a generalized switch parametrized by $\epsilon$ as described in Theorem \ref{thm.generalized.switch}. Then, for any real number $\theta$ such that $|\theteps|\leq \theta^*$ we have
	\begin{align*}
	\E{\left(e^{\theteps\langle\vcl,\left(\vqbar^\peps\right)^\plus\rangle} -1\right)\left(e^{-\theteps \langle\vcl,\vubar^\peps\rangle}-1\right)}\quad\text{is $o(\epsilon^2)$}
	\end{align*}
\end{lemma}

We present the proof of Lemma \ref{lemma.gen.switch.expo.qu} in Appendix \ref{app.gen.switch.expo.qu.proof}.

Before continuing, we need to prove that the MGF of $\epsilon \langle\vcl,\vqbar^\peps\rangle$ exists in an interval around 0. The proof is presented in Appendix \ref{app.gen.switch.MGF}. Then, following the steps sketched in Step 1 in Section \ref{sec:details.MGFmethod} we obtain \eqref{eq.MGFmethod.fraction}. 

When we applied the MGF method to the single server queue and to the load balancing system, we used the fact that the service rate vector is independent of the queue length vector to obtain \eqref{eq.ssq.fraction} and \eqref{eq.load.balancing.fraction}, respectively. However, in the case of the generalized switch this is no longer true. To overcome this difficulty we use the following lemma.

\begin{lemma}\label{lemma.gen.switch.qs}
	Consider a generalized switch operating under MaxWeight algorithm parametrized by $\epsilon$, as described in Theorem \ref{thm.generalized.switch}. Then, for any $\theta\in\bR$ we have
	\begin{align*}
	\E{\left(e^{\theteps\langle\vcl,\vqbar^\peps\rangle}-1\right) \left(e^{\theteps\left(\Bbar-\langle\vcl,\vsbar^\peps\rangle \right)}-1\right)}\quad\text{is $o(\epsilon^2)$.}
	\end{align*}
\end{lemma}

We present the proof in Appendix \ref{app.gen.switch.qs}. Working with the left hand side of \eqref{eq.MGFmethod.fraction} we obtain
\begin{align*}
& \E{e^{\theteps \langle\vcl, \vqbar\rangle}\left(1- e^{\theteps \langle\vcl,\vabar-\vsbar\rangle} \right) } \\
\stackrel{(a)}{=}& \E{e^{\theteps \langle\vcl,\vqbar\rangle}\left(1-e^{\theteps\left(\langle\vcl,\vabar\rangle-\Bbar \right)} \right) } + \E{e^{\theteps\langle\vcl,\vqbar\rangle}\left(e^{\theteps\left(
		\langle\vcl,\vabar-\vsbar\rangle \right)}-e^{\theteps\left(\langle\vcl,\vabar\rangle-\Bbar \right)} \right)} \\
\stackrel{(b)}{=}& \E{e^{\theteps\langle\vcl,\vqbar\rangle}}\left(1-\E{e^{\theteps\left( \langle\vcl,\vabar\rangle-\Bbar \right)}}\right) - \E{e^{\theteps\left(\langle\vcl,\vabar\rangle-\Bbar\right)}}\left(1-\E{e^{\theteps (\Bbar-\langle\vcl,\vsbar\rangle)} } \right) \\
&\; + \E{e^{\theteps \left(\langle\vcl,\vabar\rangle-\Bbar \right)}}\E{\left(e^{\theteps\langle\vcl,\vqbar\rangle} -1\right)\left(e^{\theteps\left(\Bbar-\langle\vcl,\vsbar\rangle\right)} -1\right)} \\
\stackrel{(c)}{=}& \E{e^{\theteps\langle\vcl,\vqbar\rangle}}\left(1-\E{e^{\theteps\left( \langle\vcl,\vabar\rangle-\Bbar \right)}}\right) - \E{e^{\theteps\left(\langle\vcl,\vabar\rangle-\Bbar\right)}}\left(1-\E{e^{\theteps (\Bbar-\langle\vcl,\vsbar\rangle)} } \right) + o(\epsilon^2),
\end{align*} 
where $(a)$ holds after adding and subtracting $\E{e^{\theteps \left(\langle\vcl,\vqbar\rangle + \langle\vcl,\vabar\rangle - \Bbar \right)}}$, and reorganizing terms; $(b)$ holds because $\vabar$ and $\Bbar$ are independent of the queue lengths vector $\vqbar$ and the potential service vector $\vsbar$, and after adding and subtracting $\E{e^{\theteps\left(\langle\vcl,\vabar\rangle-\Bbar \right)}}\E{e^{\theteps\left(\Bbar-\langle\vcl,\vsbar\rangle \right)}-1}$; and $(c)$ holds by Lemma \ref{lemma.gen.switch.qs} and because $\vabar$ and $\Bbar$ are bounded. Reorganizing terms we obtain
\begin{align}\label{eq.gen.switch.fraction}
\E{e^{\theteps\langle\vcl,\vqbar\rangle}}=& \dfrac{1- \E{e^{-\theteps\langle\vcl,\vubar\rangle}}+ \E{e^{\theteps\langle\vcl,\vabar\rangle}} \E{e^{-\theteps\langle\vcl,\vsbar\rangle}- e^{-\theteps\Bbar}} +o(\epsilon^2)}{1-\E{e^{\theteps\left(\langle\vcl,\vabar\rangle-\Bbar\right)}}}.
\end{align}

\subsubsection*{Step 2. Bound unused service and take heavy-traffic limit.}

The right hand side of \eqref{eq.gen.switch.fraction} yields a $\frac{0}{0}$ form in the limit as $\epsilon\dto 0$. Then, we take Taylor expansion of each of its terms, using Lemma \ref{lemma.taylor}. Similar to the case of the load balancing system, in this step we need to obtain bounds on $\E{\langle\vcl,\vubar\rangle}$. In this case we use the following lemma. 

\begin{lemma} \label{lemma.gen.switch.linear.test.func.}
	Consider a generalized switch parametrized by $\epsilon$ as described in Theorem \ref{thm.generalized.switch}. Then,
	\begin{align*}
	\E{\langle\vc,\vubar^\peps\rangle}+\bl- \E{\langle\vcl,\vsbar^\peps\rangle}=\epsilon.
	\end{align*}
\end{lemma}

\proof{Proof of Lemma \ref{lemma.gen.switch.linear.test.func.}.}

We set to zero the drift of $V_1(\vq)=\langle\vcl,\vq\rangle$. We obtain
\begin{align}
0=& \E{\langle\vcl,\vqbar^\plus\rangle-\langle\vcl,\vqbar\rangle} \nonumber\\
=& \E{\langle \vcl, \vqbar+\vabar-\vsbar+\vubar\rangle-\langle\vcl,\vqbar\rangle} \nonumber\\
=& \E{\langle\vcl,\vabar\rangle}-\E{\langle\vcl,\vsbar\rangle}+\E{\langle\vcl,\vubar\rangle}.\label{eq.gen.switch.lin.function}
\end{align}

Now, observe that
\begin{align}
\E{\langle\vcl,\vabar\rangle}=& \langle \vcl,\vrl-\epsilon \vcl\rangle \nonumber\\
=& \langle\vcl,\vrl\rangle -\epsilon\left\|\vcl \right\|^2 \nonumber\\
\stackrel{(a)}{=}& \bl-\epsilon , \label{eq.gen.switch.Eca}
\end{align}
where $(a)$ holds because $\vrl\in\cFl$ and because $\left\|\vcl \right\|=1$.

Then, using \eqref{eq.gen.switch.Eca} in \eqref{eq.gen.switch.lin.function} and rearranging terms we obtain the result. \Halmos
\endproof

Now we expand each term in the right hand side of \eqref{eq.gen.switch.fraction}. For the first term in the numerator, we have
\begin{align}
1-\E{e^{-\theteps\langle\vcl,\vubar\rangle}}=& 1-\E{f_{\epsilon,-\langle\vcl,\vubar\rangle}(\theta)} \nonumber\\
=& \theteps\E{\langle\vcl,\vubar\rangle}- \dfrac{(\theteps)^2}{2}\E{\left(\langle\vcl,\vubar\rangle \right)^2}+O(\epsilon^3). \label{eq.gen.switch.num1.partial}
\end{align}

In this case the numerator has more terms than in the case of the single server queue and the load balancing system, so we will keep the first moment of the unused service in the equation in order to use Lemma \ref{lemma.gen.switch.linear.test.func.}. However, we still need to bound the second moment.

\begin{claim}\label{claim.gen.switch.u}
	Consider a generalized switch as described in Theorem \ref{thm.generalized.switch}. Then,
	\begin{align*}
	\dfrac{(\theteps)^2}{2}\E{\left(\langle\vcl,\vubar\rangle \right)^2}\;\text{is }O(\epsilon^3)
	\end{align*}
\end{claim}

We present a proof of Claim \ref{claim.gen.switch.u} in Appendix \ref{app.claim.gen.switch.u}. Then, using Claim \ref{claim.gen.switch.u} in \eqref{eq.gen.switch.num1.partial} we obtain
\begin{align}
1-\E{e^{-\theteps\langle\vcl,\vubar\rangle}}=\theteps\E{\langle\vcl,\vubar\rangle}+O(\epsilon^3). \label{eq.gen.switch.num1}
\end{align}

For the second term in the numerator, we have
\begin{align}
\E{e^{\theteps\langle\vcl,\vabar\rangle}}\E{e^{-\theteps\langle\vcl,\vsbar\rangle} -e^{-\theteps\Bbar}}
=& \E{e^{\theteps\left(\langle\vcl,\vabar\rangle-\Bbar\right)}} \E{e^{\theteps\left(\Bbar-\langle\vcl,\vsbar\rangle\right)}-1} \nonumber \\
=& \E{f_{\epsilon,\left(\langle\vcl,\vabar\rangle-\Bbar\right)}(\theta)} \E{f_{\epsilon,\left(\Bbar-\langle\vcl,\vsbar\rangle\right)}(\theta)-1} \label{eq.gen.switch.num2.partial}
\end{align}

\begin{claim}\label{claim.gen.switch.num2}
	Consider a generalized switch as described in Theorem \ref{thm.generalized.switch} and the notation introduced in Lemma \ref{lemma.taylor}. Then,
	\begin{align*}
	&\E{f_{\epsilon,\left(\langle\vcl,\vabar\rangle-\Bbar\right)}(\theta)}= 1+\theteps^2+O(\epsilon^3)\\
	\text{and}\quad&
	\E{f_{\epsilon,\left(\Bbar-\langle\vcl,\vsbar\rangle\right)}(\theta)-1}=\theteps\E{\Bbar-\langle\vcl,\vsbar\rangle}+O(\epsilon^3)
	\end{align*}
\end{claim}
We prove the claim in Appendix \ref{app.claim.gen.switch.num2}. Using Claim \ref{claim.gen.switch.num2} in \eqref{eq.gen.switch.num2.partial}, reorganizing terms and using that $\Bbar$ and $\sbar_i$ are bounded for all $i\in\{1,\ldots,n\}$, we obtain
\begin{align}
\E{e^{\theteps\langle\vcl,\vabar\rangle}}\E{e^{-\theteps\langle\vcl,\vsbar\rangle} -e^{-\theteps\Bbar}}=& \theteps\E{\Bbar-\langle\vcl,\vsbar\rangle}+O(\epsilon^3)
\end{align}

Then, the numerator of \eqref{eq.gen.switch.fraction} yields
\begin{align}
& 1- \E{e^{-\theteps\langle\vcl,\vubar^\peps\rangle}}+ \E{e^{\theteps\langle\vcl,\vabar\rangle}} \E{e^{-\theteps\langle\vcl,\vsbar\rangle}- e^{-\theteps\Bbar}} +o(\epsilon^2) \nonumber \\
=& \left(\theteps\E{\langle\vcl,\vubar\rangle}+\theteps\E{\Bbar-\langle\vcl,\vsbar\rangle} +O(\epsilon^3)\right) + o(\epsilon^2) \nonumber \\
\stackrel{(a)}{=}& \theteps\left(\E{\langle\vcl,\vubar\rangle+\Bbar-\langle\vcl,\vsbar\rangle}\right)+o(\epsilon^2) \nonumber \\
\stackrel{(b)}{=}& \theteps^2+o(\epsilon^2),\label{eq.gen.switch.numerator}
\end{align}
where $(a)$ holds because $O(\epsilon^3)$ is $o(\epsilon^2)$; and $(b)$ holds by Lemmas \ref{gen.switch.lemma.B.bl} and \ref{lemma.gen.switch.linear.test.func.}.

For the denominator, we obtain
\begin{align}
& 1-\E{e^{-\theteps\left(\Bbar-\langle\vc,\vabar\rangle\right) } } \nonumber \\
=&  1-\E{f_{\epsilon,\left(\langle\vc,\vabar\rangle-\Bbar\right)}(\theta) } \nonumber \\
=& -\theteps\E{\langle\vc,\vabar\rangle-\Bbar}- \dfrac{(\theteps)^2}{2}\E{\left(\Bbar-\langle\vc,\vabar\rangle\right)^2}+O(\epsilon^3) \nonumber \\
\stackrel{(a)}{=}& \theteps^2 - \dfrac{(\theteps)^2}{2}\left(\E{\langle\vcl,\vabar\rangle^2}+\E{\Bbar^2}-2\E{\langle\vcl,\vabar\rangle \Bbar} \right) +O(\epsilon^3) \nonumber \\
\stackrel{(b)}{=}& 
\theteps^2 - \dfrac{(\theteps)^2}{2}\left(\sum_{i=1}^n\sum_{j=1}^n \Cov{a^\peps_i,a^\peps_j} + \sigma_{B_\ell}^2+\left(\E{\langle\vcl,\vabar\rangle}-\E{\Bbar}\right)^2\right)+O(\epsilon^3) \nonumber \\
\stackrel{(c)}{=}& \theteps^2 -\dfrac{(\theteps)^2}{2}\left(\sum_{i=1}^n \sum_{j=1}^n \Cov{a^\peps_i,a^\peps_j} +\sigma_{B_\ell}^2+\epsilon^2\right)+ O(\epsilon^3) \label{eq.gen.switch.denominator}
\end{align}
where $(a)$ holds by \eqref{eq.gen.switch.Eca} and expanding the square; $(b)$ holds by definition of variance and covariance, because $\vabar$ and $\Bbar$ are independent, and reorganizing terms; and $(c)$ holds by \eqref{eq.gen.switch.Eca}.

Using \eqref{eq.gen.switch.numerator} and \eqref{eq.gen.switch.denominator} in \eqref{eq.gen.switch.fraction} we obtain
\begin{align*}
\E{e^{\theteps\langle\vcl,\vqbar\rangle}}=& \dfrac{\theteps^2+o(\epsilon^2)}{\theteps^2 -\dfrac{(\theteps)^2}{2}\left(\ds\sum_{i=1}^n \sum_{j=1}^n \Cov{a^\peps_i,a^\peps_j} +\sigma_{B_\ell}^2+\epsilon^2\right)+ O(\epsilon^3)} \\
=& \dfrac{1+o(1)}{1- \dfrac{\theta}{2} \left(\ds\sum_{i=1}^n \sum_{j=1}^n \Cov{a^\peps_i,a^\peps_j} +\sigma_{B_\ell}^2+\epsilon^2\right)+O(\epsilon)}.
\end{align*}

Then, taking the heavy-traffic limit yields
\begin{align*}
\lim_{\epsilon\dto 0}\E{e^{\theteps\langle\vc,\vqbar\rangle}}= \dfrac{1}{1-\frac{\theta}{2}\left(\sum_{i=1}^n \sum_{j=1}^n \Cov{a_i,a_j} +\sigma_{B_\ell}^2\right)},
\end{align*}
which is the MGF of an exponential random variable with mean $\frac{1}{2}\left(\sum_{i=1}^n \sum_{j=1}^n \Cov{a_i,a_j} +\sigma_{B_\ell}^2 \right)$. This implies that $\vqparbar^\peps= \langle\vcl,\vqbar^\peps\rangle\vcl \Rightarrow \overline{\Upsilon}\vcl$, where $\overline{\Upsilon}$ is an exponential random variable with mean $\frac{1}{2}\left(\sum_{i=1}^n \sum_{j=1}^n \Cov{a_i,a_j} +\sigma_{B_\ell}^2\right)$.

Then, we conclude that $\epsilon\vqbar^\peps=\epsilon\vqparbar^\peps+\epsilon\vqperpbar^\peps$ converges in distribution to $\overline{\Upsilon}\vcl$ as $\epsilon\dto 0$. This proves Theorem \ref{thm.generalized.switch}. \Halmos
\endproof

\section{Future work}\label{sec:future.work}
The current paper develops the MGF method, which we believe can be used to study more general set of queueing systems. We outline a few of such future directions in this section.

In this paper we assumed that the number of arrivals and services in one time slot are bounded. We believe that this assumption is not required, and it is sufficient to assume that the first two moments of the arrival and service sequences exist. Relaxing these assumptions is an immediate future work. We will explore two paths for this generalization. One is the use of Characteristic Functions or one-sided Laplace transforms instead of MGF, since they always exist for nonnegative random variables. The main challenge in this approach is to establish the SSC under unbounded arrivals and service sequences. In the current paper, we used the SSC established by \cite{atilla}, which is based on the results from \cite{hajek_drift}, where the existence of all the moments of the arrival and service processes is assumed. We will explore ways to relax this assumption. The second approach that we will pursue is the MGF truncation arguments, similar to the ones introduced by \cite{braverman_truncation} for Markov Decision Processes. The main idea of their method is to take second order Taylor expansion of the value function in order to solve the Bellman equations. We believe this can give us insight to work with the second order Taylor expansion of the MGF.

Another question for future research is to use the MGF method to study the rate of convergence to the heavy-traffic limit. In addition to obtaining the results on the heavy-traffic limiting behavior, the Drift method also gives upper and lower bounds that are applicable in all traffic \citep{atilla,MagSri_SSY16_Switch,QUESTA_switch}.  These bounds give the rate of convergence to the heavy-traffic limit.
Since the MGF method is a natural generalization of the Drift method, it may be used to obtain results on rate of convergence too, which is a topic for future study.

The next set of future work is on developing the MGF method for its use in systems that do not satisfy the CRP condition, and this will be the culmination of the present work because the main motivation in developing the MGF method is to study systems when the CRP condition is not met. We believe that the MGF method is a promising approach to obtain the heavy-traffic distribution of the queue lengths when CRP condition does not hold, even though the Drift method is known to fail in this case \citep{Hurtado_gen-switch_temp}, because of the following reason. The queue lengths process is a multi-dimensional Discrete Time Markov Chain (DTMC) (or a continuous Markov Chain in some cases). For a positive recurrent and irreducible DTMC, it is known that the stationary distribution exists and is unique. One first establishes positive recurrence of the DTMC using Foster-Lyapunov Theorem. This has an added benefit that one typically obtains as a consequence a (possibly loose) upper bound on an expression of them form $E[\epsilon\sum_i\qbar_i]$.
If $P$ is the transition matrix, then the stationary distribution is a unique solution of the equation, $\pi=\pi P$. Clearly, solving for the stationary distribution in general is hard. However, we know that it is unique and is characterized by this equation. If we take two-sided Laplace transform of the equation $\pi=\pi P$ we obtain an equation which is same as the one we obtain by setting the drift of the exponential test function to zero. Since Laplace transform is invertible, solving this equation uniquely characterizes the stationary distribution through its MGF. However, as shown in Section \ref{sec:single.server.queue}, even for the single server queue it is challenging to obtain a solution for this equation in all traffic (see Equation \eqref{eq.ssq.fraction}). Therefore, using the MGF approach, we seek to solve it in the heavy-traffic limit. To do this, one first needs to prove tightness of the sequence of the stationary distributions as the heavy-traffic parameter $\epsilon$ goes to zero.  Tightness follows directly from the bound on $E[\epsilon\sum_i\qbar_i]$ that one obtains from the Foster-Lyapunov Theorem. Therefore, we expect that the MGF drift equation that we have in the heavy-traffic limit must have a unique solution. Typically, since the system is tractable in steady-state, we expect to solve this equation explicitly to get the joint stationary distribution in steady-state. Even in cases when this equation may not be solved explicitly, one may be able to obtain moments from this equation. For instance, one may be able to obtain the moment bounds computed by \cite{MagSri_SSY16_Switch}, \cite{QUESTA_switch} and \cite{Weina_bandwidth_journal} from such an equation.

Two systems of special interest that do not satisfy the CRP condition are the  bandwidth-sharing network operating under proportional scheduling and the input-queued crossbar switch operating under MaxWeight. The bandwidth-sharing network \citep{massoulierobertsbandwidth} operating under the so-called proportional scheduling algorithm is a good model for studying flow level dynamics in data centers.
If the arrivals are Poisson and job-sizes are exponential, it is known that the stationary distribution in heavy-traffic is product of exponentials \citep{kang2009state,yeyaobandwidth2012}.
The bandwidth sharing network is one of the simplest systems that does not satisfy the CRP condition because of this product form structure.
It is also known that the stationary distribution of the corresponding RBM in the diffusion limit is insensitive to the job size distribution as long as it belongs to the class of phase-type distributions, which are known to be dense in the space of distributions \citep{zwart_bandwidth_diffusion}. However, the interchange of limits step was not shown by \cite{zwart_bandwidth_diffusion}, so their result does not show if the stationary distribution of the original system in heavy-traffic is also insensitive. Recently, the Drift method was used to complete this limit-interchange step \citep{Weina_bandwidth_journal}.
We will use the MGF method to directly study the stationary distribution in heavy-traffic under phase-type arrivals using the MGF method to show insensitivity, and to show that the stationary distribution is indeed the product of exponentials.

The input-queued cross bar switch is an idealized model of a data center network. It can be modeled as an $n\times n$ matrix of queues
where the rows represent the input ports and the columns represent the output ports. Therefore, the dimension of the state space is $n^2$. \cite{MagSri_SSY16_Switch} studied an input-queued cross-bar switch operating under MaxWeight and proved that SSC occurs onto a $(2n-1)$-dimensional cone. Moreover, the expected sum of the scaled queue lengths in heavy-traffic was obtained using the Drift method, resolving an open conjecture. Characterizing the higher moments and the distribution (marginals and joint) of scaled queue lengths are still open questions. The MGF method is developed in this paper with the goal of answering these questions given the limitation of the Drift method to solve these problems \citep{Hurtado_gen-switch_temp}.

\section{Conclusion}\label{sec:conclusion}
In this paper we introduced transform methods to compute the steady-state distribution of the scaled queue lengths in heavy-traffic. We focused on two-sided Laplace transform, which is also known as Moment Generating Function (MGF). We motivated the method with a single server queue and we applied it in queueing systems that satisfy the CRP condition, such as load balancing systems and the generalized switch. The main idea in the MGF method is to set the drift on an exponential test function to zero. The key step is in getting a handle on the unused service, and the paper illustrates how the unused service is handled in two different types of queueing systems.
Further developing the MGF method to study system when the CRP condition is not satisfied such as the bandwidth sharing network and the input-queued switch forms future work.

\newpage

\renewcommand{\theHsection}{A\arabic{section}}
\begin{APPENDIX}{}
	
\section{Proof of Lemma \ref{lemma.taylor}}\label{app.taylor}
\proof{Proof of Lemma \ref{lemma.taylor}.}

Fix $\Theta>0$ and $x\in\bR$. Then, from Taylor approximation of $f_{\epsilon,x}(\theta)=e^{\theteps x}$ at $\theta=0$ we have
\begin{align*}
	e^{\theteps x}\leq 1+\theteps x+\dfrac{(\theteps)^2}{2}x^2 + \dfrac{(\tilde{\theta}\epsilon)^3}{3!}x^3\quad\forall \theta\in[-\Theta,\Theta],\,\forall x\in \bR,
\end{align*}
where $\tilde{\theta}$ is a real number between 0 and $\theta$. Then, for all $0\leq x\leq K$ we have
\begin{align*}
	e^{\theteps x}\leq 1+\theteps x+\dfrac{(\theteps)^2}{2}x^2 + \dfrac{(\tilde{\theta}\epsilon)^3}{3!}K^3.
\end{align*}
Since $\tilde{\theta}$ is between 0 and $\theta$, and $|\theta|\leq \Theta$ we have
\begin{align*}
	\left|\dfrac{(\tilde{\theta}\epsilon)^3}{3!}K^3 \right|=\dfrac{|\tilde{\theta}|^3\epsilon^3}{3!}K^3\leq \dfrac{(\Theta\epsilon)^3}{3!}K^3,
\end{align*}
which is finite for every $\epsilon$. Then,
\begin{align*}
	e^{\theteps x}\leq 1+\theteps x+\dfrac{(\theteps)^2}{2}x^2 + \dfrac{(\Theta\epsilon)^3}{3!}K^3.
\end{align*}
Therefore,
\begin{align*}
	\left| e^{\theteps x}-1-\theteps x-\dfrac{(\theteps)^2}{2}x^2\right|\leq C_1 \epsilon^3,
\end{align*}
where $C_1=\dfrac{\Theta^3 K^3}{3!}$ is a finite constant.

Now, since $X^\peps\leq K_{\max}$ with probability 1, we have
\begin{align*}
	\E{e^{\theteps X^\peps}}\leq 1+\theteps \E{X^\peps}+\dfrac{(\theteps)^2}{2}\E{\left(X^\peps\right)^2}+ \dfrac{\Theta \epsilon^3 K_{\max}}{3!},
\end{align*}
which proves the lemma. \Halmos
\endproof
	
\section{Details of the proofs in Section \ref{sec:load.balancing}}
	In this section we provide the details of the proofs of the lemmas stated in Section \ref{sec:load.balancing}.
	
	\subsection{Proof of SSC in the load balancing system operating under JSQ}\label{app.jsq.ssc}
	In this section we present an insight of the proof of SSC as developed in \cite{atilla}. They prove the result for the case where the servers are independent, but it also holds in the case where they are not. We first state the result.

\begin{proposition}\label{prop.load.balancing.jsq.ssc}
	Consider a load balancing system as described in Corollary \ref{cor.load.balancing.jsq}. Then, for each $m=1,2,\ldots$ there exists a finite constant $M_m$ such that 
	\begin{align*}
		\E{\left\|\vqperpbar^\peps \right\|^m}\leq M_m.
	\end{align*}
\end{proposition}

This proof is based on a lemma that was first proved by \cite{hajek_drift}. The original statement is more general than what we need here, so we present the specific result that we will use, as stated by \cite{atilla}.
	
\begin{lemma}\label{hajek.lemma}
	For an irreducible and aperiodic Markov Chain $\{X(k):\,k\geq 1\}$ over a countable state space $\mathcal{X}$, suppose $Z:\mathcal{X}\to \bR_+$ is a nonnegative valued Lyapunov function. The drift of $Z$ at $x$ is
		\ba \Delta Z(x)\stackrel{\triangle}{=}\big[Z\big(X(k+1)\big)-Z\big(X(k)\big) \big]\ind{X(k)=x} \ea
	Thus, $\Delta Z(x)$ is a random variable that measures the amount of change in the value of $Z$ in one step, starting from state $x$. This drift is assumed to satisfy the following conditions:
	\begin{enumerate}[label=(C\arabic*)]
		\item There exists $\eta>0$ and $\kappa<\infty$ such that
		\ba \E{\left.\Delta Z(x)\,\right|\,X(k)=x}\leq -\eta\quad\text{for all $x\in\mathcal{X}$ with $Z(x)\geq \kappa$} \ea
		\item There exists $D<\infty$ such that
		\ba |\Delta Z(x)|\leq D\quad\text{ with probability 1 for all $x\in\mathcal{X}$} 	\ea
	\end{enumerate}
	Then, there exist $\theta^*>0$ and $C^*<\infty$ such that
	\ba \limsup_{k\to\infty} \E{e^{\theta^* Z(X(k))}}\leq C^* \ea
	If we further assume that the Markov chain $\{X(k):\,k\geq 1\}$ is positive recurrent, then $Z(X(k))$ converges in distribution to a random variable $\overline{Z}$ for which
	\ba \E{e^{\theta^* \overline{Z}}}\leq C^* \ea
\end{lemma} 

\proof{Proof of Proposition \ref{prop.load.balancing.jsq.ssc}.}
\cite{atilla} use the Lyapunov function $Z(\vq)=\|\vqperp^\peps\|$ and they prove that
\ba \E{\Delta Z(\vq)\,|\, \vq(k)=\vq}\leq -\delta+\dfrac{n(\max\{\amax,\smax \})^2+2n\smax^2}{2\|\vqperp^\peps\|},\ea
where $\delta$ is a fixed constant in $(0,\mu_{\min})$. The proof is based on the fact that $\left\|\vx\right\|=\sqrt{\left\|\vx \right\|^2}$, that square root is a concave function and that JSQ sends all arrivals to the shortest queue in each time slot. This verifies condition (C1) of Lemma \ref{hajek.lemma}.

To verify condition (C2), they prove that for all $\vq\in\bR^n_+$
\ba |\Delta Z(\vq)|\leq 2\sqrt{n}\max\{\amax,\smax \}, \ea 
using triangle inequality and boundedness of the arrival and service processes.

Also, for $\epsilon>0$ the Markov Chain $\{\vq(k):\,k\geq 1\}$ is positive recurrent. Also, since projection is nonexpansive we have $\left\|\vqperp^\peps(k) \right\|\leq \left\| \vq^\peps(k)\right\|$, which implies that $\{\vq_\perp (k):\,k\geq 1\}$ is positive recurrent. Therefore, by Lemma \ref{hajek.lemma} there exists $\theta^*>0$ and $C^*>0$ such that 
\ba \E{e^{\theta^* \|\vqperpbar^\peps\| }}\leq C^* \ea
Finally, since $ \left\|\vqperpbar^\peps\right\|\geq 0$ and $f(x)=e^x$ is a nonnegative increasing function, we obtain that $\E{e^{\theta \|\vqperpbar^\peps\| }}\leq C^*$ for all $\theta\in[-\theta^*,\theta^*]$. This implies that for each $m=1,2,\ldots$ 
\begin{align*}
	\E{\left\|\vqbar_\perp^\peps\right\|^m}\leq M_m
\end{align*}
\Halmos

\endproof
	
	\subsection{Existence of MGF of $\epsilon\sum_{i=1}^n \qbar_i^\peps$ in the load balancing system operating under JSQ}\label{app.jsq.mgf}
	We first state the result formally.

\begin{lemma}\label{lemma:mgf-jsq}
	Consider a load balancing system operating under JSQ, parametrized by $\epsilon\in(0,\mu_\Sigma)$ as described in Corollary \ref{cor.load.balancing.jsq}. Then, for each $\epsilon\in(0,\mu_\Sigma)$ there exists $\Theta>0$ such that $\E{e^{\theteps \sum_{i=1}^n \qbar_i^\peps}}<\infty$ for all $\theta\in[-\Theta,\Theta]$.
\end{lemma}

\proof{Proof of Lemma \ref{lemma:mgf-jsq}.}
We omit the dependence on $\epsilon$ of the variables for ease of exposition. First observe that if $\theta\leq 0$, then $\E{e^{\theteps\sum_{i=1}^n \qbar_i}}<\infty$ trivially because $\vqbar\geq \vzero$ by definition of queue length.

In the rest of this proof we assume $\theta>0$. Observe that the function $f(x)=e^{\theteps x}$ is convex. Then, by Jensen's inequality we have that, for all $\vq\geq \vzero$
\begin{align*}
	e^{\frac{\theteps}{n} \sum_{i=1}^n q_i}\leq \dfrac{1}{n}\sum_{i=1}^n e^{\theteps q_i}.
\end{align*}

Hence, it suffices to show that $\sum_{i=1}^n \E{e^{\theteps\qbar_i}}<\infty$ for $\theta<\Theta$. We use Foster-Lyapunov theorem \cite[Proposition 6.13]{hajekrandomprocbook} with Lyapunov function $V(\vq)=\sum_{i=1}^n e^{\theteps q_i}$.

Using Lemma \ref{lemma.ssq.expoqu} for each of the $n$ queues and rearranging terms we obtain that, for each $i\in[n]$ and $k\geq 1$
\begin{align*}
	e^{\theteps q_i(k+1)} =& 1-e^{-\theteps u_i(k)} + e^{\theteps\left(q_i(k)+a_i(k)-s_i(k) \right)}
\end{align*}

Then, using the notation $\Evq{\cdot}\defn \E{\cdot \,|\,\vq(k)=\vq}$, we obtain
\begin{align*}
	\Evq{V\left(\vq(k+1)\right)- V\left(\vq(k) \right)}=& \sum_{i=1}^n \Evq{e^{\theteps q_i(k+1)}-e^{\theteps q_i(k)}} \\
	=& \sum_{i=1}^n \left(1-\Evq{e^{-\theteps u_i(k)}}\right) + \sum_{i=1}^n e^{\theteps q_i}\left(\Evq{e^{\theteps \left(a_i(k)-s_i(k)\right)}}-1 \right).
\end{align*}

Observe that, since $\Evq{e^{-\theteps u_i(k)}}\geq 0$ we have
\begin{align*}
	\sum_{i=1}^n \left(1-\Evq{e^{-\theteps u_i(k)}}\right)\leq n.
\end{align*}
Then, it suffices to show that for some $\Theta$ and some $\eta>0$, we have
\begin{align*}
	\sum_{i=1}^n e^{\theteps q_i}\left(\Evq{e^{\theteps \left(a_i(k)-s_i(k)\right)}}-1 \right)\leq -\eta \qquad\forall \theta\in(0,\Theta].
\end{align*}

Given $\vq(k)=\vq$, let $i^*\in\argmin_{i\in\{1,\ldots,n\}} \left\{q_i(k) \right\}$ be the queue where arrivals in time slot $k$ are routed. Then,
\begin{align*}
	\sum_{i=1}^n e^{\theteps q_i}\left(\Evq{e^{\theteps \left(a_i(k)-s_i(k)\right)}}-1 \right) =& e^{\theteps q_{i^*}}\left(\E{e^{\theteps (a(k)-s_{i^*}(k))}}-1\right) + \sum_{\substack{ i=1 \\ \;\; i\neq i^*}}^n e^{\theteps q_i}\left(\E{e^{-\theteps s_i(k)}} -1\right) \\
	=& e^{\theteps q_{i^*}}\theta M'_{a-s_{i^*}}(\xi_{i^*}) + \sum_{\substack{ i=1 \\ \;\; i\neq i^*}}^n e^{\theteps q_i}\left(- \theta M'_{s_i}(\xi_i) \right),
\end{align*}
where we used the notation $M_X(\theta)=\E{e^{\theteps X}}$ and $\xi_i,\xi_{i^*}$ are numbers between 0 and $\theta$ for all $i\neq i^*$. The second equality holds by Taylor expansion up to first order of $M_{a-s_{i^*}}(\theta)$ and $M_{s_i}(\theta)$ for all $i\neq i^*$, around $\theta=0$.

Also, observe $M'_{a-s_{i^*}}(0)=\epsilon(\lambda-\mu_{i^*})$ and $M'_{s_i}(0)=\epsilon\mu_i$ for all $i\neq i^*$, and MGF is continuous at $\theta=0$ \cite[p. 78]{mood}. Then, for each $i\in\{1,\ldots,n\}$ there exists $\Theta_i$ such that 
\begin{align*}
	& M'_{s_i}(\xi_i)\geq \dfrac{\epsilon \mu_i}{2}\quad \forall |\theta|<\Theta_i \quad\text{for each }i\neq i^*\\ \quad\text{and}\quad& \left|M'_{a-s_{i^*}}(\xi_{i^*}) \right|\leq \left|\dfrac{\epsilon \left(\lambda-\mu_{i^*}\right)}{2} \right| \quad\forall |\theta|<\Theta_{i^*}.
\end{align*}

Let $\Theta=\min_{i=1,\ldots,n} \Theta_i$. Then, for all $\theta\in\left(0,\Theta\right]$ we have
\begin{align*}
	\Evq{V\left(\vq(k+1)\right)- V\left(\vq(k) \right)}\leq&n+ e^{\theteps q_{i^*}}\left(\dfrac{\theteps(\lambda-\mu_{i^*})}{2}\right) - \sum_{\substack{ i=1 \\ \;\; i\neq i^*}}^n e^{\theteps q_i}\left(\dfrac{\theteps \mu_i}{2}\right) \\
	\stackrel{(a)}{=}& n+\dfrac{\theteps}{2}\sum_{i=1}^n \lambda_i\left(e^{\theteps q_{i^*}} - e^{\theteps q_i} \right) + \dfrac{\theteps}{2}\sum_{i=1}^n e^{\theteps q_i}\left(\lambda_i-\mu_i\right) \\
	\stackrel{(b)}{\leq}& n+\sum_{i=1}^n e^{\theteps q_i}\left(\dfrac{\theteps(\lambda_i-\mu_i)}{2}\right) \\
	\stackrel{(c)}{=}&n -\dfrac{\theteps^2}{2n}\sum_{i=1}^n e^{\theteps q_i}
\end{align*}
where $\lambda_i\defn \mu_i-\frac{\epsilon}{n}$. Here, $(a)$ holds by adding and subtracting $\sum_{i=1}^n e^{\theteps q_i}\left(\frac{\theteps\lambda_i}{2}\right)$, realizing that $\lambda=\sum_{i=1}^n \lambda_i$ and rearranging terms; $(b)$ holds because $q_{i^*}\leq q_i$ for all $i$ by definition of $i^*$; and $(c)$ holds because $\mu_i-\lambda_i=\frac{\epsilon}{n}$. This proves the lemma. \Halmos
\endproof
	
	\subsection{Proof of Lemma \ref{lemma.load.balancing.expoqu}}\label{app.load.balancing.expo.qu}
	To prove Lemma \ref{lemma.load.balancing.expoqu} we use the following result.
\begin{lemma}\label{lemma.load.balancing.to.expoqu}
	Consider the load balancing system indexed by $\epsilon$ described in Theorem \ref{thm.load.balancing}. Then, for any $\alpha\in\bR$ and for all $k\geq 1$ we have
	\begin{align*}
		\sum_{i=1}^n u_i^\peps(k)\left(e^{\frac{\alpha}{n}\sum_{j=1}^n q_j^\peps(k+1)}-1 \right)= \sum_{i=1}^n u_i^\peps(k)\left(e^{-\alpha q_{\perp i}^\peps(k+1)}-1 \right),
	\end{align*}
	where $q_{\perp i}^\peps(k)$ is the $i\tth$ element of $\vqperp^\peps(k)$, for each $i\in\{1,\ldots,n\}$.
\end{lemma}

\proof{Proof of Lemma \ref{lemma.load.balancing.to.expoqu}.}

If $\alpha=0$, the equation trivially holds. So now assume $\alpha\neq 0$. Since $q_i(k+1) u_i(k)=0$ for all $i\in\{1,\ldots,n\}$, we have
\begin{align*}
	u_i(k)(e^{-\alpha q_i(k+1)}-1)=0\quad\forall i\in\{1,\ldots,n\}.
\end{align*}

Then, summing over $i\in\{1,\ldots,n\}$ we obtain
\begin{align*}
	\sum_{i=1}^n u_i(k)\left(e^{-\alpha q_i(k+1)}-1\right)=0.
\end{align*}

By definition of $\vqpar(k)$ and $\vqperp(k)$ we have $\vq(k)=\vqpar(k)+\vqperp(k)$, so
\begin{align*}
	\sum_{i=1}^n u_i(k)(e^{-\alpha \left(q_{\parallel i}(k+1) +q_{\perp i}(k+1)\right)}-1)=0.
\end{align*}

But $\vqpar(k+1) = \left(\frac{1}{n}\sum_{j=1}^n q_j(k+1)\right)\vone$ so $q_{\parallel i}(k+1)=q_{\parallel 1}(k+1)$ for all $i\in\{1,\ldots,n\}$. Then, reorganizing terms we obtain
\begin{align*}
	\sum_{i=1}^n u_i(k) e^{-\alpha q_{\perp i}(k+1)}=e^{\alpha q_{\parallel 1}(k+1)}\sum_{i=1}^n u_i(k).
\end{align*}

By definition of $\vqpar(k)$ we obtain
\begin{align*}
	\sum_{i=1}^n u_i(k) e^{-\alpha q_{\perp i}(k+1)}=e^{\frac{\alpha}{n} \sum_{j=1}^n q_j(k+1)}\sum_{i=1}^n u_i(k).
\end{align*}

Finally, subtracting $\sum_{i=1}^n u_i(k)$ in both sides we obtain
\begin{align*}
	\sum_{i=1}^n u_i(k)\left(e^{\frac{\alpha}{n} \sum_{j=1}^n q_j(k+1)}-1 \right)=\sum_{i=1}^n u_i(k) \left(e^{-\alpha q_{\perp i}(k+1)}-1\right).
\end{align*}
\Halmos
\endproof

In the proof of Lemma \ref{lemma.load.balancing.expoqu} we use Lemma \ref{lemma.load.balancing.to.expoqu} and the following facts:

\begin{enumerate}[label=(\roman*)]
	\item\label{facts.functions.fraction}  The function $g(x)=\frac{e^x-1}{x}$ is nonnegative and nondecreasing for all $x\in\bR$
	\item\label{facts.functions.ineq} Suppose $0\leq x\leq y$. Then, for all $\theta\in\bR$ we have $\ds e^{\theta x}-1\leq (\theta x)\left(\dfrac{e^{\theta y}-1}{\theta y}\right)$
	\item\label{facts.functions.exp} For all $x\in\bR_+$, $\ds \dfrac{e^x-1}{x}<e^x$
\end{enumerate}

All these facts can be shown using calculus techniques, so we omit the proof. Now we prove Lemma \ref{lemma.load.balancing.expoqu}.

\proof{Proof of Lemma \ref{lemma.load.balancing.expoqu}.} 
First observe that if $\theta=0$ the statement trivially holds. If $\theta\neq 0$, by properties of expectation and absolute value we obtain
\begin{align}
	& \left|\E{\left(e^{\theteps \sum_{i=1}^n \qbar_i^\plus}-1 \right)\left(e^{-\theteps \sum_{i=1}^n \ubar_i}-1 \right)} \right| \nonumber \\
	\leq& \E{\left|\left(e^{\theteps \sum_{i=1}^n \qbar_i^\plus}-1 \right)\left(e^{-\theteps \sum_{i=1}^n \ubar_i}-1 \right)\right|} \nonumber  \\
	\stackrel{(a)}{=}& |\theta|\epsilon\E{\left| \left(\sum_{i=1}^n \ubar_i\right) \left(e^{\theteps \sum_{i=1}^n \qbar_i^\plus}-1 \right)\left(\dfrac{e^{-\theteps \sum_{i=1}^n \ubar_i}-1}{-\theteps \sum_{i=1}^n \ubar_i} \right)\right|\ind{\sum_{i=1}^n \ubar_i\neq 0} } \nonumber  \\
	\stackrel{(b)}{\leq}& |\theta|\epsilon \left(\dfrac{e^{|\theta|\epsilon n S_{\max}}-1}{|\theta|\epsilon n S_{\max}} \right) \E{\left|\sum_{i=1}^n \ubar_i\left(e^{\theteps \sum_{j=1}^n \qbar_j^\plus}-1 \right) \right|} \nonumber  \\
	\stackrel{(c)}{\leq}& |\theta|\epsilon \left(\dfrac{e^{|\theta|\epsilon n S_{\max}}-1}{|\theta|\epsilon n S_{\max}} \right) \E{\sum_{i=1}^n \ubar_i \left|e^{-\theteps n \qbar_{\perp i}}-1 \right| } \nonumber  \\
	\stackrel{(d)}{\leq}& |\theta|\epsilon \left(\dfrac{e^{|\theta|\epsilon \smax}-1}{|\theta|\epsilon \smax} \right) \E{\sum_{i=1}^n \ubar_i^p}^{\tfrac{1}{p}}  \E{\sum_{i=1}^n \left|e^{-\theta\epsilon n\qbar_{\perp i}} -1\right|^{\tfrac{p}{p-1}}}^{\frac{p-1}{p}} \nonumber \\
	\stackrel{(e)}{\leq}& |\theta| \epsilon^{1+\frac{1}{p}}\smax^{\frac{p-1}{p}} \left(\dfrac{e^{|\theta|\epsilon \smax}-1}{|\theta|\epsilon \smax} \right)\E{\sum_{i=1}^n \left|e^{-\theta\epsilon n\qbar_{\perp i}} -1\right|^{\tfrac{p}{p-1}}}^{\frac{p-1}{p}} \nonumber \\ 
	=& \theta^2 \epsilon^{2+\frac{1}{p}}\smax^{\frac{p-1}{p}} n \left(\dfrac{e^{|\theta|\epsilon \smax}-1}{|\theta|\epsilon \smax} \right) \left(\sum_{i=1}^n \E{\left|\dfrac{e^{-\theteps n \qbar_{\perp i}}-1}{-\theteps n \qbar_{\perp i}} \right|^\frac{p}{p-1} |\qbar_{\perp i}|^{\frac{p}{p-1}} \ind{\qbar_{\perp i}\neq 0}} \right)^{\frac{p-1}{p}},
 \label{eq.jsq.partial.eq.expoq.u.oeps2} 
\end{align}
where $p> 1$. Here $(a)$ holds because if $\sum_{i=1}^n \ubar_i=0$ then $e^{-\theteps \sum_{i=1}^n \ubar_i}-1=0$, and by multiplying and dividing everything by $\left|\theteps \sum_{i=1}^n \ubar_i\right|$; $(b)$ holds by the fact \ref{facts.functions.fraction} stated above, because $\ubar_i\leq S_{\max}$ for all $i\in\{1,\ldots,n\}$ and because $0\leq \ind{\sum_{i=1}^n\ubar_i\neq 0}\leq 1$; $(c)$ holds by triangle inequality and Lemma \ref{lemma.load.balancing.to.expoqu}; $(d)$ holds by Hölder's inequality; and $(e)$ holds because $\ubar_i\leq S_{\max}$ for all $i\in\{1,\ldots,n\}$, because $\sum_{i=1}^n \E{\ubar_i}=\epsilon$ and because $x^{\frac{1}{p}}$ is an increasing function for $x\geq 0$.

By L'Hospital's rule we have
\begin{align*}
	\lim_{\epsilon\dto 0} \dfrac{e^{|\theta|\epsilon n S_{\max}}-1}{|\theta|\epsilon n S_{\max}}=1
\end{align*}

Then, the last step is to prove that the last expression in  \eqref{eq.jsq.partial.eq.expoq.u.oeps2} is $O(1)$. To do that we show the following claim at the end of this section.

\begin{claim}\label{claim.jsq.last.step.expo.qu}
	Consider a load balancing system as described in Lemma \ref{lemma.load.balancing.expoqu}. Then, there exists $\theta_{\max}>0$ finite such that for all $|\theta|<\theta_{\max}$ we have
	\begin{align*}
		\left(\sum_{i=1}^n \E{\left|\dfrac{e^{-\theteps n \qbar_{\perp i}}-1}{-\theteps n \qbar_{\perp i}} \right|^\frac{p}{p-1} |\qbar_{\perp i}|^{\frac{p}{p-1}} \ind{\qbar_{\perp i}\neq 0}} \right)^{\frac{p-1}{p}} \text{ is } o(\epsilon^2).
	\end{align*}
	An expression for $\theta_{\max}$ is provided in \eqref{eq:claim4-thetamax}.
\end{claim}

Therefore,
\begin{align*}
	\E{\left(e^{\theteps \sum_{i=1}^n \qbar_i^\plus}-1 \right)\left(e^{-\theteps \sum_{i=1}^n \ubar_i}-1 \right)}\;\text{is }o(\epsilon^2)
\end{align*}
\Halmos
\endproof

Now we prove the claim.

\proof{Proof of Claim \ref{claim.jsq.last.step.expo.qu}.}
	By Hölder's inequality, for each $i\in\{1,\ldots,n\}$
	\begin{align*}
		\E{\left|\dfrac{e^{-\theteps \qbar_{\perp i}}-1}{-\theteps \qbar_{\perp i}} \right|^{\frac{p}{p-1}} |\qbar_{\perp i}|^{\frac{p}{p-1}} \ind{\qbar_{\perp i}\neq 0}} \leq \E{\left|\dfrac{e^{-\theteps \qbar_{\perp i}}-1}{-\theteps \qbar_{\perp i}} \right|^{\left(\frac{p}{p-1}\right)\left(\frac{\tilde{p}}{\tilde{p}-1}\right)} \ind{\qbar_{\perp i}\neq 0}}^{\frac{\tilde{p}-1}{\tilde{p}}} \E{ |\qbar_{\perp i}|^{\left(\frac{p}{p-1}\right)\tilde{p}} }^{\frac{1}{\tilde{p}}} ,
	\end{align*}
	where $\tilde{p}> 1$. On one hand, we can choose $p$ large so that $\frac{p}{p-1}\approx 1$, and $\tilde{p}> 1$ such that $\left(\frac{p}{p-1}\right)\tilde{p}= 2$. 
	Then, $\E{|\qbar_{\perp i}|^{\left(\frac{p}{p-1}\right)\tilde{p}} }^{\frac{1}{\tilde{p}}}$ is $o\left(\epsilon^2\right)$ by SSC. 
	
	Also, by SSC we know that $\epsilon|\qbar_{\perp i}|$ converges to zero in the mean-square sense and, therefore, in distribution. Then, by the continuous mapping theorem \cite[Theorem 10.4 in Section 5]{gut2012probability} we have that 
	\begin{align*}
		 \left(\dfrac{e^{-\theteps|\qbar_{\perp i}|}-1}{-\theteps|\qbar_{\perp i}|}\right)^{\left(\frac{p}{p-1}\right)\left(\frac{\tilde{p}}{\tilde{p}-1}\right)} \Longrightarrow 1.
	\end{align*}
	It remains to prove that $\frac{e^{-\theteps|\qbar_{\perp i}|}-1}{-\theteps|\qbar_{\perp i}|}$ is bounded to conclude that its expected value also converges to 1. In fact, we have
	\begin{align*}
		-\theta\epsilon|\qbar_{\perp i}|\leq |\theta|\epsilon|\qbar_{\perp i}|\leq |\theta|\epsilon\|\vqbar_{\perp}\|
	\end{align*}
	and $|\theta|\epsilon\|\vqbar_{\perp}\|\geq 0$. Then, by the facts \ref{facts.functions.fraction} and \ref{facts.functions.exp} stated above we obtain
	\begin{align*}
	0\leq \dfrac{e^{-\theteps|\qbar_{\perp i}|}-1}{-\theteps|\qbar_{\perp i}|}\ind{\qbar_{\perp i}\neq 0} \leq \dfrac{e^{|\theta|\epsilon\|\vqbar_{\perp}\|}-1}{|\theta|\epsilon\|\vqbar_{\perp}\|}\ind{\qbar_{\perp i}\neq 0} \leq e^{|\theta|\epsilon\|\vqbar_{\perp}\|}\ind{\qbar_{\perp i}\neq 0} \leq e^{|\theta|\epsilon\|\vqbar_{\perp}\|}
	\end{align*}
	Therefore,
	\begin{align*}
		\E{\left(\dfrac{e^{-\theteps|\qbar_{\perp i}|}-1}{-\theteps|\qbar_{\perp i}|}\right)^{\left(\frac{p}{p-1}\right)\left(\frac{\tilde{p}}{\tilde{p}-1}\right)}\ind{\qbar_{\perp i}\neq 0}}\leq& \E{e^{|\theta|\left(\frac{p}{p-1}\right)\left(\frac{\tilde{p}}{\tilde{p}-1}\right)\epsilon\|\vqbar_\perp \|}} \\
		\stackrel{(a)}{\leq}& \E{e^{|\theta|\left(\frac{p}{p-1}\right)\left(\frac{\tilde{p}}{\tilde{p}-1}\right)\epsilon\|\vqbar \|}} \\
		\stackrel{(b)}{\leq}& \E{e^{|\theta|\left(\frac{p}{p-1}\right)\left(\frac{\tilde{p}}{\tilde{p}-1}\right)\epsilon\sum_{i=1}^n \qbar_i}} \\
		\stackrel{(c)}{<}&\infty
	\end{align*}
	where $(a)$ holds because projection is nonexpansive; $(b)$ holds because norm-1 is greater than Euclidean norm; and $(c)$ holds by assumption of Theorem \ref{thm.load.balancing} for $|\theta|\left(\frac{p}{p-1}\right)\left(\frac{\tilde{p}}{\tilde{p}-1}\right)\leq\Theta$. Then, the claim holds with
	\begin{align}\label{eq:claim4-thetamax}
		\theta_{\max}=\Theta\left(\frac{\tilde{p}-1}{2}\right),
	\end{align}
	where we used that $\left(\frac{p}{p-1}\right)\tilde{p}= 2$.	This completes the proof.
\Halmos
\endproof

\section{Details of the proofs in Section \ref{sec:generalized.switch}}
	In this section we provide the details of the proofs of the lemmas stated in Section \ref{sec:generalized.switch}, that we use in the proof of Theorem \ref{thm.generalized.switch}.
	
	\subsection{Proof of Lemma \ref{lemma.gen.switch.expo.qu}}\label{app.gen.switch.expo.qu.proof}
	To prove Lemma \ref{lemma.gen.switch.expo.qu} we use the following lemma, which is similar to Lemma \ref{lemma.load.balancing.to.expoqu}. 

\begin{lemma}\label{lemma.gen.switch.to.expo.qu}
	Consider a generalized switch parametrized by $\epsilon$, as described in Theorem \ref{thm.generalized.switch}. Then, for any $\alpha\in\bR$ and for all $k\geq 1$ we have
	\begin{align*}
	\sum_{i=1}^n c_i^\pl u_i^\peps(k)e^{-\frac{\alpha}{c_i^\pl}\qbar^\peps_{\perp i}(k+1)}= \langle\vcl,\vu^\peps(k)\rangle e^{\alpha\langle\vcl,\vq^\peps(k+1)\rangle}
	\end{align*}
\end{lemma}

\proof{Proof of Lemma \ref{lemma.gen.switch.to.expo.qu}.}

First observe that if $\alpha=0$ the lemma trivially holds. Now we prove the lemma for $\alpha\neq 0$. From Equation \eqref{eq.qu} we know that $q_i(k+1)u_i(k)=0$ for all $i\in\{1,\ldots,n\}$. Then, for all $\beta\in\bR$ we have
\begin{align*}
u_i\left(e^{-\beta q_i(k+1)}-1 \right)=0\quad\forall i\in\{1,\ldots,n\},
\end{align*}
and this equation implies
\begin{align*}
c^\pl_i u_i\left(e^{-\beta q_i(k+1)}-1 \right)=0\quad\forall i\in\{1,\ldots,n\}.
\end{align*}

Without loss of generality, we assume $c^\pl_i>0$ for all $i\in\{1,\ldots,n\}$ because otherwise the last equation holds trivially. Let $\alpha\in\bR$ and for each $i\in\{1,\ldots,n\}$ let $\alpha_i\in\bR$ be such that $\alpha=\alpha_ic_i^\pl$ for all $i\in\{1,\ldots,n\}$. Then,
\begin{align*}
c^\pl_i u_i\left(e^{-\alpha_i q_i(k+1)}-1 \right)=0\quad\forall i\in\{1,\ldots,n\}.
\end{align*}

Summing over all $i\in\{1,\ldots,n\}$ we obtain
\begin{align*}
0=& \sum_{i=1}^n c_i^\pl u_i(k)\left(e^{-\alpha_i q_i(k+1)}-1\right) \\
=& \sum_{i=1}^n c_i^\pl u_i(k)\left(e^{-\alpha_i q_{\parallel i}(k+1)-\alpha_i q_{\perp i}(k+1)}-1\right) \\
\stackrel{(a)}{=}& \sum_{i=1}^n c_i^\pl u_i(k)\left(e^{-\alpha_i\langle\vcl,\vq(k+1)\rangle c^\pl_i-\alpha_i q_{\perp i}(k+1)}-1\right) \\
\stackrel{(b)}{=}& \sum_{i=1}^n c^\pl_i u_i(k)\left(e^{-\alpha\langle\vcl,\vq(k+1)\rangle-\frac{\alpha}{c^\pl_i}q_{\perp i}(k+1)} -1\right) \\
\stackrel{(c)}{=}& e^{-\alpha\langle\vcl,\vq(k+1)\rangle}\sum_{i=1}^n c^\pl_i u_i(k)e^{-\frac{\alpha}{c_i^\pl}q_{\perp i}(k+1)}-\langle\vcl,\vu(k)\rangle
\end{align*}
where $(a)$ holds by definition of $\vqpar(k)$; $(b)$ holds by definition of $\alpha$; and $(c)$ holds by expanding the product and reorganizing terms. Therefore, we have
\begin{align*}
\langle\vcl,\vu(k)\rangle= e^{-\alpha\langle\vcl,\vq(k+1)\rangle}\sum_{i=1}^n c^\pl_i u_i(k)e^{-\frac{\alpha}{c_i^\pl}q_{\perp i}(k+1)}.
\end{align*}

Multiplying both sides by $e^{\alpha\langle\vcl,\vq(k+1)\rangle}$ we obtain
\begin{align*}
\langle\vcl,\vu(k)\rangle e^{\alpha\langle\vcl,\vq(k+1)\rangle}= \sum_{i=1}^n c_i^\pl u_i(k) e^{-\frac{\alpha}{c_i^\pl}q_{\perp i}(k+1)},
\end{align*}
which proves the lemma.\Halmos
\endproof

Now we prove Lemma \ref{lemma.gen.switch.expo.qu}.

\proof{Proof of Lemma \ref{lemma.gen.switch.expo.qu}.}

First observe that if $\theta=0$ the lemma holds trivially. Now assume $\theta\neq 0$. Since $\vcl\geq 0$ and $\ubar_i\leq \sbar_i\leq \smax$ for all $i\in\{1,\ldots,n\}$, we have
\begin{align*}
	0\leq \langle\vcl,\vubar\rangle\leq \smax\langle\vcl,\vone \rangle.
\end{align*}

Then, from facts \ref{facts.functions.fraction} and \ref{facts.functions.ineq} stated in Appendix \ref{app.load.balancing.expo.qu} we have
\begin{align}\label{eq.gen.switch.lemma.expo.qu.two}
	\left|e^{-\theteps\langle\vcl,\vubar\rangle} \right|\leq& \left|\theteps\langle\vcl,\vubar\rangle\right|\left(\dfrac{e^{-\theteps \smax\langle\vcl,\vone\rangle}-1}{-\theteps \smax\langle\vcl,\vone\rangle}\right).
\end{align}

Now, by properties of expected value, we have
\begin{align*}
	& \left|\E{\left(e^{\theteps\langle\vcl,\vqbar^\plus\rangle}-1\right)\left(e^{-\theteps\langle\vcl,\vubar\rangle}-1\right)} \right|\\
	\leq& \E{\left|e^{\theteps\langle\vcl,\vqbar^\plus\rangle}-1 \right| \left|e^{-\theteps\langle\vcl,\vubar\rangle}-1\right|} \\
	\stackrel{(a)}{\leq}& |\theteps| \left(\dfrac{e^{-\theteps \smax\langle\vcl,\vone\rangle}-1}{-\theteps \smax\langle\vcl,\vone\rangle}\right) \E{\left|\langle\vcl,\vubar\rangle\left(e^{\theteps\langle\vcl,\vqbar^\plus\rangle}-1 \right) \right|} \\
	\stackrel{(b)}{=}& |\theteps| \left(\dfrac{e^{-\theteps \smax\langle\vcl,\vone\rangle}-1}{-\theteps \smax\langle\vcl,\vone\rangle}\right) \E{\left|\sum_{i=1}^n c^\pl_i\ubar_i\left(e^{-\left(\frac{\theteps}{c^\pl_i}\right)\qbar^\plus_{\perp i}} \right) \right|} \\
	\stackrel{(c)}{\leq}& |\theteps| \left(\dfrac{e^{-\theteps \smax\langle\vcl,\vone\rangle}-1}{-\theteps \smax\langle\vcl,\vone\rangle}\right)\E{\sum_{i=1}^n c_i\ubar_i\left|e^{-\left(\frac{\theteps}{c^\pl_i}\right)\qbar^\plus_{\perp i}} -1\right| } \\
	\stackrel{(d)}{\leq}& |\theteps| \left(\dfrac{e^{-\theteps \smax\langle\vcl,\vone\rangle}-1}{-\theteps \smax\langle\vcl,\vone\rangle}\right) \E{\sum_{i=1}^n \left(c^\pl_i \ubar_i\right)^p}^{\frac{1}{p}} \E{\sum_{i=1}^n \left|e^{-\left(\frac{\theteps}{c^\pl_i}\right)\qbar^\plus_{\perp i}} -1 \right|^{\tfrac{p}{p-1}}}^{\tfrac{p-1}{p}},
\end{align*}
where $p> 1$. Here $(a)$ holds by Equation \eqref{eq.gen.switch.lemma.expo.qu.two}; $(b)$ holds by Lemma \ref{lemma.gen.switch.to.expo.qu} with $\alpha=\theteps$; $(c)$ holds by triangle inequality; and $(d)$ holds by Hölder's inequality.

But
\begin{align*}
	\E{\sum_{i=1}^n \left(c_i^\pl \ubar_i \right)^p}\leq& (c_{\max}\smax)^{p-1} \E{\sum_{i=1}^n c_i^\pl \ubar_i} \leq (c_{\max}\smax)^{p-1} \epsilon
\end{align*}
where $c_{\max}=\max_{i}c^\pl_i$ and the last equality holds by the following reason. By Lemma \ref{lemma.gen.switch.linear.test.func.} we have
\begin{align*}
	\E{\langle\vcl,\vu\rangle}=& \epsilon-\bl+\E{\langle\vcl,\vsbar\rangle}.
\end{align*}
Also, by definition of the capacity region in  \eqref{eq.gen.switch.capacity.region} and because $\vsbar$ depends on the channel state, we have that $\E{\langle\vcl,\vsbar\rangle}\in\cC$. Then, 
\begin{align}\label{eq.gen.switch.b-cs.pos}
	-\bl+\E{\langle\vcl,\vsbar\rangle}\leq 0.
\end{align}

Therefore,
\begin{align*}
	& \left|\E{\left(e^{\theteps\langle\vcl,\vqbar^\plus\rangle}-1\right)\left(e^{-\theteps\langle\vcl,\vubar\rangle}-1\right)} \right|\\
	\leq&|\theta|\epsilon^{1+\frac{1}{p}} (c_{\max}\smax)^{\frac{p-1}{p}} \left(\dfrac{e^{-\theteps \smax\langle\vcl,\vone\rangle}-1}{-\theteps \smax\langle\vcl,\vone\rangle}\right)\E{\sum_{i=1}^n \left|e^{-\left(\frac{\theteps}{c^\pl_i}\right)\qbar^\plus_{\perp i}} -1 \right|^{\frac{p}{p-1}} }^{\frac{p-1}{p}}.
\end{align*}
The rest of the argument is similar to the last steps in the proof of Lemma \ref{lemma.load.balancing.expoqu}. However, in this case we do not need to use existence of the MGF of $\epsilon\sum_{i=1}^n\qbar_i$ because we know $\E{e^{\theteps \|\vq_\perp\|}}$ is bounded for $\theteps\leq \theta^*$ from SSC.\Halmos

\endproof
	
	\subsection{Existence of MGF of $\epsilon \left\|\vqbar\right\|$ in the generalized switch}\label{app.gen.switch.MGF}
	We prove the following lemma.

\begin{lemma}\label{lemma.gen.switch.MGFq}
	Consider a generalized switch parametrized by $\epsilon$ as described in Theorem \ref{thm.generalized.switch}. Then, for each $\epsilon>0$ there exists $\Theta>0$ such that $\E{e^{\theteps\langle\vcl,\vqbar\rangle}}<\infty$ for all $\theta\in[-\Theta,\Theta]$.
\end{lemma}

\proof{Proof of Lemma \ref{lemma.gen.switch.MGFq}.}
First observe that if $\theta=0$ the lemma holds trivially. Therefore, in this proof we assume $\theta\neq 0$. We use Foster-Lyapunov theorem \cite[Proposition 6.13]{hajekrandomprocbook} with Lyapunov function $V(\vq)=e^{\theteps\langle\vcl,\vq\rangle}$. In this proof we use the notation
\begin{align*}
\Evq{\;\cdot\;}\defn\E{\;\cdot\;|\vqbar=\vq} \quad\text{and}\quad \Et{\;\cdot\;}\defn \E{\;\cdot\;|T(k)=t}
\end{align*}

The drift of $V(\vqbar)$ conditioned on $\vqbar=\vq$ is
\begin{align*}
	& \Evq{e^{\theteps\langle\vcl,\vqbar^\plus\rangle}-e^{\theteps\langle\vcl,\vqbar\rangle}} \\
	\stackrel{(a)}{=}& \Evq{e^{\theteps\langle\vcl,\vqbar+\vabar-\vsbar\rangle}- e^{-\theteps\langle\vcl,\vubar\rangle}+1 -e^{\theteps\langle\vcl,\vqbar\rangle}}+o(\epsilon^2) \\
	=& \Evq{e^{\theteps\langle\vcl,\vqbar+\vabar\rangle}e^{-\theteps\langle\vcl,\vsbar\rangle}- e^{-\theteps\langle\vcl,\vubar\rangle}+1 -e^{\theteps\langle\vcl,\vqbar\rangle}}+o(\epsilon^2) \\
	\stackrel{(b)}{=}& \Evq{e^{\theteps\left(\langle\vcl,\vqbar+\vabar\rangle-\Bbar \right)} -e^{-\theteps\langle\vcl,\vubar\rangle}+1-e^{\theteps\langle\vcl,\vqbar\rangle}}+o(\epsilon^2) \\
	&\quad + \Evq{e^{\theteps\langle\vcl,\vqbar+\vabar-\vsbar\rangle}-e^{\theteps\left(\langle\vcl,\vqbar+\vabar\rangle- \Bbar \right)}} \\
	\stackrel{(c)}{=}& \Evq{e^{\theteps\left(\langle\vcl,\vqbar+\vabar\rangle- \Bbar \right)} -e^{-\theteps\langle\vcl,\vubar\rangle}+1-e^{\theteps\langle\vcl,\vqbar\rangle}}+o(\epsilon^2)  \\
	& \quad +\E{e^{\theteps\langle\vcl,\vabar\rangle}} \Evq{e^{\theteps\langle\vcl,\vqbar\rangle}\left(e^{-\theteps\langle\vcl,\vsbar\rangle}-e^{-\theteps \Bbar} \right)}
\end{align*}
where $(a)$ holds expanding the product and rearranging terms in Lemma \ref{lemma.gen.switch.expo.qu}; $(b)$ holds after adding and subtracting $\Evq{e^{\theteps\left(\langle\vcl,\vqbar+\vabar\rangle-\Bbar \right)}}$, and reorganizing terms; $(c)$ holds because the arrival process is independent of the queue lengths and  services processes.

But
\begin{align*}
	\Evq{e^{\theteps\langle\vcl,\vqbar\rangle}\left(e^{-\theteps\langle\vcl,\vsbar\rangle}-e^{-\theteps \Bbar}\right)}=& \Evq{e^{\theteps\left(\langle\vcl,\vqbar\rangle-\Bbar \right)}\left(e^{\theteps\left(\Bbar-\langle\vcl,\vsbar\rangle \right)} -1\right)} \\
	=& \Evq{\Et{e^{-\theteps \Bbar}} \Et{e^{\theteps\langle\vcl,\vqbar\rangle}\left(e^{\theteps\left(\Bbar-\langle\vcl,\vsbar\rangle \right)}-1\right)} } \\
	=& \Evq{\Et{e^{-\theteps \Bbar}} \Et{e^{\theteps\left(\Bbar-\langle\vcl,\vsbar\rangle\right)}-1}}+o(\epsilon^2),
\end{align*}
where the last equality holds by Lemma \ref{lemma.gen.switch.expo.qu} and because the random variable $\Bbar$ is bounded (since it takes finitely many values).

Rearranging terms we obtain
\begin{align}
	\Evq{e^{\theteps\langle\vcl,\vqbar^\plus\rangle} -e^{\theteps\langle\vcl,\vqbar\rangle}} 
	=& 1-\Evq{e^{-\theteps\langle\vcl,\vubar\rangle}}+\E{e^{\theteps\langle\vcl,\vabar\rangle}}\Evq{e^{-\theteps\langle\vcl,\vsbar\rangle} - e^{\theteps \Bbar}}+o(\epsilon^2) \label{eq.gen.switch.MGF.top} \\
	&\quad+ e^{\theteps\langle\vcl,\vq\rangle}\E{e^{\theteps\left(\langle\vcl,\vabar\rangle-\Bbar\right)}-1}. \label{eq.gen.switch.MGF.bottom}
\end{align}

Observe that the right hand side of Equation \eqref{eq.gen.switch.MGF.top} is bounded because $\ubar_i\leq\sbar_i\leq\smax$ and $\abar_i\leq\amax$ with probability 1 for all $i\in\{1,\ldots,n\}$. Also, $\Bbar$ is bounded because it takes a finite number of values.

Then, it suffices to show that for $\delta>0$ there exists $\Theta>0$ such that
\begin{align*}
\E{e^{\theteps\left(\langle\vcl,\vabar\rangle-\Bbar\right)}-1}<-\delta\qquad\forall \theta\in\left[-\Theta,0\right)\cup\left(0,\Theta\right]
\end{align*}

This result can be easily using continuity of MGF at $\theta=0$ and Taylor expansion of $\E{e^{\theteps\left(\langle\vcl,\vabar\rangle-\Bbar\right)}}$ with respect to $\theta$, up to first order around $\theta=0$. We omit the details for brevity.  \Halmos
\endproof
	
	\subsection{Proof of Lemma \ref{lemma.gen.switch.qs}}\label{app.gen.switch.qs}
	\proof{Proof of Lemma \ref{lemma.gen.switch.qs}.}

First observe that if $\theta=0$ the proof holds trivially. Now assume $\theta\neq 0$. 

In this proof we use a geometric vision of MaxWeight algorithm. Before presenting the technical details we present an intuitive overview of the proof. Recall that, given the channel state, MaxWeight algorithm maximizes $\langle\vq(k),\vx\rangle$ over the set of feasible service rate vectors. Then, MaxWeight solves an optimization problem with linear objective function. Equivalently, MaxWeight finds a vector $\vx^*$ which is an optimal solution of
	\begin{align}
		\begin{aligned}\label{gs.eq.MW.LP}
		\max \quad& \langle\vq(k),\vx\rangle \\
		s.t.\quad & \vx\in ConvexHull\left(\cS^\pt\right)
		\end{aligned}
	\end{align}
and sets $\vs(k)$ as one of these optimal solutions. To make the optimization problem linear, we use $ConvexHull\left(\cS^\pt\right)$ as the feasible region instead of $\cS^\pt$. However, this does not change the problem because the objective function is linear and, therefore, an optimal solution of \eqref{gs.eq.MW.LP} is at an extreme point, i.e. at a point in $\cS^\pt$. 

The gradient of the objective function is $\vq(k)$. Then, depending on its direction, the optimal solution(s) $\vx^*$ will belong to a different facet or vertex of $ConvexHull\left(\cS^\pt \right)$. In Figure \ref{fig:example.theta.proof} we present pictorial examples where we show the optimal solution(s) when the vector of queue lengths goes in three different directions. 

\begin{figure}[t]
	\centering
	\begin{subfigure}[t]{0.28\textwidth}
		\includegraphics[width=\linewidth]{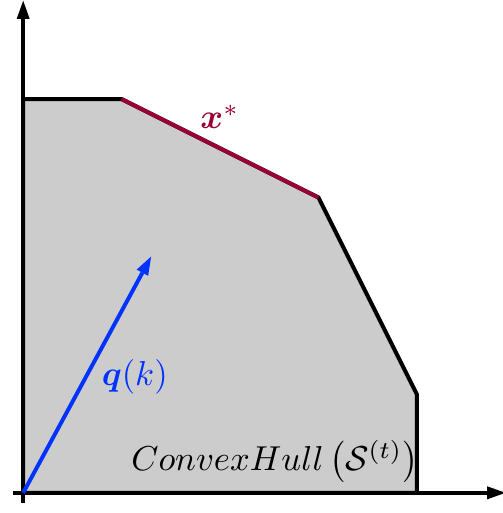}
		\caption{Example 1: Multiple solutions, since $\vq(k)$ is perpendicular to the second facet from left to right.}
	\end{subfigure}
	~
	\begin{subfigure}[t]{0.28\textwidth}
		\includegraphics[width=\linewidth]{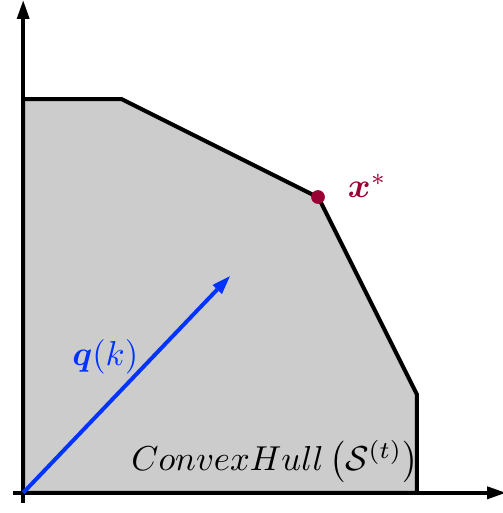}
		\caption{Example 2: Unique solution.}
	\end{subfigure}
	~
	\begin{subfigure}[t]{0.28\textwidth}
		\includegraphics[width=\linewidth]{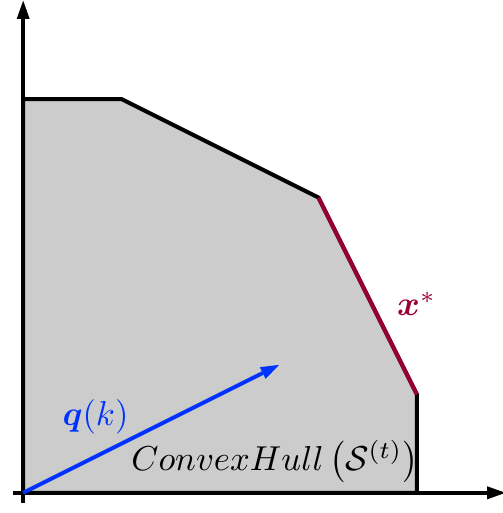}
		\caption{Example 3: Another example of multiple solutions.}
	\end{subfigure}
	{\caption{Example of optimal solutions depending on the queue lengths vector.}\label{fig:example.theta.proof}}
\end{figure}

Recall that $\vqpar(k)$ goes in the same direction as $\vcl$. Also, if $\epsilon$ is small we expect that $\vq(k)\approx \vqpar(k)$ by SSC. Then, if $\epsilon$ is small we expect that any optimal solution $\vx^*$ to the linear program \eqref{gs.eq.MW.LP} satisfies $\langle\vcl,\vx^*\rangle=\btl$ with high probability. 

Now we present the technical details. 
We start with a definition. Let $t\in\cT$ and suppose that the channel state is $\Tbar=t$. Then, let $\nu^\pt\in\left(0,\frac{\pi}{2} \right]$ be an angle such that $\langle\vcl,\vsbar\rangle=\btl$ if $\dfrac{\left\|\vqparbar \right\|}{\left\|\vqbar \right\|}\geq \cos\left(\nu^\pt\right)$. Let $\nuq$ be the angle between $\vqparbar$ and $\vqbar$ and define $\numin\defn \min_{t\in \cT}\nu^\pt$. Therefore, since $\vsbar$ is scheduled using MaxWeight algorithm, if channel state is $t$ we have
\begin{align}\label{eq.gs.def.angle}
	\btl\neq \langle\vcl,\vsbar\rangle\quad\text{implies}\quad \nu^\pt<\nuq.
\end{align}

In this proof we use the notation $\Et{\,\cdot\,}=\E{\,\cdot\;|\Tbar=t}$. By definition of conditional expectation we have
\begin{align*}
	\E{\left(e^{\theteps\langle\vcl,\vqbar^\peps\rangle}-1\right) \left(e^{\theteps\left(\Bbar-\langle\vcl,\vsbar^\peps\rangle \right)}-1\right)} =& \sum_{t\in\cT}\psi_t \Et{\left(e^{\theteps\langle\vcl,\vqbar^\peps\rangle}-1\right) \left(e^{\theteps\left(\btl-\langle\vcl,\vsbar^\peps\rangle \right)}-1\right)},
\end{align*}
where 
\begin{align*}
	& \Et{\left(e^{\theteps\langle\vcl,\vqbar^\peps\rangle}-1\right) \left(e^{\theteps\left(\btl-\langle\vcl,\vsbar^\peps\rangle \right)}-1\right)} \\
	\stackrel{(a)}{=}& \Et{\left(e^{\theteps \left\|\vqparbar \right\|}-1\right)\left(e^{\theteps\left(\btl-\langle\vcl,\vsbar\rangle \right)}-1\right)\ind{\btl\neq \langle\vcl,\vsbar\rangle}} \\
	\stackrel{(b)}{\leq}& \Et{\left(e^{\theteps \left\|\vqparbar \right\|}-1\right)\left(e^{\theteps\left(\btl-\langle\vcl,\vsbar\rangle \right)}-1\right)\ind{\nuq>\nu^\pt}} \\
	=& \Et{\left(e^{\theteps\left\|\vqperpbar \right \|\cot(\nuq)}-1\right)\left(e^{\theteps\left(\btl-\langle\vcl,\vsbar\rangle \right)}-1\right)\ind{\nuq>\nu^\pt}} \\
	\stackrel{(c)}{\leq}& \Et{\left(e^{\theteps\left\|\vqperpbar \right \|\cot(\nu^\pt)}-1\right)\left(e^{\theteps\left(\btl-\langle\vcl,\vsbar\rangle \right)}-1\right)\ind{\nuq>\nu^\pt}} \\
	\stackrel{(d)}{=}& \Et{\left(e^{\theteps\left\|\vqperpbar \right \|\cot(\nu^\pt)}-1\right)\left(e^{\theteps\left(\btl-\langle\vcl,\vsbar\rangle \right)}-1\right)} \\
	\stackrel{(e)}{\leq}& \Et{\left(e^{\theteps \left\|\vqperpbar \right\|\cot(\nu^\pt)} -1\right)^p}^{\frac{1}{p}} \Et{\left(e^{\theteps\left(\btl-\langle\vcl,\vsbar\rangle \right)} -1\right)^{\frac{p}{p-1}} }^\frac{p-1}{p},
\end{align*}
where $p>1$. Here $(a)$ holds by definition of indicator function and because $\vqparbar=\langle\vcl,\vqbar\rangle\vcl$ by definition of projection; $(b)$ holds by \eqref{eq.gs.def.angle}; $(c)$ and $(d)$ holds because $\cot(\nu)$ is decreasing for $\nu\in\left(0,\frac{\pi}{2}\right]$; $(d)$ holds by \eqref{eq.gs.def.angle} and by definition of indicator function; and $(e)$ holds by Hölder's inequality.

Using an argument similar to the one at the end of Lemma \ref{lemma.load.balancing.expoqu}, it can be proved that 
\begin{align*}
	0\leq \Et{\left(e^{\theteps \left\|\vqperpbar \right\|\cot(\nu^\pt)}-1\right)^p}^{\frac{1}{p}}
\end{align*}
converges to a constant as $\epsilon\dto 0$. On the other hand,
\begin{align*}
	& \Et{\left(e^{\theteps\left(\btl-\langle\vcl,\vsbar\rangle \right)} -1\right)^{\frac{p}{p-1}}}\\ 
	=& 	\Et{\left(e^{\theteps\left(\btl-\langle\vcl,\vsbar\rangle \right)} -1\right)^{\frac{p}{p-1}} \ind{\btl\neq\langle\vcl,\vsbar\rangle}} \\
	=& \E{\left(\dfrac{e^{\theteps\left(\btl-\langle\vcl,\vsbar\rangle \right)} -1}{\theteps\left(\btl-\langle\vcl,\vsbar\rangle \right)}\right)^{\frac{p}{p-1}} \left(\theteps\left(\btl-\langle\vcl,\vsbar\rangle \right)\right)^{\frac{p}{p-1}} \ind{\btl\neq\langle\vcl,\vsbar\rangle}} \\
	\leq& \left(\dfrac{e^{\theteps\left(\Bbar_{\max}-\langle\vcl, \smax \vone\rangle \right)} -1}{\theteps\left(\Bbar_{\max}-\langle\vcl,\smax \vone\rangle \right)}\right)^{\frac{p}{p-1}} \left(\theteps\left(\Bbar_{\max}-\langle\vcl,\smax \vone\rangle \right)\right)^{\frac{p}{p-1}}\Prob{\btl\neq \langle\vcl,\vsbar\rangle}
\end{align*}
where $\Bbar_{\max}=\max_{t\in\cT}\btl$. \cite{atilla} prove that $\Prob{\btl\neq \langle\vcl,\vsbar\rangle}=K\epsilon$ for a finite constant $K$, and their proof also holds here. Therefore, 
\begin{align*}
	\E{\left(e^{\theteps\left(\btl-\langle\vcl,\vsbar\rangle \right)} -1\right)^{\frac{p}{p-1}}}\text{ is }O\left(\epsilon^{1+ \frac{p}{p-1}}\right)
\end{align*}

This completes the proof.\Halmos
\endproof

\section{Proof of the claims in Sections \ref{sec:load.balancing.theorem} and \ref{sec:generalized.switch.theorem}}\label{app.claims}
	In this appendix we show the proof of all the claims that we did in the proofs of our Theorems.

\subsection{Proof of Claim \ref{claim.load.balancing}}\label{app.claim.load.balancing}

\proof{Proof of Claim \ref{claim.load.balancing}.}
	We have
	\begin{align*}
		0\leq \dfrac{(\theteps)^2}{2}\E{\left(\sum_{i=1}^n \ubar_i\right)^2}\stackrel{(a)}{\leq}& \epsilon^2 \left(\dfrac{n\smax\theta^2}{2}\right)\E{\sum_{i=1}^n \ubar_i} \\
		\stackrel{(b)}{=}& \epsilon^3 \left(\dfrac{n\smax\theta^2}{2}\right)
	\end{align*}
	where $(a)$ holds because, by definition of unused service, we have $\ubar_i\leq \sbar_i\leq \smax$ and all terms are nonnegative; and $(b)$ holds by Lemma \ref{lemma.load.balancing.Eu}.
	
	Therefore,
	\begin{align*}
		\dfrac{(\theteps)^2}{2}\E{\left(\sum_{i=1}^n \ubar_i\right)^2}\; \text{is }O(\epsilon^3). 
	\end{align*}
	\Halmos
\endproof

\subsection{Proof of Claim \ref{claim.gen.switch.u}}\label{app.claim.gen.switch.u}

Now we prove Claim \ref{claim.gen.switch.u}.
\proof{Proof of Claim \ref{claim.gen.switch.u}.}
We have
\begin{align*}
	0\leq \dfrac{(\theteps)^2}{2}\E{\left(\langle\vcl,\vubar\rangle \right)^2}\stackrel{(a)}{\leq}& \epsilon^2 \left(\dfrac{\langle \vcl,\smax\vone\rangle \theta^2}{2}\right)\E{\langle\vcl,\vubar\rangle} \\
	\stackrel{(b)}{\leq}&\epsilon^3 \left(\dfrac{\langle \vcl,\smax\vone\rangle \theta^2}{2}\right)
\end{align*}
where $(a)$ holds because $\ubar_i\leq \sbar_i\leq \smax$ and $\vcl\geq 0$;  and $(b)$ holds by Lemma \ref{lemma.gen.switch.linear.test.func.}, because $\E{\langle\vcl,\vsbar\rangle-\Bbar}\leq 0$.

Therefore, 
	\begin{align*}
		\dfrac{(\theteps)^2}{2}\E{\left(\langle\vcl,\vubar\rangle \right)^2}\; \text{is }O(\epsilon^3). 
	\end{align*}
\Halmos
\endproof

\subsection{Proof of Claim \ref{claim.gen.switch.num2}}\label{app.claim.gen.switch.num2}

Now we prove Claim \ref{claim.gen.switch.num2}. 
\proof{Proof of Claim \ref{claim.gen.switch.num2}.}
For the first expression, from Lemma \ref{lemma.taylor} we have
\begin{align*}
	\E{f_{\epsilon,\left(\langle\vcl,\vabar\rangle-\Bbar\right)}(\theta)}=& 1+\theteps\E{\langle\vcl,\vabar\rangle-\Bbar}+ \dfrac{(\theteps)^2}{2}\E{\left(\langle\vcl,\vabar\rangle-\Bbar\right)^2}+O(\epsilon^3) \\
	=& 1+\theteps^2+\dfrac{(\theteps)^2}{2}\E{\left(\langle\vcl,\vabar\rangle-\Bbar\right)^2}+O(\epsilon^3),
\end{align*}
where the last equality holds by Equation \eqref{eq.gen.switch.Eca}. Also,
\begin{align*}
	0\leq \dfrac{(\theteps)^2}{2}\E{\left(\langle\vcl,\vabar\rangle-\Bbar\right)^2}\stackrel{(a)}{\leq}& \epsilon^2 \left(\dfrac{(\langle\vcl,\amax\vone\rangle+B_{\max})\theta^2}{2}\right)\E{\langle\vcl,\vabar\rangle-\Bbar} \\
	\stackrel{(b)}{=}& \epsilon^3 \left(\dfrac{(\langle\vcl,\amax\vone\rangle+B_{\max})\theta^2}{2}\right)
\end{align*}
where $(a)$ holds because $\abar_i\leq \amax$ with probability 1 for all $i\in\{1,\ldots,n\}$, $\vcl\geq 0$, $\Bbar$ is bounded by a constant that we denote $B_{\max}$ and because all quantities are nonnegative; and $(b)$ holds by Equation \eqref{eq.gen.switch.Eca}. Then, 
\begin{align*}
	\dfrac{(\theteps)^2}{2}\E{\left(\langle\vcl,\vabar\rangle-\Bbar\right)^2}\;\text{is }O(\epsilon^3).
\end{align*}
Therefore,
\begin{align*}
	\E{f_{\epsilon,\left(\langle\vcl,\vabar\rangle-\Bbar\right)}(\theta)}= 1+\theteps^2+O(\epsilon^3).
\end{align*}
This proves the first equation of the claim.

For the second expression, using Lemma \ref{lemma.taylor} we obtain
\begin{align*}
	\E{f_{\epsilon,\left(\Bbar-\langle\vcl,\vsbar\rangle\right)}(\theta)}-1=& \theteps\E{\Bbar-\langle\vcl,\vsbar\rangle}+ \dfrac{(\theteps)^2}{2}\E{\left(\Bbar-\langle\vcl,\vsbar\rangle\right)^2}+O(\epsilon^3).
\end{align*}

But
\begin{align*}
	0\leq \dfrac{(\theteps)^2}{2}\E{\left(\Bbar-\langle\vcl,\vsbar\rangle\right)^2} \stackrel{(a)}{\leq}& \epsilon^2\left(\dfrac{\left(B_{\max}+\langle\vcl,\smax\vone\rangle \right)\theta^2}{2}\right)\E{\Bbar-\langle\vcl,\vsbar\rangle} \\
	\stackrel{(b)}{\leq}& \epsilon^3\left(\dfrac{\left(B_{\max}+\langle\vcl,\smax\vone\rangle \right)\theta^2}{2}\right)
\end{align*}
where $(a)$ holds because $\sbar_i\leq \smax$ with probability 1 for all $i\in\{1,\ldots,n\}$, $\vcl\geq 0$, $\Bbar\leq B_{\max}$ and all quantities are nonnegative (see Equation \eqref{eq.gen.switch.b-cs.pos} to see why $\E{\Bbar-\langle\vcl,\vsbar\rangle}\geq 0$); and $(b)$ holds by Lemma \ref{lemma.gen.switch.linear.test.func.} and because $\E{\langle\vcl,\vubar\rangle}\geq 0$ since $\vubar\geq 0$ and $\vcl\geq 0$. 

Then,
\begin{align*}
	\dfrac{(\theteps)^2}{2}\E{\left(\Bbar-\langle\vcl,\vsbar\rangle\right)^2}\; \text{is }O(\epsilon^3).
\end{align*}

Therefore,
\begin{align*}
	\E{f_{\epsilon,\left(\Bbar-\langle\vcl,\vsbar\rangle\right)}(\theta)}-1= \theteps\E{\Bbar-\langle\vcl,\vsbar\rangle}+O(\epsilon^3).
\end{align*}
\Halmos
\endproof

\end{APPENDIX}

\ACKNOWLEDGMENT{We thank Professor Jim Dai for his meaningful feedback on the proof of Lemma 4. \\ 
	We acknowledge the support from iDDA at the Chinese University of Hong King, Shenzhen during the Summer 2018. This research was partially supported by the NSF grant NSF-CCF: 1850439. Daniela Hurtado-Lange has partial funding from ANID/DOCTORADO BECAS CHILE/2018 - 72190413}

\bibliographystyle{informs2014}
\bibliography{biblio-ok}

\begin{thebibliography}{68}
\providecommand{\natexlab}[1]{#1}
\providecommand{\url}[1]{\texttt{#1}}
\providecommand{\urlprefix}{URL }

\bibitem[{Bertsimas et~al.(2001)Bertsimas, Gamarnik, \protect\BIBand{}
  Tsitsiklis}]{bertsimas_momentbound}
Bertsimas D, Gamarnik D, Tsitsiklis JN (2001) Performance of multiclass
  markovian queueing networks via piecewise linear {L}yapunov functions.
  \emph{Ann. Appl. Probab.} 11(4):1384--1428.

\bibitem[{Braverman(2018)}]{braverman2018jsq}
Braverman A (2018) Steady-state analysis of the join the shortest queue model
  in the {H}alfin-{W}hitt regime. \emph{arXiv preprint arXiv:1801.05121} .

\bibitem[{Braverman \protect\BIBand{} Dai(2017)}]{braverman2017stein2}
Braverman A, Dai J (2017) Stein's method for steady-state diffusion
  approximations of {M/Ph/n+ M} systems. \emph{The Annals of Applied
  Probability} 27(1):550--581.

\bibitem[{Braverman et~al.(2017{\natexlab{a}})Braverman, Dai, \protect\BIBand{}
  Feng}]{braverman2017stein}
Braverman A, Dai J, Feng J (2017{\natexlab{a}}) Stein's method for steady-state
  diffusion approximations: {A}n introduction through the {E}rlang-{A} and
  {E}rlang-{C} models. \emph{Stochastic Systems} 6(2):301--366.

\bibitem[{Braverman et~al.(2017{\natexlab{b}})Braverman, Dai, \protect\BIBand{}
  Miyazawa}]{braverman_BAR}
Braverman A, Dai J, Miyazawa M (2017{\natexlab{b}}) Heavy traffic approximation
  for the stationary distribution of a {G}eneralized {J}ackson {N}etwork: The
  {BAR} approach. \emph{Stochastic Systems} 7(1):143--196.

\bibitem[{Braverman et~al.(2018)Braverman, Gurvich, \protect\BIBand{}
  Huang}]{braverman_truncation}
Braverman A, Gurvich I, Huang J (2018) On the {T}aylor expansion of value
  functions. \emph{arXiv preprint arXiv:1804.05011} .

\bibitem[{Chen \protect\BIBand{} Ye(2012)}]{chen2012po2_journal}
Chen H, Ye H (2012) Asymptotic optimality of balanced routing. \emph{Operations
  Research} 60(1):163--179.

\bibitem[{Dai(2018)}]{jim_achievement_lecture}
Dai J (2018) Steady-state approximations: Achievement lecture. \emph{Abstracts
  of the 2018 ACM International Conference on Measurement and Modeling of
  Computer Systems}, 1--1 (ACM).

\bibitem[{Dai \protect\BIBand{} Lin(2008)}]{dai2008max_pressure}
Dai J, Lin W (2008) Asymptotic optimality of maximum pressure policies in
  stochastic processing networks. \emph{The Annals of Applied Probability}
  18(6):2239--2299.

\bibitem[{Ephremides et~al.(1980)Ephremides, Varaiya, \protect\BIBand{}
  Walrand}]{ephremides1980simple}
Ephremides A, Varaiya P, Walrand J (1980) A simple dynamic routing problem.
  \emph{IEEE Transactions on Automatic Control} 25(4):690--693.

\bibitem[{Eryilmaz \protect\BIBand{} Srikant(2012)}]{atilla}
Eryilmaz A, Srikant R (2012) Asymptotically tight steady-state queue length
  bounds implied by drift conditions. \emph{Queueing Systems} 72(3-4):311--359,
  ISSN 0257-0130.

\bibitem[{Eschenfeldt \protect\BIBand{} Gamarnik(2018)}]{Gamarnik_JSQ}
Eschenfeldt P, Gamarnik D (2018) Join the {S}hortest {Q}ueue with many servers.
  {T}he heavy-traffic asymptotics. \emph{Mathematics of Operations Research} .

\bibitem[{Foschini \protect\BIBand{} Salz(1978)}]{JSQ_HT_optimality}
Foschini G, Salz J (1978) A basic dynamic routing problem and diffusion.
  \emph{IEEE Transactions on Communications} 26(3):320--327.

\bibitem[{Gamarnik \protect\BIBand{} Zeevi(2006)}]{gamarnik2006validity}
Gamarnik D, Zeevi A (2006) Validity of heavy traffic steady-state
  approximations in {G}eneralized {J}ackson {N}etworks. \emph{The Annals of
  Applied Probability} 56--90.

\bibitem[{Gaunt \protect\BIBand{} Walton(2020)}]{Walton_SteinHT}
Gaunt R, Walton N (2020) Stein's method for the single server queue in heavy
  traffic. \emph{Statistics \& Probability Letters} 156:108566.

\bibitem[{Gupta \protect\BIBand{} Shroff(2010)}]{gupta2010delay}
Gupta G, Shroff N (2010) Delay analysis for wireless networks with single hop
  traffic and general interference constraints. \emph{IEEE/ACM Transactions on
  Networking (TON)} 18(2):393--405.

\bibitem[{Gurvich(2014)}]{gurvich2014diffusion}
Gurvich I (2014) Diffusion models and steady-state approximations for
  exponentially ergodic {M}arkovian queues. \emph{The Annals of Applied
  Probability} 24(6):2527--2559.

\bibitem[{Gurvich et~al.(2013)Gurvich, Huang, \protect\BIBand{}
  Mandelbaum}]{gurvich2013excursion}
Gurvich I, Huang J, Mandelbaum A (2013) Excursion-based universal
  approximations for the {E}rlang-{A} queue in steady-state. \emph{Mathematics
  of Operations Research} 39(2):325--373.

\bibitem[{Gut(2012)}]{gut2012probability}
Gut A (2012) \emph{Probability: {A} Graduate Course}, volume~75 (Springer
  Science \& Business Media).

\bibitem[{Hajek(1982)}]{hajek_drift}
Hajek B (1982) Hitting-time and occupation-time bounds implied by drift
  analysis with applications. \emph{Advances in Applied Probability} 502--525.

\bibitem[{Hajek(2015)}]{hajekrandomprocbook}
Hajek B (2015) \emph{Random Processes for Engineers} (Cambridge university
  press).

\bibitem[{Harrison(1988)}]{harrison1988brownian}
Harrison J (1988) Brownian models of queueing networks with heterogeneous
  customer populations. \emph{Stochastic Differential Systems, Stochastic
  Control Theory and Applications}, 147--186 (Springer).

\bibitem[{Harrison(1998)}]{har_state_space}
Harrison J (1998) Heavy traffic analysis of a system with parallel servers:
  Asymptotic optimality of discrete review policies. \emph{Ann. App. Probab.}
  822--848.

\bibitem[{Harrison \protect\BIBand{} L{\'o}pez(1999)}]{harlop_state_space}
Harrison J, L{\'o}pez M (1999) Heavy traffic resource pooling in
  parallel-server systems. \emph{Queueing Systems} 339--368.

\bibitem[{Harrison \protect\BIBand{} Patel(1992)}]{harrison1992performance}
Harrison P, Patel N (1992) \emph{Performance Modeling of Communication Networks
  and Computer Architectures} (Addison-Wesley Longman Publishing Co., Inc.).

\bibitem[{Huang \protect\BIBand{}
  Gurvich(2018)}]{huang_gurvich2018universality}
Huang J, Gurvich I (2018) Beyond heavy-traffic regimes: {U}niversal bounds and
  controls for the single-server queue. \emph{Operations Research}
  66(4):1168--1188.

\bibitem[{Hurtado-Lange \protect\BIBand{}
  Maguluri(2019)}]{Hurtado_gen-switch_temp}
Hurtado-Lange D, Maguluri ST (2019) Heavy-traffic analysis of queueing systems
  with no complete resource pooling. Technical Report
  \url{https://arxiv.org/pdf/1904.10096.pdf}.

\bibitem[{Kang et~al.(2009)Kang, Kelly, Lee, \protect\BIBand{}
  Williams}]{kang2009state}
Kang W, Kelly F, Lee N, Williams R (2009) State space collapse and diffusion
  approximation for a network operating under a fair bandwidth sharing policy.
  \emph{The Annals of Applied Probability} 1719--1780.

\bibitem[{Kingman(1961)}]{kingman1961_charfunction}
Kingman J (1961) The single server queue in heavy traffic. \emph{Mathematical
  Proceedings of the Cambridge Philosophical Society}, volume~57, 902--904
  (Cambridge University Press).

\bibitem[{Kingman(1962{\natexlab{a}})}]{kingman1962_brownian}
Kingman J (1962{\natexlab{a}}) On queues in heavy traffic. \emph{Journal of the
  Royal Statistical Society. Series B (Methodological)} 383--392.

\bibitem[{Kingman(1962{\natexlab{b}})}]{kingman}
Kingman J (1962{\natexlab{b}}) Some inequalities for the queue {GI/G/1}.
  \emph{Biometrika} 315--324.

\bibitem[{K{\"o}llerstr{\"o}m(1974)}]{kollerstrom1974heavy}
K{\"o}llerstr{\"o}m J (1974) Heavy traffic theory for queues with several
  servers. i. \emph{Journal of Applied Probability} 11(3):544--552.

\bibitem[{Lehoczky(1996)}]{lehoczky1996real}
Lehoczky J (1996) Real-time queueing theory. \emph{Real-Time Systems Symposium}
  186.

\bibitem[{Lehoczky(1997)}]{lehoczky1997using}
Lehoczky J (1997) Using real-time queueing theory to control lateness in
  real-time systems. \emph{ACM SIGMETRICS Performance Evaluation Review}
  25(1):158--168.

\bibitem[{Li et~al.(2018)Li, Kong, \protect\BIBand{}
  Wang}]{li2018loadbalancing}
Li B, Kong X, Wang L (2018) Optimal load-balancing for high-density wireless
  networks with flow-level dynamics. \emph{Proceedings of the Eighteenth ACM
  International Symposium on Mobile Ad Hoc Networking and Computing} 316--317.

\bibitem[{Lindley(1952)}]{lindley_equation}
Lindley D (1952) The theory of queues with a single server. \emph{Mathematical
  Proceedings of the Cambridge Philosophical Society}, volume~48, 277--289
  (Cambridge University Press).

\bibitem[{Liu \protect\BIBand{} Ying(2019)}]{liu2018simple}
Liu X, Ying L (2019) A simple steady-state analysis of load balancing
  algorithms in the sub-{H}alfin-{W}hitt regime. \emph{ACM SIGMETRICS
  Performance Evaluation Review} 46(2):15--17.

\bibitem[{Lu et~al.(2011)Lu, Xie, Kliot, Geller, Larus, \protect\BIBand{}
  Greenberg}]{lu_JIQ}
Lu Y, Xie Q, Kliot G, Geller A, Larus J, Greenberg A (2011)
  Join-{I}dle-{Q}ueue: A novel load balancing algorithm for dynamically
  scalable web services. \emph{Performance Evaluation} 68(11):1056--1071.

\bibitem[{Lukacs(1970)}]{lukacs1970characteristic}
Lukacs E (1970) \emph{Characteristic Functions} (Griffin).

\bibitem[{Maguluri et~al.(2018)Maguluri, Burle, \protect\BIBand{}
  Srikant}]{QUESTA_switch}
Maguluri ST, Burle S, Srikant R (2018) Optimal heavy-traffic queue length
  scaling in an incompletely saturated switch. \emph{Queueing Systems}
  88(3-4):279--309.

\bibitem[{Maguluri \protect\BIBand{} Srikant(2016)}]{MagSri_SSY16_Switch}
Maguluri ST, Srikant R (2016) Heavy traffic queue length behavior in a switch
  under the {M}ax{W}eight algorithm. \emph{Stoch. Syst.} 6(1):211--250,
  \urlprefix\url{http://dx.doi.org/10.1214/15-SSY193}.

\bibitem[{Maguluri et~al.(2014)Maguluri, Srikant, \protect\BIBand{}
  Ying}]{magsriyin_itc12_journal}
Maguluri ST, Srikant R, Ying L (2014) Heavy traffic optimal resource allocation
  algorithms for cloud computing clusters. \emph{Performance Evaluation}
  81:20--39.

\bibitem[{Marshall(1968)}]{marshall_ineq_drift}
Marshall K (1968) Some inequalities in queuing. \emph{Operations research}
  16(3):651--668.

\bibitem[{Massouli{\'e} \protect\BIBand{}
  Roberts(2000)}]{massoulierobertsbandwidth}
Massouli{\'e} L, Roberts J (2000) Bandwidth sharing and admission control for
  elastic traffic. \emph{Telecommunication systems} 15(1-2):185--201.

\bibitem[{McKeown et~al.(1996)McKeown, Anantharam, \protect\BIBand{}
  Walrand}]{mckeown96walrand}
McKeown N, Anantharam V, Walrand J (1996) Achieving 100\% throughput in an
  input queued switch. \emph{Proceedings of IEEE INFOCOM}, 296--302.

\bibitem[{Meyn(2008)}]{Mey_08}
Meyn S (2008) Stability and asymptotic optimality of generalized {M}ax{W}eight
  policies. \emph{SIAM J. Control and Optimization} To appear.

\bibitem[{Mitzenmacher(1996)}]{mitzenmacher_po2}
Mitzenmacher M (1996) Load balancing and density dependent jump {M}arkov
  processes. \emph{focs}, 213 (IEEE).

\bibitem[{Mitzenmacher(2001)}]{mitzenmacher_po2_2}
Mitzenmacher M (2001) The power of two choices in randomized load balancing.
  \emph{IEEE Transactions on Parallel and Distributed Systems}
  12(10):1094--1104.

\bibitem[{Mood(1950)}]{mood}
Mood A (1950) \emph{Introduction to the Theory of Statistics.} (McGraw-hill).

\bibitem[{Skorokhod(1961)}]{Skorohod_map}
Skorokhod A (1961) Stochastic equations for diffusion processes in a bounded
  region. \emph{Theory of Probability \& Its Applications} 6(3):264--274.

\bibitem[{Srikant \protect\BIBand{} Ying(2014)}]{srikantleibook}
Srikant R, Ying L (2014) \emph{Communication Networks: An Optimization, Control
  and Stochastic Networks Perspective} (Cambridge University Press), ISBN
  9781107036055.

\bibitem[{Stolyar(2004)}]{stolyar2004maxweight}
Stolyar A (2004) Max{W}eight scheduling in a generalized switch: State space
  collapse and workload minimization in heavy traffic. \emph{Annals of Applied
  Probability} 1--53.

\bibitem[{Stolyar(2017)}]{stolyar_JIQ}
Stolyar A (2017) Pull-based load distribution among heterogeneous parallel
  servers: The case of multiple routers. \emph{Queueing Systems}
  85(1-2):31--65.

\bibitem[{Tassiulas \protect\BIBand{} Ephremides(1992)}]{TasEph_92}
Tassiulas L, Ephremides A (1992) Stability properties of constrained queueing
  systems and scheduling policies for maximum throughput in multihop radio
  networks. \emph{IEEE Transactions on Automatic Control} 37(12):1936--1948.

\bibitem[{van~der Boor et~al.(2018)van~der Boor, Borst, van Leeuwaarden,
  \protect\BIBand{} Mukherjee}]{load_balancing_survey}
van~der Boor M, Borst S, van Leeuwaarden J, Mukherjee D (2018) Scalable load
  balancing in networked systems: A survey of recent advances. \emph{arXiv
  preprint arXiv:1806.05444} .

\bibitem[{Vlasiou et~al.(2014)Vlasiou, Zhang, \protect\BIBand{}
  Zwart}]{zwart_bandwidth_diffusion}
Vlasiou M, Zhang J, Zwart B (2014) Insensitivity of proportional fairness in
  critically loaded bandwidth sharing networks. \emph{arXiv preprint
  arXiv:1411.4841} .

\bibitem[{Vvedenskaya et~al.(1996)Vvedenskaya, Dobrushin, \protect\BIBand{}
  Karpelevich}]{dobrushin_po2}
Vvedenskaya N, Dobrushin R, Karpelevich F (1996) Queueing system with selection
  of the shortest of two queues: An asymptotic approach. \emph{Problems of
  Information Transmission} 32(1):15--27.

\bibitem[{Wang et~al.(2017)Wang, Maguluri, \protect\BIBand{}
  Javidi}]{javidi_optical_ht}
Wang CH, Maguluri ST, Javidi T (2017) Heavy traffic queue length behavior in
  switches with reconfiguration delay. \emph{INFOCOM 2017-IEEE Conference on
  Computer Communications, IEEE}, 1--9 (IEEE).

\bibitem[{Wang et~al.(2018)Wang, Maguluri, Srikant, \protect\BIBand{}
  Ying}]{Weina_bandwidth_journal}
Wang W, Maguluri ST, Srikant R, Ying L (2018) Heavy-traffic delay insensitivity
  in connection-level models of data transfer with proportionally fair
  bandwidth sharing. \emph{SIGMETRICS Perform. Eval. Rev.} 45(3):232--245, ISSN
  0163-5999, \urlprefix\url{http://dx.doi.org/10.1145/3199524.3199565}.

\bibitem[{Weber(1978)}]{weber1978optimal}
Weber R (1978) On the optimal assignment of customers to parallel servers.
  \emph{Journal of Applied Probability} 15(2):406--413.

\bibitem[{Williams(1998)}]{Williams_state_space}
Williams R (1998) Diffusion approximations for open multiclass queueing
  networks: Sufficient conditions involving state space collapse.
  \emph{Queueing Systems Theory and Applications} 27 -- 88.

\bibitem[{Williams(2000)}]{Williams_CRP}
Williams R (2000) On dynamic scheduling of a parallel server system with
  complete resource pooling. \emph{Fields Institute Communications}
  28(49-71):5--1.

\bibitem[{Winston(1977)}]{winston_JSQ_1977}
Winston W (1977) Optimality of the shortest line discipline. \emph{Journal of
  Applied Probability} 14(1):181–189,
  \urlprefix\url{http://dx.doi.org/10.1017/S0021900200104772}.

\bibitem[{Ye \protect\BIBand{} Yao(2012)}]{yeyaobandwidth2012}
Ye HQ, Yao D (2012) A stochastic network under proportional fair resource
  control—diffusion limit with multiple bottlenecks. \emph{Operations
  Research} 60(3):716--738.

\bibitem[{Ying(2016)}]{Lei_Stein_SIGM16}
Ying L (2016) On the approximation error of mean-field models.
  \emph{Proceedings of the 2016 ACM SIGMETRICS International Conference on
  Measurement and Modeling of Computer Science}, 285--297, SIGMETRICS '16 (New
  York, NY, USA: ACM), ISBN 978-1-4503-4266-7,
  \urlprefix\url{http://dx.doi.org/10.1145/2896377.2901463}.

\bibitem[{Ying(2017)}]{Lei_steinHT_SIGM17}
Ying L (2017) Stein's method for mean field approximations in light and heavy
  traffic regimes. \emph{Proc. ACM Meas. Anal. Comput. Syst.} 1(1):12:1--12:27,
  ISSN 2476-1249, \urlprefix\url{http://dx.doi.org/10.1145/3084449}.

\bibitem[{Ying et~al.(2017)Ying, Srikant, \protect\BIBand{}
  Kang}]{ying_power_more_than2}
Ying L, Srikant R, Kang X (2017) The power of slightly more than one sample in
  randomized load balancing. \emph{Mathematics of Operations Research}
  42(3):692--722.

\bibitem[{Zhou et~al.(2018)Zhou, Tan, \protect\BIBand{}
  Shroff}]{zhou2018flexible}
Zhou X, Tan J, Shroff N (2018) Flexible load balancing with multi-dimensional
  state-space collapse: Throughput and heavy-traffic delay optimality.
  \emph{Performance Evaluation} 127:176--193.

\end{thebibliography}

\end{document}